\documentclass[12pt]{article}
\usepackage{amsmath,amssymb,amsthm, amscd, pb-diagram}

\newtheorem{Thm}[subsection]{Theorem}
\newtheorem{Prop}[subsection]{Proposition}
\newtheorem{Lem}[subsection]{Lemma}
\newtheorem{Exam}[subsection]{Example}

\newtheorem{Cor}[subsection]{Corollary}
\newtheorem{Rem}[subsection]{Remark}

\newtheorem{Paragraph}[subsection]{}

\begin{document}

\begin{center}
{\bf Towards a characterization of toric hyperk\"{a}hler varieties among symplectic singularities}
\end{center}
\vspace{0.4cm}

\begin{center}
{\bf Yoshinori Namikawa}      
\end{center}
\vspace{0.2cm}

\begin{abstract}
Let $(X, \omega)$ be a conical symplectic variety of dimension $2n$ which has a projective symplectic resolution. Assume that 
$X$ admits an effective Hamiltonian action of an $n$-dimensional algebraic torus $T^n$, compatible with the conical $\mathbf{C}^*$-action. 
A typical example of $X$ is a toric hyperk\"{a}hler variety $Y(A, 0)$. In this article, we prove that this property characterizes $Y(A, 0)$ with 
$A$ unimodular. More precisely, if $(X, \omega)$ is such a conical symplectic variety, then there is a $T^n$-equivariant (complex analytic) 
isomorphism $\varphi: (X, \omega) \to (Y(A,0), \omega_{Y(A,0)})$ under which both moment maps are identified. Moreover, $\varphi$ sends the center $0_X$ of $X$ to the center $0_{Y(A,0)}$ of $Y(A,0)$. \vspace{0.3cm}

\end{abstract}

\begin{center}
{\bf Introduction}. 
\end{center}

In \cite{Del} Delzant characterized a compact (real) symplectic manifold $(M, \omega)$ of dimension $2n$ with an effective Hamiltonian $(S^1)^n$-action as a projective toric manifold. We prove a holomorphic analogue of Delzant's result in this article. 
   
Let $(X, \omega)$ be an affine symplectic variety of dimension $2n$ in the sense of \cite{Be}. 
The symplectic form $\omega$ on the smooth locus $X_{reg}$ of $X$ determines a Poisson structure on $X_{reg}$ and it uniquely 
extends to a Poisson structure $\{\:, \:\}$ on $X$. 
For a function $f \in \Gamma (X, \mathcal{O}_X)$, we define the Hamiltonian vector field by $H_f := \{f, \cdot\}$. This correspondence 
determines a map $H: \Gamma (X, \mathcal{O}_X) \to \Gamma (X, \Theta_X)$.    
Let $G$ be an algebraic group acting on $(X, \omega)$. For $a \in \mathfrak{g}$, such an action 
determines a vector field $\zeta_a$ on $X$, and this correspondence gives rise to a map $\zeta: \mathfrak{g} \to \Gamma (X, \Theta_X)$. 
The $G$-action is called Hamiltonian if $\zeta$ factorizes as 
$$\mathfrak{g} \stackrel{\mu^*}\to \Gamma (X, \mathcal{O}_X) \stackrel{H}\to \Gamma (X, \Theta_X).$$ 
Here $\mu^*$ is a $G$-equivariant map which satisfies $$\{\mu^*a, \mu^*b\} = \mu^*([a, b]), \:\:\: a, b \in \mathfrak{g}.$$   
Let $T^n$ be an $n$-dimensional algebraic torus.  
In this article we consider an affine symplectic variety $(X, \omega)$ of dimension $2n$ with an effective Hamiltonian $T^n$-action. 
A typical example of such a variety is a {\em toric hyperk\"{a}hler variety} (or a {\em hypertoric varietiy}) studied by \cite{Go}, \cite{BD}, 
\cite{HS}, \cite{Ko}, \cite{Pr} and others.     
Let $N$ be a positive integer with $n \le N$ and let $B$ be an integer valued $N \times n$-matrix such that $B$ determines an injection 
$\mathbf{Z}^n \to \mathbf{Z}^N$. We assume that each row vector of $B$ is primitive and $\mathrm{Coker}(B)$ is torsion-free. 
Then we have an exact sequence $$0 \to \mathbf{Z}^n \stackrel{B}\to \mathbf{Z}^N \stackrel{A}\to \mathbf{Z}^{N-n} \to 0$$
with an integer valued $(N-n) \times N$-matrix $A$. The exact sequence yields the exact sequence of algebraic tori 
$$1 \to T^{N-n} \to T^N \to T^n \to 1$$ in such a way that the induced exact sequence of character groups
$$0 \to \mathrm{Hom}_{alg.gp}(T^n, \mathbf{C}^*) \to \mathrm{Hom}_{alg.gp}(T^N, \mathbf{C}^*) \to \mathrm{Hom}_{alg.gp}(T^{N-n}, \mathbf{C}^*) \to 0$$ coincides with the given exact sequence. Let $(\mathbf{C}^{2N}, \omega_{st})$ be the pair of a $2N$-dimensional affine 
space with coordintes $(z_1, ..., z_N, w_1, ..., w_N)$ and a symplectic form $$\omega_{st} := \sum_{1 \le i \le N} dw_i \wedge dz_i.$$     
Then $T^N$ acts on $(\mathbf{C}^{2N}, \omega_{st})$ by $z_i \to t_iz_i \: (1 \le i \le N)$ and $w_i \to t_i^{-1}w_i \: (1 \le i \le N)$. 
This action is Hamiltonian and induces a Hamiltonian $T^{N-n}$-action on $(\mathbf{C}^{2N}, \omega_{st})$. Let $\mu: \mathbf{C}^{2N} \to \mathbf{C}^{N-n}$ be the moment map for this Hamiltonian $T^{N-n}$-action such that $\mu(0) = 0$. Then the toric hyperk\"{a}hler variety $Y(A, 0)$ is defined 
as $\mu^{-1}(0) /\hspace{-0.1cm}/ T^{N-n}$. By construction $Y(A,0)_{reg}$ admits a symplectic 2-form $\omega_{Y(A,0)}$ and 
$(Y(A, 0), \omega_{Y(A,0)})$ is an affine symplectic variety. 
Moreover, $Y(A, 0)$ admits a Hamiltonian $T^n$-action.   
If we define $Y(A, \alpha) := \mu^{-1}(0) /\hspace{-0.1cm}/_{\alpha} T^{N-n}$ for a generic 
$\alpha \in \mathrm{Hom}_{alg.gp}(T^{N-n}, \mathbf{C}^*)$, then  
$Y(A, \alpha)$ has only quotient singularities and 
the induced map $Y(A, \alpha) \to Y(A, 0)$ is a projective crepant partial resolution. When $A$ is unimodular (or equivalently, $B$ is 
unimodular), this map gives a projective symplectic resolution. When $N = n$ and $B$ is an isomorphism, we understand that 
$(Y(A, 0), \omega_{Y(A, 0)}) = (\mathbf{C}^{2N}, \omega_{st})$.     

An affine symplectic variety $(X, \omega)$ is conical if it admits a $\mathbf{C}^*$-action such that the coordinate ring $R$ of $X$ 
is positively graded, i.e. $R = \oplus_{i \ge 0}R_i$, $R_0 = \mathbf{C}$ and $\omega$ is homogeneous with respect to the 
$\mathbf{C}^*$-action. In this situation, $X$ has a unique fixed point $0_X$ corresponding to the maximal ideal $\oplus_{i > 0}R_i$ 
of $R$.  A toric hyperk\"{a}hler variety is a conical symplectic variety because the scaling $\mathbf{C}^*$-action on $\mathbf{C}^{2N}$: 
$$z_i \to tz_i, \:\:\: w_i \to tw_i \:\:\: (1 \le i \le N)$$ induces a conical $\mathbf{C}^*$-action on $Y(A, 0)$.  

It would be an interesting problem to characterize toric hyperk\"{a}hler varieties among conical symplectic varieties. In fact, it is a special case of a more general conjecture by Arbo and Proudfoot [\cite{A-P}, Conjecture 5.8].  
The main result of this article is the following, which can be regarded as a weak version of their conjecture: \vspace{0.2cm}

{\bf Theorem 5.8.} {\em Let $(X, \omega)$ be a conical symplectic variety of dimension $2n$ which has a projective symplectic resolution. 
Assume that $X$ admits an effective Hamiltonian action of an $n$-dimensional algebraic torus $T^n$, compatible with the conical 
$\mathbf{C}^*$-action. Then there is a $T^n$-equivariant (complex analytic) isomorphism 
$\varphi: (X, \omega) \to (Y(A, 0), \omega_{Y(A,0)})$ 
which makes the following diagram commutative}  
\begin{equation} 
\begin{CD} 
(X, \omega) @>{\varphi}>> (Y(A, 0), \omega_{Y(A,0)}) \\  
@V{\mu}VV @V{\bar{\mu}}VV  \\  
(\mathfrak{t}^n)^* @>{id}>> (\mathfrak{t}^n)^*.
\end{CD} 
\end{equation}   
{\em Here $A$ is unimodular and the vertical maps are moment maps for the $T^n$-actions.} 
{\em Moreover,  we have $\varphi (0_X) = 0_{Y(A, 0)}$.} 
\vspace{0.2cm}

At this moment $\varphi$ is not necessarily an algebraic isomorphism, but only a complex analytic isomorphism. We note that, in general,    
$\varphi$ cannot be chosen to be $\mathbf{C}^*$-equivariant even when $wt(\omega) = 2$ (cf. Remark \ref{(5.9)}).  
However, we have many different choices of the conical $\mathbf{C}^*$-actions on $(X, \omega)$, compatible with the $T^n$-actions. 
A natural question is the following (cf. also Question 5.10):   
\vspace{0.2cm}

{\bf Question}. {\em If necessary, after replacing the original conical $\mathbf{C}^*$-action on $(X, \omega)$ by a different conical $\mathbf{C}^*$-action on $(X, \omega)$, can we take $\varphi$ in a $\mathbf{C}^*$-equivariant way ?}
\vspace{0.2cm}

If this question is affirmative, then the map $\varphi$ is automatically an algebraic isomorphism.  
 
If $M$ is a complete hyperk\"{a}hler manifold of real dimension $4n$ with an effective tri-Hamiltonian $(S^1)^n$-action, and $M$ has 
Euclidean volume growth, then Bielawski \cite{B} has given a similar characterization. But our result does not assume the existence of a 
hyperk\"{a}hler metric.  
  
In the remainder we shall explain the strategy for proving Theorem (5.8) and how we use the conical $\mathbf{C}^*$-action on 
$X$. Our proof is based on Losev's work \cite{Lo}. As explained in \S 3, the moment map $\mu: X \to (\mathfrak{t}^n)^*$ is 
surjective and it coincides with the GIT quotient map of the $T^n$-action. Then we associate $\mu$ with an effective divisor $H \subset  
(\mathfrak{t}^n)^*$ called the discriminant divisor (cf. \ref{(3.3)}). $H$ has the form  
$$H = m_1H_1 + ... + m_kH_k + H_{k+1} + ... + H_d$$ with $m_i > 1\: (1 \le i \le k)$, where each $H_i$ is a hyperplane of  $(\mathfrak{t}^n)^*$ passing through the origin and 
defined as $$H_i := \{\eta \in (\mathfrak{t}^n)^*\: \vert \: \langle \mathbf{b}_i, \eta \rangle = 0\}$$ with some 
primitive elements $\mathbf{b}_i \in \mathrm{Hom}_{alg.gp}(\mathbf{C}^*, T^n) \subset \mathfrak{t}^n$. 

Take a point $\eta \in (\mathfrak{t}^n)^*$ and let $U$ be a sufficiently small open neighborhood of $\eta \in (\mathfrak{t}^n)^*$. We are interested in the local structure of $\mu$ 
around $\mu^{-1}(U)$. 

As proved in \ref{(3.4)}, there is a closed subset $F_X$ of $(\mathfrak{t}^n)^*$ with $\mathrm{Codim}_{(\mathfrak{t}^n)^*}F_X \geq 2$ 
such that $\eta \in (\mathfrak{t}^n)^* - F_X$ has a singular fiber $\mu^{-1}(\eta)$ if and only if $\eta \in H$.    
When $\eta \in (\mathfrak{t}^n)^* - F_X - H_1 - \cdot\cdot\cdot - H_k$, $\mu^{-1}(U)$ is smooth. Losev \cite{Lo} (cf. also Example \ref{(1.2)}, Theorem \ref{(1.3)}) has already given a normal form of the map $\mu^{-1}(U) \to U$. On the other hand,  $\mu^{-1}(U)$ has only  $A_{m_i-1}$-singularities when $\eta \in H_i - F_X$ for some $i$ with $1 \le i \le k$. In this case, we exhibit a normal form of the map $\mu^{-1}(U) \to U$ in Theorem \ref{(2.11)}.   
We put $N := \sum_{1 \le i \le k} m_i + d -k$ and prepare $N$ primitive vectors $$\mathbf{b}_1, ..., \mathbf{b}_1, ..., \mathbf{b}_k, ..., 
\mathbf{b}_k, \mathbf{b}_{k+1}, \mathbf{b}_{k+2}, ..., \mathbf{b}_d.$$ Here $\mathbf{b_i}$ appears in $m_i$ times when $i \le k$ and once when 
$i \geq k+1$. These vectors determine a map $B: \mathbf{Z}^n \to \mathbf{Z}^N$. 

Our initial plan was to construct a toric hyperk\"{a}hler variety 
$Y(A, 0)$ from this $B$ and compare $X$ with $Y(A, 0)$.  In fact, $\mu: X \to (\mathfrak{t}^n)^*$ and $\bar{\mu}: 
Y(A, 0) \to (\mathfrak{t}^n)^*$ have the same discriminant divisor $H \subset (\mathfrak{t}^n)^*$. For $Y(A, 0)$ we similarly 
define a closed subset $F_{Y(A,0)} \subset (\mathfrak{t}^n)^*$ and we put $F := F_X \cup F_{Y(A,0)}$ and define $(\mathfrak{t}^n)^{*, 0} := 
(\mathfrak{t}^n)^* - F$, $X^0 := \mu^{-1}((\mathfrak{t}^n)^{*, 0})$, and $Y(A, 0)^0 := \bar{\mu}^{-1}((\mathfrak{t}^n)^{*, 0})$. 
Then, as explained above, both $\mu$ and $\bar{\mu}$ have the same local form around  each $\eta \in 
(\mathfrak{t}^n)^{*, 0}$. As in \cite{Lo}, let ${\mathcal Aut}^{X^0}$ be the sheaf on $(\mathfrak{t}^n)^{*, 0}$ of the Hamiltonian automorphsms of 
$(X^0, \omega\vert_{X^0})$. Then the cohomology group $H^1((\mathfrak{t}^n)^{*, 0}, {\mathcal Aut}^{X^0})$ classifies the isomorphism classes of Hamiltonian $T^n$-spaces over $(\mathfrak{t}^n)^{*, 0}$ with a fixed discriminant divisor. By the exact sequence 
$$0 \to \mathbf{C} \oplus \mathrm{Hom}_{alg.gp}(T^n, \mathbf{C})^* \to \mathcal{O}_{(\mathfrak{t}^n)^{*, 0}} \to 
{\mathcal Aut}^{X^0} \to 0$$ we can compute $H^1((\mathfrak{t}^n)^{*, 0}, {\mathcal Aut}^{X^0})$. However, we may possibly have $H^1((\mathfrak{t}^n)^{*, 0}, \mathcal{O}_{(\mathfrak{t}^n)^{*, 0}}) \ne 0$ since $\mathrm{Codim}_{(\mathfrak{t}^n)^*}F \geq 2$, and this means that $$H^1((\mathfrak{t}^n)^{*, 0}, {\mathcal Aut}^{X^0}) \ne 0.$$ Unfortunately we have no idea how to see 
that $X^0$ and $Y(A, 0)^0$ both determine the same class\footnote{An exceptional case is when $X$ itself is smooth.  
In this case there is a $T^n$-equivarant isomorphism $(X, \omega) \cong (\mathbf{C}^{2n}, \omega_{st})$.} in $H^1((\mathfrak{t}^n)^{*, 0}, {\mathcal Aut}^{X^0})$.  
    
In order to compare $X$ and $Y(A, 0)$ in a general case, we need more information. Here we use the assumption that $(X, \omega)$ has a projective symplectic resolution $\pi: \tilde{X} \to X$. We assume that $X$ is singular. Then we have $r:= b_2(\tilde{X}) > 0$. 
Let $$f: (\tilde{\mathcal X}, \omega_{\tilde{\mathcal X}/\mathbf{C}^d}) \to \mathbf{C}^r$$ be the universal Poisson deformation of 
$(\tilde{X}, \omega_{\tilde X})$ (cf. \cite{Na 1}, \cite{Na 2}). 
We put $\mathcal{X} := \mathrm{Spec}\: \Gamma (\tilde{\mathcal X}, \mathcal{O}_{\tilde{\mathcal X}})$. 
Then we get a Poisson deformation $\bar{f}: (\mathcal{X}, \omega_{\mathcal{X}/\mathbf{C}^1}) \to \mathbf{C}^r$ of $\bar{f}^{-1}(0) = X$. 
There is a projective birational map $\Pi: \tilde{\mathcal X} \to \mathcal{X}$ over $\mathbf{C}^r$ and for a general point 
$t \in \mathbf{C}^r$, the map $\Pi_t: \tilde{\mathcal X}_t \to \mathcal{X}_t$ is an isomorphism. In particular, $\mathcal{X}_t$ is smooth
\footnote{This fact actually ensures that $B$ is unimodular. Since $\mathrm{Coker}(B)$ is torsion free,  
we can define $A$ to be the map from $\mathbf{Z}^N \to \mathrm{Coker}(B) = \mathbf{Z}^{N-n}$.} 
for such a point $t \in \mathbf{C}^d$.  
Since $\bar{f}$ is a family of symplectic varieties with Hamiltonian $T^n$-actions, we have 
a relative moment map $\mu_{\mathcal X}: \mathcal{X} \to (\mathfrak{t}^n)^* \times \mathbf{C}^d$ in such a way that 
$(\mu_{\mathcal X})_0: \mathcal{X}_0 \to (\mathfrak{t}^n)^* \times \{0\}$ coincides with $\mu$.   
Take a general line $\mathbf{C}^1 \subset \mathbf{C}^d$ passing through $0$ and pull back $\bar{f}$ to this line. Then we have 
a Poisson deformation $\mathcal{Z} \to \mathbf{C}^1$ of $X$. The fibers $\mathcal{Z}_t$ are smooth for all $t \ne 0$. 
Now the relative moment map $\mu_{\mathcal X}$ is restricted to the relative moment map $$\mu_{\mathcal Z}: \mathcal{Z} \to 
(\mathfrak{t}^n)^* \times \mathbf{C}^1.$$ There is an effective divisor $\mathcal{H}$ of $(\mathfrak{t}^n)^* \times \mathbf{C}^1$ such 
that $\mathcal{H}_t$ is the discriminant divisor for the moment map $\mu_{{\mathcal Z}_t}$ for each $t \in \mathbf{C}^1$. 
We call $\mathcal{H}$ the discriminant divisor of $\mu_{\mathcal Z}$. On the other hand, for the toric hyperk\"{a}hler $Y(A, 0)$, we similarly construct a Poisson deformation $\mathcal{Z}' \to \mathbf{C}^1$ of $Y(A, 0)$ so that the discriminant divisor $\mathcal{H}'$ for 
the relative moment map $\mu_{{\mathcal Z}'}: \mathcal{Z}' \to 
(\mathfrak{t}^n)^* \times \mathbf{C}^1$ satisfies  $$\mathcal{H}' = \mathcal{H}.$$ 
Recall that we have defined a closed subset $F$ of $(\mathfrak{t}^n)^*$ for $X$ and $Y(A,0)$.  
Now we regard $F$ as a subset of $(\mathfrak{t}^n)^* \times \{0\}$. 
We put $S := (\mathfrak{t}^n)^* \times \mathbf{C}^1$, $S^0 := (\mathfrak{t}^n)^* \times \mathbf{C}^1 -F$, 
$\mathcal{Z}^0 := \mu_{\mathcal Z}^{-1}(S^0)$ and $({\mathcal Z}')^0 := \mu_{{\mathcal Z}'}^{-1}(S^0)$. 
We compare $\mathcal{Z}^0 \to S^0$ with $({\mathcal Z}')^0 \to S^0$. This time we have $$\mathrm{Codim}_S F \geq 3.$$
These two spaces turn out to be isomorphic as $T^n$-Hamiltonian spaces over $S^0$ (Corollary \ref{(5.4)}). 
Note that $\mathcal{Z}$ and $\mathcal{Z}'$ are both Stein normal varieties. Since $\mathrm{Codim}_{\mathcal Z}(\mathcal{Z} - \mathcal{Z}^0) 
\geq  2$ and $\mathrm{Codim}_{{\mathcal Z}'}(\mathcal{Z}' - (\mathcal{Z}')^0) \geq 2$, this implies that $\mathcal{Z}$ and $\mathcal{Z}'$ 
are isomorphic as $T^n$-Hamiltonian spaces over $S$ (Theorem \ref{(5.5)}). If we restrict this isomorphism over 
$(\mathfrak{t}^n)^* \times \{0\} \subset S$, then we obtain Theorem \ref{(5.8)}.  

When $X$ does not have a projective symplectic resolution, we would need a different approach. For 
example, if $X$ has only {\bf Q}-factorial terminal singularities, $X$ is rigid in Poisson deformation; hence there is  
no substitute for $\mathcal{Z}$. Moreover, the definiton of a toric hyperk\"{a}hler variety should be slightly relaxed as in \cite{BD} so that 
$\mathrm{Coker}(B)$ is not necessarily torsion free.

{\bf Acknowledgement}. We thank the referee for valuable suggestions to improve the present article. In particular, the referee informed us of 
the conjecture of Arbo and Proudfoot.    
\vspace{0.3cm}

\section{}  Let $(M, \omega)$ be a complex symplectic manifold of dimension $2n$. We assume that $(M, \omega) $ admits a Hamiltonian action  of an algebraic torus $T$. Let $\mathfrak{t}$ be the Lie algebra of $T$. For each $a \in \mathfrak{t}$, 
the torus action determines a vector field $\zeta_a$ on $M$. 
By definition, there is a moment map $$\mu: M \to \mathfrak{t}^*,$$ which is $T$-equivariant and satisfies 
$$\omega_x(v, \zeta_a(x)) = \langle d\mu_x(v), a \rangle, \:\: x \in M, \:\: v \in T_xM, \:\: a \in \mathfrak{t}.$$ 
Here $d\mu_x$ is the tangential map $T_xM \to \mathfrak{t}^*$ induced by $\mu$. 

\begin{Lem}\label{(1.1)} 

(1) Every $T$-orbit $T\cdot x$ is contained in a fiber of $\mu$.

(2) $T\cdot x$ is an isotropic submanifold of $M$.
\end{Lem}

{\em Proof}. (1): Since $\mu$ is $T$-equivariant and $T$ acts trivially on $\mathfrak{t}^*$, the first statement is clear.

(2): Define a map $T \to M$ by $t \to t\cdot x$. Then it induces a map $\mathfrak{t} \to T_xM$. We denote by $\mathfrak{t}_*x$ its image. By (1) we have $$\mathfrak{t}_*x \subset \mathrm{Ker}(d\mu_x).$$
We show that $$\mathrm{Ker}(d\mu_x) = (\mathfrak{t}_*x)^{\perp_{\omega}}.$$ Here $(\mathfrak{t}_*x)^{\perp_{\omega}}$ is the orthogonal complement of $\mathfrak{t}_*x$ with respect to $\omega$. 
In fact, by the property of the moment map, we have $v \in \mathrm{Ker}(d\mu_x)$ if and only if $\omega_x(v, \zeta_a(x)) = 0$ for all
$a \in \mathfrak{t}$. Notice that $\mathfrak{t}_*x$ is the subspace of $T_xM$ generated by $\{\zeta_a(x)\}$. 
Therefore $v \in \mathrm{Ker}(d\mu_x)$ if and only if $v \in (\mathfrak{t}_*x)^{\perp_{\omega}}$. $\square$  \vspace{0.2cm} 

We assume in addition that 

i) $M$ is a Stein manifold of dimension $2n$, 

ii) $\dim T = n$ and $T$ acts effectively on $M$. 

In this situation, the moment map $\mu$ has been extensively studied by Losev \cite{Lo}. 
The starting point is the next example. 

\begin{Exam}\label{(1.2)} {\rm Consider two algebraic tori  $(\mathbf{C}^*)^k$ with coordinates $t:= (t_1, ..., t_k)$ and $(\mathbf{C}^*)^{n-k}$ with 
coordinates $\theta := (\theta_1, ..., \theta_{n-k})$. We denote by $T_0$ the first torus and denote by $T_1$ the second torus. Put 
$T := T_0 \times T_1$. 
Choose characters $\chi_1, ..., \chi_k \in \mathrm{Hom}_{alg.gp}(T_0, \mathbf{C}^*)$ in such a way that they form a basis of 
$\mathrm{Hom}_{alg.gp}(T_0, \mathbf{C}^*)$. Determine a $k$-dimensional $T_0$-representation $$V = \bigoplus_{1 \le i \le k} \mathbf{C}v_i, \:\:\:  \mathrm{by} \:\:\: t \cdot v_i := \chi_i(t)v_i.$$ Let $V^*$ be the dual representation of $V$ and let $v^1, ..., v^k$ be the dual basis. 
Let $\mathfrak{t}_1$ be the Lie algebra of $T_1$. Regarding $\theta_1, ..., \theta_{n-k}$ as linear functions on $T_1$, define a basis $\beta_1, ..., \beta_{n-k}$ of $\mathfrak{t}_1^*$ by 
$$\beta_i := \frac{d\theta_i}{\theta_i}.$$ Let $$\beta^1, ..., \beta^{n-k} \in \mathfrak{t}_1$$ be the dual basis. 
Then $T_1$ naturally acts on $T_1$ itself. Since $\beta_1, ..., \beta_{n-k}$ are $T_1$-invariant forms, $T_1$ acts trivially on 
$\mathfrak{t}_1^*$. Therefore $T$ acts on $T^*(V \times T_1) :=  V \times V^* \times T_1 \times \mathfrak{t}_1^*$. Notice that  
$$(v^1, ..., v^k, v_1, ..., v_k, \theta_1, ..., \theta_{n-k}, \beta^1, ..., \beta^{n-k})$$ are coordinates of $V \times V^* \times T_1 \times \mathfrak{t}_1^*$.  
Define a $T$-invariant 1-form $\alpha$ on $T^*(V \times T_1)$ by 
$$\alpha := \sum_{1 \le i \le k}v_idv^i + \sum_{1 \le j \le n-k}\beta^j \frac{d\theta_j}{\theta_j}.$$  
Then $$\omega := d\alpha = \sum_{1 \le i \le k}dv_i\ \wedge dv^i + \sum_{1 \le j \le n-k}d\beta^j \wedge \frac{d\theta_j}{\theta_j}$$ 
is a $T$-invariant symplectic 2-form on $T^*(V \times T_1)$.  
The $T$-action on $(T^*(V \times T_1), \omega)$ is actually a Hamiltonian action. Take 
$$\frac{d\chi_1}{\chi_1}, ..., \frac{d\chi_k}{\chi_k}, \beta_1, ..., \beta_{n-k}$$ as a basis of $\mathfrak{t}^*$. 
Then the moment map $\mu: V \times V^* \times T_1 \times \mathfrak{t}_1^* \to \mathfrak{t}^*$ is given by 
$$\mu (v^1, ..., v^k, v_1, ..., v_k, \theta_1, ..., \theta_{n-k}, \beta^1, ..., \beta^{n-k}) = (v^1v_1, ..., v^kv_k, \beta^1, ..., \beta^{n-k}) + \lambda$$ 
with a constant $\lambda \in \mathfrak{t}^*$.  
We regard $\mathfrak{t}_0$ (resp. $\mathfrak{t}_0^*$) as the $\mathbf{C}$-vector space of $T_0$-invariant vector fields (resp. $T_0$-invariant 1-forms).  Now $\frac{d\chi_1}{\chi_1}, ..., \frac{d\chi_k}{\chi_k}$ form a basis of $\mathfrak{t}_0^*$. One can embed the 
$\mathbf{Z}$-module $\mathrm{Hom}_{alg.gp}(T_0, \mathbf{C}^*)$ into $\mathfrak{t}_0^*$ by $\chi_i \to \frac{d\chi_i}{\chi_i}$. 
Let $\chi^1, ..., \chi^k \in \mathrm{Hom}_{alg.gp}(\mathbf{C}^*, T_0)$ be the dual basis of $\chi_1, ..., \chi_k$. 
Then we can embed the $\mathbf{Z}$-module $\mathrm{Hom}_{alg.gp}(\mathbf{C}^*, T_0)$ sending $\chi^1, ..., \chi^k$ 
to the dual basis of $\frac{d\chi_1}{\chi_1}, ..., \frac{d\chi_k}{\chi_k}$.  In this way $\chi^i$ is regarded as an element of $\mathfrak{t}_0$. In other words, $\chi^i$ is a linear function on $\mathfrak{t}_0^*$. By the surjection $\mathfrak{t}^* \to \mathfrak{t}_0^*$, $\chi^i$ is also regarded as a linear function on $\mathfrak{t}^*$.  
Now we put $$H_i := \{\beta \in \mathfrak{t}^* \:\vert \: \chi^i (\beta) = 0\}$$ and define 
$$D := \bigcup_{1 \le i \le k} (\lambda + H_i).$$ Then the moment map $\mu: T^*(V \times T_1) \to \mathfrak{t}^*$ has singular fibers exactly over $D \subset \mathfrak{t}^*$.  
Finally we add two observations. First every fiber of $\mu$ contains only finitely many $T$-orbits. Second a smooth fiber of $\mu$ consists of only one $T$-orbit with trivial stabilizer group.     $\square$  } 
\end{Exam}

A main point of \cite{Lo} is that, when $M$ satisfies (i) and (ii), the moment map $\mu$ is locally isomorphic to Example \ref{(1.2)}. 
Let $\tau: M \to N$ be the GIT quotient of $M$ by $T$ (cf. \cite{Sn}).        

\begin{Thm}\label{(1.3)}(\cite{Lo})

(1) The map $\mu$ factors through $N$: 
$$ M \stackrel{\tau}\to N \stackrel{\nu}\to \mathfrak{t}^*.$$ Moreover, $\nu$ is an etale map.

(2) There is a divisor $\mathcal{D}$ of $N$ such that $\tau: M \to N$ has singular fibers exactly over $\mathcal{D}$. 
For each $y \in N$, there is an open neighborhood $U$ of $y$ such that $\nu\vert_U : U \to \nu(U)$ is an isomorphism, $\mu\vert_{\tau^{-1}(U)}:  \tau^{-1}(U) \to \mathfrak{t}^*$ is isomorphic to Example \ref{(1.2)} localized around $\lambda := \nu(y)  
\in \mathfrak{t}^*$, and that $\mathcal{D} \cap U$ coincides with $D$. $\square$ 

\end{Thm}

\section{}
Let $(X, \omega)$ be an affine symplectic variety of dimension $2n$ with an effective Hamiltonian action of an $n$-dimensional 
algebraic torus $T^n$.  Let $\mu: X \to (\mathfrak{t}^n)^*$ be the moment map.  
We assume that $(X, \omega)$ has a projective symplectic resolution 
$\pi: (\tilde{X}, \omega_{\tilde X}) \to (X, \omega)$. By \cite{Ka} a singular symplectic variety $X$ is stratified into a finite number of symplectic 
leaves $Y$. Elements of $T^n$ acting trivially on $Y$ form a subgroup $H$ of $T^n$. Then $T_Y := T^n/H$ acts on $Y$. 
In Theorem \ref{(2.2)} we show that $\dim T_Y = \frac{1}{2} \dim Y$ and such an action is Hamiltonian. The moment map 
$\mu_Y$ for the $T_Y$-action is nothing but the restriction $\mu\vert_Y$ of the original moment map $\mu$ to $Y$. Moreover, 
every connected component of a general fiber of $\mu_Y$ is a closed $T^n$-orbit by Corollary \ref{(2.3)}.  
 
We next look at a symplectic leaf $Y$ of codimension $2$. It is well known that $X$ has Klein singularities along $Y$, but we prove here that such Klein singularities are of type $A$ in Corollary \ref{(2.5)}. Let $x \in Y$ be a point on a general fiber of $\mu_Y$ in the sense of Corollary \ref{(2.3)}. Then the stabilizer group $T^n_x \subset T^n$ coincides with $H$. We show that $H$ is connected in Lemma \ref{(2.9)}; hence $H = \mathbf{C}^*$. Therefore the complex analytic germ $(X, x)$ is a symplectic singularity with a $\mathbf{C}^*$-action. Proposition \ref{(2.8)} gives a description of $(X, x)$. In Theorem \ref{(2.11)} we give a normal form of the moment map $\mu: X \to (\mathfrak{t}^n)^*$ around $\mu (x) \in (\mathfrak{t}^n)^*$. We use Moser's argument for Darboux lemma and an analytic version of Luna's fundamental lemma to prove these results.            

\begin{Prop}\label{(2.1)} 

(1) Each fiber of $\mu: X \to (\mathfrak{t}^n)^*$ contains only finitely many $T^n$-orbits.

(2) The moment map  $\mu$ is a dominating map and each fiber has dimension $n$. 
\end{Prop}

{\em Proof}. (1): The $T^n$-action on $(X, \omega)$ extends to a Hamiltonian $T^n$-action on $(\tilde{X}, \omega_{\tilde X})$. The composite map $\tilde{\mu} := \mu \circ \pi$ is a moment map for the $T^n$-action on $(\tilde{X}, \omega_{\tilde X})$.
It is enough to show that each fiber of $\tilde{\mu}$ contains only finitely many $T^n$-orbits. 
By Sumihiro's theorem, $\tilde{X}$ is covered by a finitely many $T^n$-invariant affine open set $U_i$: 
$$\tilde X = \cup_{i \in I}U_i.$$
Restrict the map $\tilde{\mu}$ to $U_i$: 
$$\mu_i : U_i \to (\mathfrak{t}^n)^*.$$ Then $\mu_i$ is a moment map for the $T^n$-action on $(U_i, \omega_{\tilde X}\vert_{U_i})$. 
By Theorem \ref{(1.3)}, (1), $\mu_i$ factorizes as $$U_i \stackrel{\tau_i}\to U_i/\hspace{-0.1cm}/T^n \stackrel{\nu_i}\to (\mathfrak{t}^n)^*$$ and $\nu_i$ is an etale map. 
In particular, for $t \in (\mathfrak{t}^n)^*$, the fiber $\nu_i^{-1}(t)$ consists of finite points, say, $t_1, ..., t_m$. By the local description of 
$\tau_i$ (cf. Example \ref{(1.2)}, Theorem \ref{(1.3)}),  each fiber
$\tau_i^{-1}(t_j)$ contains only finitely many $T^n$-orbits. Therefore, each fiber of $\mu_i$ contains only finitely many $T^n$-orbits.
Since the index set $I$ of the open covering is finite, we see that each fiber of $\tilde{\mu}$ contains only finitely many $T^n$-orbits. 

(2): We first show that any fiber of $\tilde{\mu}$ has dimension $\le n$. In fact, suppose that some fiber has dimension $> n$. 
By Lemma \ref{(1.1)}, (1) this fiber is a union of $T^n$-orbits. Since each $T^n$-orbit is an isotropic submanifold of $\tilde{X}$ by Lemma \ref{(1.1)}, (2), it has dimension $\le n$. Hence the fiber has infinitely many $T^n$-orbits, which contradicts (1). 
 
If $\mu$ is not a dominating map, then $\tilde{\mu}$ is not a dominating map. Then every fiber of $\tilde{\mu}$ has dimension 
$ > n$. This is a contradiction. Therefore $\mu$ is a dominating map. Now we see that every fiber of $\mu$ has at least dimension $n$.    
If some fiber $\mu^{-1}(\eta)$ has dimension $> n$, then $\dim \tilde{\mu}^{-1}(\eta) > n$. This is a contradiction. Hence, every fiber of 
$\mu$ has dimension $n$.  $\square$   
\vspace{0.2cm}

A symplectic variety is stratified into a finite number of symplectic leaves \cite{Ka}. 
Let $Y \subset X$ be a symplectic leaf of dimension $2n-2r$. The symplectic form 
$\omega$ determines a Poisson structure on $X$. Then this Poisson structure is restricted to a Poisson structure on $Y$ and determines 
a symplectic form $\omega_Y$ on $Y$. The torus $T^n$ acts on $(Y, \omega_Y)$. We set 
$$H := \{t \in T^n\: \vert \: t \: \mathrm{acts} \: \mathrm{on} \: Y \: \mathrm{trivially}\} $$ and put $T_Y := T^n/H$. By definition, $T_Y$ acts effectively on $Y$. Then we have:

\begin{Thm}\label{(2.2)}

(1) The action of $T_Y$ on $(Y, \omega_Y)$ is a Hamiltonian action and we have a commutative diagram of 
moment maps  
\begin{equation} 
\begin{CD} 
Y @>>> X \\ 
@V{\mu_Y}VV @V{\mu}VV \\ 
\mathfrak{t}_Y^* @>>> (\mathfrak{t}^n)^*.     
\end{CD} 
\end{equation}  

(2) $\dim T_Y = n -r$. 

(3) $\mu_Y$ is a dominating map.  

\end{Thm}

{\em Proof}.  (1): The moment map $\mu$ is regarded as a moment map for the Hamiltonian $T^n$-action on the Poisson variety $(X, \{\:, \:\})$. Since $(Y,  \{\:, \:\}_Y)$ is a Poisson subvariety of $(X, \{\:, \:\})$, the $T^n$-action on $(Y, \omega_Y)$ is also Hamiltonian and we get 
the commutative diagram of moment maps. 

(2):  We first notrice that $\dim T_Y \le 1/2 \cdot \dim Y (= n -r)$. In fact, since $T_Y$ acts effectively on $Y$, the stabilizer group 
of a general $T_Y$-orbit in $Y$ is trivial. Assume that $\dim T_Y > 1/2 \cdot \dim Y$. Then a general $T_Y$-orbit has dimension $> 1/2 \cdot 
\dim Y$. But each $T_Y$-orbit is an isotropic submanifold of $(Y, \omega)$. This is a contradiction. 

We next assume that $\dim T_Y < n - r$. In this case each fiber of $\mu_Y$ has dimension $> n-r$. By Lemma (1.1), each $T_Y$-orbit 
of $Y$ is an isotropic submanifold of $Y$, which is contained in a fiber of $\mu_Y$.
This means that a fiber of $\mu_Y$ contains infinitely 
many $T_Y$-orbits (hence $T^n$-orbits). But this contradicts Proposition \ref{(2.1)}; hence we get (2). 

(3): We take a $T_Y$-invariant affine open subset $Y^0$ of $Y$ and apply Theorem \ref{(1.3)} by putting $M = Y^0$. Then 
$\mu\vert_{Y^0}$ is a dominating map. Hence $\mu_Y$ is also a dominating map.  $\square$

\begin{Cor}\label{(2.3)}
Every connected component of a  general fiber of $\mu_Y$ is a closed $T^n$-orbit.
\end{Cor}

{\em Proof}. If we take $t \in \mathrm{Im}(\mu_Y)$ general, then, for any $y \in \mu_Y^{-1}(t)$, 
the closure $\overline{T^n\cdot y}$ is contained in $Y$. In fact, choose $t$ so that $t \notin \mathfrak{t}_{Y'}^*$ for any 
symplectic leaf $Y'$ with $Y' \subset \bar{Y}$. This is possible because $\dim \mathfrak{t}_{Y'}^* < \dim \mathfrak{t}_Y^*$ and 
$\mu_Y$ is dominating.  Take $z \in \overline{T^n\cdot y}$.
If $z \notin Y$, then $z$ is contained in a smaller 
symplectic leaf $Y'$ such that $Y' \subset \bar{Y}$. By Theorem \ref{(2.2)}, a quotient torus $T_{Y'}$ of $T^n$ acts effectively on $Y'$ with the 
moment map $\mu_{Y'}$. Then $\mu (z) = \mu_{Y'}(z) \in \mathfrak{t}_{Y'}^*$. Since $z \in \overline{T^n\cdot y}$, we have 
$\mu (y) = \mu(z) \in \mathfrak{t}_{Y'}^*$. This contradicts the choice of $t$. Therefore $\overline{T^n\cdot y} \subset Y$.   

We can also take $t \in \mathrm{Im}(\mu_Y)$ so that 
$\mu_Y^{-1}(t)$ is smooth. This is possible because $Y$ and $\mathfrak{t}_Y^*$ are both smooth. 
Take a point $y \in \mu_Y^{-1}(t)$. We prove that $T^n \cdot y (= T_Y \cdot y)$ is a closed orbit and coincides with a connected component of $\mu_Y^{-1}(t)$.  
The symplectic leaf $Y$ is covered by $T_Y$-invariant smooth affine open subsets $Y_i, \: (i \in I)$. Then the moment map $\mu_{Y_i}$ for the $T_Y$-action on $Y_i$ is nothing but the restriction of $\mu_Y$ to $Y_i$.    
Choose an $i \in I$ so that  $y \in \mu_{Y_i}^{-1}(t)$. Then 
$\mu_{Y_i}^{-1}(t)$ is smooth because it is a non-empty open subset of $\mu_Y^{-1}(t)$.   
Then we see, by Theorem \ref{(1.3)} and Example \ref{(1.2)}, that each connected component of $\mu_{Y_i}^{-1}(t)$ consists of a $T_Y$-orbit with trivial stabilizer. 
In particular, $T^n \cdot y$ is a connected component of $\mu_{Y_i}^{-1}(t)$ and $\dim (T^n\cdot y) = \dim T_Y$. Moreover, $T^n\cdot y$ is a closed orbit in $X$. In fact, suppose to the contrary that $T^n\cdot y \ne \overline{T^n\cdot y}$. Then, since $\overline{T^n\cdot y} \subset Y$, any point $z \in \overline{T^n\cdot y} - T^n\cdot y$ is contained in some other $Y_j$. Since $\mu_Y(z) = t$, we have $z \in \mu_{Y_j}^{-1}(t)$. Then, again by Theorem \ref{(1.3)} and Example \ref{(1.2)}, $T^n\cdot z$ coincides with a connected component of $\mu_{Y_j}^{-1}(t)$ and $\dim (T^n\cdot z) = \dim T_Y$. On the other hand, since $T^n\cdot z \subset \overline{T^n\cdot y} - T^n\cdot y$, we must have 
$\dim (T^n \cdot z) < \dim (T^n \cdot y)$. This is a contradiction. Therefore $T^n\cdot y$ is closed in 
$X$. Note that $T^n\cdot y$ is an open subset of  $\mu_{Y_i}^{-1}(t)$ because it is a connected component of $\mu_{Y_i}^{-1}(t)$. Since $\mu_{Y_i}^{-1}(t)$ is an open subset of $\mu_Y^{-1}(t)$, $T^n\cdot y$ is open in $\mu_Y^{-1}(t)$. Hence $T^n\cdot y$ is a connected component of $\mu_Y^{-1}(t)$.   $\square$
\vspace{0.2cm}

In the remainder we assume that $Y \subset X$ is a symplectic leaf of codimension $2$.   
By Theorem \ref{(2.2)}, $\dim T_Y = n -1$ and $\dim H = 1$. By Corollary \ref{(2.3)}, a connected component of a general fiber of $\mu_Y$ is a closed $T_Y$ orbit with trivial stabilizer group. Let $x \in Y$ be a point on a general fiber of 
$\mu_Y$. Let $T^n_x \subset T^n$ be the stabilizer group of $x$. Then $T^n_x = H$. We shall prove that $H$ can be written as 
$$ H = G \times \mathbf{C}^*$$ with a finite abelian group $G$. First note that the identity component $H^0$ of $H$ is a 1-dimensional torus because $\dim H = 1$. We then have a commutative diagram of exact sequences of algebraic groups: 
\begin{equation} 
\begin{CD} 
1 @>>> H @>>> T^n @>>> T_Y @>>> 1 \\ 
@. @AAA  @A{id}AA @AAA @. \\ 
1 @>>> \mathbf{C}^*  @>>> T^n @>>> (\mathbf{C}^*)^{n-1} @>>> 1.      
\end{CD} 
\end{equation}  
Here the vertical map on the left hand side is the inclusion map of the identity component $\mathbf{C}^* (= H^0)$ into $H$, and 
$(\mathbf{C}^*)^{n-1} = T^n/\mathbf{C}^*$. On the second row there is an injective homomorphism $(\mathbf{C}^*)^{n-1} \to T^n$ 
which splits the exact sequence. Then the composite $(\mathbf{C}^*)^{n-1} \to T^n \to T_Y$ is a surjection. Let $G$ be its kernel. 
Then $G$ is a finite abelian subgroup of $T^n$ and we get $H = G \times \mathbf{C}^*$.    

We first look at the $H^0$-action on the complex analytic germ $(X, x)$. 
Since $X$ has Klein singularities along $Y$, we have an isomorphism  
$$\phi: (X, x) \cong (\mathbf{C}^2/\Gamma, 0) \times (\mathbf{C}^{2n-2}, 0),$$ 
where $\Gamma$ is a finite subgroup of $SL(2, \mathbf{C})$. Let $(z_1, z_2)$ be the coordinates of $\mathbf{C}^2$ and let 
$(t_1, ..., t_{2n-2})$ be the coordinates of $\mathbf{C}^{2n-2}$. 
The symplectic form $\omega_{\mathbf{C}^2} := dz_1 \wedge dz_2$ on $\mathbf{C}^2$ descends to 
a symplectic form on $\mathbf{C}^2/\Gamma - \{0\}$, which we denote by $\omega_{\mathbf{C}^2/\Gamma}$. On the other hand, 
we put $$\omega_{st} := dt_1 \wedge dt_n + ... + dt_{n-1} \wedge dt_{2n-2}.$$ Recall that $(X, x)$ admits a symplectic form $\omega$ 
(on the regular part). 
By Darboux lemma (cf. [\cite{Na 2}, Lemma 1.3]) we can take $\phi$ in such a way that $$\omega = \phi^*(\omega_{\mathbf{C}^2/\Gamma} + \omega_{st}).$$ Consider the orbit $T^nx \subset X$. Then $T^n x$ is contained in the symlpectic leaf $Y$.   
Notice that $(Y, x)$ is identified with the subvariety $$\{0\} \times (\mathbf{C}^{2n-2}, 0) \subset 
(\mathbf{C}^2/\Gamma, 0) \times (\mathbf{C}^{2n-2}, 0)$$ by $\phi$. Therefore $(T^nx, x)$ is identified with an ($n - 1$)-dimensional 
subvariety of $\{0\} \times (\mathbf{C}^{2n-2}, 0)$.   
Let $x \in U \subset X$ be a sufficiently small open neighborhood of $x$. Then the universal covering of $U_{reg}$ induces a 
finite Galois cover $\Pi: (Z, z) \to (X, x)$ with Galois group $\Gamma$. Note that $Z$ is smooth. If we put $\omega_Z := \Pi^*\omega$, 
then $\omega_Z$ is a symplectic form on $Z$. We can lift $\phi$ to a $\Gamma$-equivariant isomorphism $\tilde{\phi} : (Z, z) \to (\mathbf{C}^2, 0) \times 
(\mathbf{C}^{2n-2}, 0)$ so that $\omega_Z = \tilde{\phi}^*(\omega_{\mathbf{C}^2} + \omega_{st})$ and the following diagram 
commutes 
\begin{equation} 
\begin{CD} 
(Z, z)  @>{\tilde{\phi}}>> (\mathbf{C}^2, 0) \times (\mathbf{C}^{2n-2}, 0) \\  
@V{\Pi}VV  @VVV \\ 
(X, x)  @>{\phi}>> (\mathbf{C}^2/\Gamma, 0) \times (\mathbf{C}^{2n-2}, 0).      
\end{CD} 
\end{equation}  
Let us consider the subvariety $\{0\} \times (\mathbf{C}^{2n-2}, 0) \subset (\mathbf{C}^2, 0) \times (\mathbf{C}^{2n-2}, 0)$, which is isomorphically mapped onto the subvariety $\{0\} \times (\mathbf{C}^{2n-2}, 0) \subset 
(\mathbf{C}^2/\Gamma, 0) \times (\mathbf{C}^{2n-2}, 0)$ by the vertical map on the right hand side.  
Now we can find a subvariety $W \subset \tilde{\phi}^{-1}(\{0\} \times (\mathbf{C}^{2n-2}, 0))$ so that  
$\Pi (W) = T^nx$ and $\Pi\vert_W: (W, z) \to (T^nx, x)$ is an isomorphism.    

We take a group extension   
$$ 1 \to \Gamma \to \tilde{H} \to H^0 \to 1$$ so that $\tilde{H}$ acts on $(Z, z)$. 
Let $\tilde{H}^0$ be the identity component of 
$\tilde{H}$. Then $\tilde{H}^0$ is a 1-dimensional algebraic torus and the induced map $\tilde{H}^0 \to \mathbf{C}^*$ is a surjection 
of 1-dimensional algebraic tori. The $\Gamma$-action on $(Z, z)$ and the $\tilde{H}^0$-action on $(Z, z)$ are compatible. Hence 
$\Gamma \times \tilde{H}^0$ acts on $(Z, z)$. The tangent space $T_zZ$ is a $\Gamma \times \tilde{H}^0$-representation. 
We take a $\Gamma \times \tilde{H}^0$-equivariant isomorphism 
$$\varphi: (Z, z) \cong  (T_zZ, 0)$$ in such a way that its tangent map $d\varphi_z: T_zZ \to T_zZ$ is the identity map. 
Note that $T_zZ$ admits the symplectic form $\omega_Z(z)$ and $\Gamma \times \tilde{H}^0$ preserves $\omega_Z(z)$. 
Let us consider the subspace $T_zW \subset T_zZ$. Since ${H}^0$ acts trivially on $T^nx$, $\tilde{H}^0$ acts trivially on $W$; hence  
acts trivially on $T_zW$. $\Gamma$ also acts trivially on $T_zW$.  
Identify $\tilde{H}^0$ with $\mathbf{C}^*$. 

\begin{Lem}\label{(2.4)} 
We have $$T_zZ = \mathbf{C}(1) \oplus \mathbf{C}(-1) \oplus \mathbf{C}(0)^{\oplus 2n -2}$$ as a 
$\mathbf{C}^*$-representation. Here $\mathbf{C}(i)$ is a weight $i$ eigenspace for $i \in \mathbf{Z}$. 
\end{Lem}

{\em Proof}. As remarked above, $T_zW$ is a trivial $\Gamma \times \tilde{H}^0$-module. We prove that 
$T_zZ/(T_zW)^{\perp_{\omega_Z(z)}}$ is also 
a trivial $\Gamma \times \tilde{H}^0$-module. In fact, take $v \in T_zZ$ and consider an element $hv - v$ for $h \in \Gamma \times 
\tilde{H}^0$. 
For $w \in T_zW$, we have $$\omega_Z(z) (hv - v, w) = \omega_Z(z)(hv, w) - \omega_Z(z)(v, w) $$ 
$$ = \omega_Z(z) (v, h^{-1}w) - \omega_Z(z)(v, w) = \omega_Z(z)(v, w) - \omega_Z(z)(v, w) = 0.$$
This means that $$\overline{hv} = \overline{v} \in T_zZ/(T_zW)^{\perp_{\omega_Z(z)}}.$$ 
We have a direct sum decomposition of $\Gamma \times \tilde{H}^0$-modules 
$$T_zZ = (T_zW)^{\perp_{\omega_Z(z)}}/T_zW  \:\: \oplus \:\: T_zW \:\: \oplus \:\: T_zZ/(T_zW)^{\perp_{\omega_Z(z)}}.$$ 
The second and the third factors are trivial $\Gamma \times \tilde{H}^0$-module. On the other hand, 
$\tilde{H}^0$ acts effectively on $T_zZ$.    
In fact, $H^0$ acts effectively on $(X, x)$. By construction $\tilde{H}^0$ also acts effectively on $(Z, z)$; hence acts effectively on $(T_zZ, 0)$. Therefore $\tilde{H}^0$ acts effectively on the first factor. By definition, the first factor $V := (T_zW)^{\perp_{\omega_Z(z)}}/T_zW $ is a symplectic vector space of dimension $2$. Hence $\tilde{H}^0 \subset Sp(V)$ is a maximal torus. All maximal tori of $Sp(V)$ are conjugate 
to each other. In particular, we see that $V = \mathbf{C}(1) \oplus \mathbf{C}(-1)$ as an $\tilde{H}^0$-module. $\square$ 
\vspace{0.2cm}

We put $V := (T_zW)^{\perp_{\omega_Z(z)}}/T_zW$.  

\begin{Cor}\label{(2.5)}
$\Gamma$ is a finite cyclic group. In particular, $(\mathbf{C}^2/\Gamma, 0)$ is a Klein singularity of type $A$. 
\end{Cor}

{\em Proof}. If $\Gamma$ is not a cyclic group, then $V$ is an irreducible $\Gamma$-module. On the other hand, since the $\mathbf{C}^*$-action and the $\Gamma$-action on $V$ are compatible,  any element $t \in \mathbf{C}^*$ determines a $\Gamma$-equivariant isomorphism of $V$. By Lemma \ref{(2.4)} this isomorphism is not of the form $\alpha I_V$ ($\alpha \in \mathbf{C}$) for a general $t \in \mathbf{C}^*$. This is a contradiction by Schur's lemma. $\square$

\begin{Cor}\label{(2.6)} 
$\Gamma \subset\tilde{H}^0$. In particular, $\tilde{H} = \tilde{H}^0$. 
\end{Cor}

{\em Proof}. By Corollary \ref{(2.5)} $\Gamma = \mathbf{Z}/m \mathbf{Z}$ for some $m > 1$.  $\mathbf{C}(1)$ and $\mathbf{C}(-1)$ are both $\Gamma$-representations. In other words, $\Gamma$ acts on $\mathbf{C}(i)$ by $\rho_i: \mathbf{Z}/m\mathbf{Z} \to GL(\mathbf{C}(i))$ for $i = 1, -1$. One can write $\rho_1(\bar{1}) = \zeta$ with a $m$-th root of unity $\zeta$. Since $\Gamma \subset Sp(V)$, one has $\rho_{-1}(\bar{1}) = \zeta^{-1}$. Moreover, since $\Gamma$ acts effectively on $V$, $\zeta$ must be a primitive $m$-th root of unity. 
Then $\Gamma$ is a subgroup of $\tilde{H}^0$. $\square$ \vspace{0.2cm}

We identify $T_zZ$ with $T_0\mathbf{C}^2 \oplus T_0\mathbf{C}^{2n-2}$ by the tangential map $d\tilde{\phi}_z$ of $\tilde{\phi}$ at $z$.
Then $(d\tilde{\phi}_z)^{-1}$ induces a $\Gamma$-equivariant injection
$$T_0\mathbf{C}^2 \to T_zZ$$ whose image coincides with the (unique) 2-dimensional non-trivial $\Gamma$-factor 
$(T_zW)^{\perp_{\omega_Z(z)}}/T_zW$ of the $\Gamma$-representation $T_zZ$. 
Since $T_zZ$ is a $\tilde{H} (= \mathbf{C}^*)$-representation, we regard $T_0\mathbf{C}^2 \oplus T_0\mathbf{C}^{2n-2}$ as a $\mathbf{C}^*$-representation by $d\tilde{\phi}_z$. Since $\omega_Z(z) = d\tilde{\phi}_z^*(\omega_{\mathbf{C}^2}(0) + \omega_{st}(0))$ and 
the $\mathbf{C}^*$-action preserves $\omega_Z(z)$, the symplectic form $\omega_{\mathbf{C}^2}(0) + \omega_{st}(0)$ is preserved by 
the $\mathbf{C}^*$-action on $T_0\mathbf{C}^2 \oplus T_0\mathbf{C}^{2n-2}$. As remarked just before Corollary \ref{(2.5)}, $V$ corresponds to 
$T_0\mathbf{C}^2$ by $d\tilde{\phi}_z$. Since $V$ is a $\mathbf{C}^*$-submodule of $T_zZ$, $V^{\perp_{\omega_Z(z)}}$ is also 
a $\mathbf{C}^*$-submodule. $V^{\perp_{\omega_Z(z)}}$ corresponds to $T_0\mathbf{C}^{2n-2}$ by $d\tilde{\phi}_z$. 
Therefore $T_0\mathbf{C}^{2n-2}$  is a $\mathbf{C}^*$-submodule of $T_0\mathbf{C}^2 \oplus T_0\mathbf{C}^{2n-2}$. 
We have $$T_0\mathbf{C}^2 = \mathbf{C}(1) \oplus \mathbf{C}(-1), \:\:\:\: T_0\mathbf{C}^{2n-2} = 
\mathbf{C}(0)^{\oplus 2n-2}.$$ 

In the argument above, we have taken a $\Gamma \times \tilde{H}^0$-equivariant isomorphism $\varphi: (Z, z) \cong (T_zZ, 0)$ such that 
$d\varphi_z = id$. We finally show that there is a $\mathbf{C}^*$-equivariant automorphism $\psi$ of $(T_zZ, 0)$ 
such that $\varphi' := \psi \circ \varphi$ satisfies  $\omega_Z = (\varphi')^*\omega_Z(z)$. 
We put $\omega_1 := \omega_Z(z)$ and $\omega_2 := (\varphi^{-1})^*\omega_Z$. They are symplectic 2-forms on the germ 
$(T_zZ, 0)$ such that $\omega_1(0) = \omega_2(0)$ because $d\varphi_z = id$. The $\Gamma \times \tilde{H}^0$-action preserves both $\omega_1$ and $\omega_2$. Since $\Gamma \subset \tilde{H}^0$, we may regard the $\Gamma \times \tilde{H}^0$-action simply as the $\tilde{H}^0$-action. 
In this situtaion, the following equivariant Darboux lemma holds. Then $\psi$ of Lemma \ref{(2.7)} is a desired one. 

\begin{Lem}\label{(2.7)} 
There is an $\tilde{H}^0$-equivariant automorphism $\psi$ of $(T_zZ, 0)$ such that 
$\omega_2 = \psi^*\omega_1$. 
\end{Lem}

{\em Proof}. The basic strategy of the proof is the same as the usual Darboux lemma. We put $u := \omega_1 - \omega_2$. 
We find a $\tilde{H}^0$-invariant 1-form $\alpha$ on a small open neighborhood $0 \in U \subset T_zZ$ such that 
1) $d\alpha = u$, and 2) $\alpha (0) = 0$. Once such an $\alpha$ exists, we can prove the lemma in the same manner as in the usual Darboux lemma. 

Let us consider the scaling action of $\mathbf{R}_{>0}$ on 
$T_zZ$ defined by $a_t(v) := tv$, $t \in \mathbf{R}_{>0}$. We assume that $tU \subset U$ for any $t \in (0, 1]$. 
Let $\xi$ be a vector field 
on $U$ determined by this action.   
When $t \to 0$, $a_t$ goes to the constant map $a_0: U \to \{0\} \subset U$. Note that $a_0^*u = u(0) = 0$.  
We now have 
$$u = a_1^*u = \int_{0}^{1} L_{\xi} a_t^*u \: dt + a_0^*u = d\int_{0}^{1} \xi \rfloor a_t^*u dt.$$
Then we can take $$\alpha = \int_{0}^{1}\xi \rfloor a_t^*u \: dt.$$ $\square$  \vspace{0.2cm}

Summing up the arguments above, we have a sequence of $\mathbf{C}^*$-equivariant isomorphisms of symplectic singularities:  
$$((Z, z), \omega_Z) \stackrel{\varphi'}\to ((T_zZ, 0), \omega_Z(z)) \stackrel{d\tilde{\phi}_z}\to 
((T_0\mathbf{C}^2 \oplus T_0\mathbf{C}^{2n-2}, 0), \omega_{\mathbf{C}^2}(0) + \omega_{st}(0)).$$
Here the $\mathbf{C}^*$-action on the leftmost  is the $\tilde{H}^0$-action on $(Z, z)$. 
We can naturally identify $((T_0\mathbf{C}^2 \oplus T_0\mathbf{C}^{2n-2}, 0), \omega_{\mathbf{C}^2}(0) + \omega_{st}(0))$ 
with $((\mathbf{C}^2, 0) \times (\mathbf{C}^{2n-2},0), \omega_{\mathbf{C}^2} + \omega_{st})$.  Therefore we have a $\mathbf{C}^*$-equivariant isomorphism 
$$((Z, z), \omega_Z) \cong ((\mathbf{C}^2, 0) \times (\mathbf{C}^{2n-2},0), \omega_{\mathbf{C}^2} + \omega_{st}).$$
As is seen in Corollary \ref{(2.6)}, $\Gamma$ is contained in $\tilde{H}^0$. 
Let $\Gamma = \mathbf{Z}/m\mathbf{Z}$. Assume that the $\Gamma$-action on $(z_1, z_2) \in \mathbf{C}^2$ is given by 
$$\rho: \mathbf{Z}/m\mathbf{Z} \to SL(2, \mathbf{C}), \:\: \rho (\bar{i}) := 
\left(\begin{array}{cccccccc} 
 \zeta^i  & 0 \\
 0 & \zeta^{-i}   
\end{array}\right)$$
with a primitive $m$-th root $\zeta$ of unity. We put $x_1 := z_1^m$, $x_2 := z_2^m$ and $x_3 := z_1z_2$. Then 
$x_1, \: x_2, \: x_3$ are regarded as a function on $\mathbf{C}^2/\Gamma$. Then $\mathbf{C}^2/\Gamma$ is embedded in $\mathbf{C}^3$ 
as the subvariety defined by $x_1x_2 = x_3^m$.  Define 
$$\omega_{\mathbf{C}^2/\Gamma} := \mathrm{Res}(\frac{dx_1 \wedge dx_2 \wedge dx_3}{x_1x_2 - x_3^m}).$$
Let us consider the $\mathbf{C}^*$-equivariant isomorphism above. By taking the quotient of both sides by $\Gamma$, we get

\begin{Prop}\label{(2.8)}
There is a $\mathbf{C}^*$-equivariant isomorphism of symplectic singularities:
$$((X, x), \: \omega) \cong ((\mathbf{C}^2/\Gamma, 0) \times 
(\mathbf{C}^{2n-2}, 0), \:\: \omega_{\mathbf{C}^2/\Gamma} + \omega_{st}).$$ 
Here the $\mathbf{C}^*$ action on the left hand side is the $H^0$-action on $(X, x)$ and the $\mathbf{C}^*$-action on the right hand side is given by $$(x_1, x_2, x_3, t_1, ..., t_{2n-2}) \to (tx_1, t^{-1}x_2, x_3, t_1, ..., t_{2n-2}), \:\: t \in \mathbf{C}^*.$$ $\square$
\end{Prop}

Now that we have proved Proposition \ref{(2.8)}, the final goal of this section is to prove Theorem \ref{(2.11)}.  
For this purpose we need to look at the $GIT$-quotient map $X \to X/\hspace{-0.1cm}/H$ aroun\d a point 
$x \in Y$. 

\begin{Lem}\label{(2.9)} 
$H$ is connected; namely $H = H^0$. 
\end{Lem} 

{\em Proof}. 
As already remarked, $H = G \times \mathbf{C}^*$ with 
a finite abelian group $G$. Let us consider $\Pi:  (Z, z) \to (X, x)$ and take a group extension 
$$1 \to \Gamma \to K \stackrel{p}\to H \to 1$$ so that $K$ acts on $(Z, z)$. There is a commutative diagram 
\begin{equation} 
\begin{CD} 
 1 @>>> \Gamma   @>>> \tilde{H} @>>> H^0 @>>> 1 \\ 
 @. @V{id}VV  @VVV   @VVV \\  
1 @>>> \Gamma @>>> K @>>>  H @>>> 1.       
\end{CD} 
\end{equation}  
Here the vertical maps are all inclusions. 
Therefore $K$ contains a 1-dimensional torus $\tilde{H} = \mathbf{C}^*$ and $p\vert_{\tilde{H}}: \tilde{H} \to H^0$ is 
nothing but the map $\mathbf{C}^* \to \mathbf{C}^*, \:\: t \to t^m$. 
We identify $(X, x)$ with $(\mathbf{C}^2/\Gamma, 0) \times (\mathbf{C}^{2n-2}, 0)$ and 
$(Z, z)$ with $(\mathbf{C}^2, 0) \times (\mathbf{C}^{2n-2}, 0)$.   
By the definition of $H$, any element of $H$ acts 
trivially on the symplectic leaf $\{0\} \times (\mathbf{C}^{2n-2}, 0)$ of $(\mathbf{C}^2/\Gamma, 0) \times (\mathbf{C}^{2n-2}, 0)$; hence $K$ acts trivially on the subspace $\{0\} \times (\mathbf{C}^{2n-2}, 0)$ of $(\mathbf{C}^2, 0) \times (\mathbf{C}^{2n-2}, 0)$. 
We can take a $K$-equivariant isomorphism 
$(Z, z) \cong (T_zZ, 0)$. We look at the $K$-action on $T_zZ = T_0\mathbf{C}^2 \oplus T_0\mathbf{C}^{2n-2}$. 
For $g \in G$ with $g \ne 1$, choose $\tilde{g} \in K$ so that  $p(\tilde g) = (g, 1)$. 
Since $\tilde{g}$ acts trivially on $T_0\mathbf{C}^{2n-2} \subset T_0\mathbf{C}^2 \oplus T_0\mathbf{C}^{2n-2}$, 
$\tilde {g}$ has a form $$\left(\begin{array}{cccccccc} 
 A  & 0 \\
 C & I_{2n-2}    
\end{array}\right)$$
with a $2 \times 2$-matrix $A$ and a $(2n-2) \times 2$-matrix $C$. 
Since $\tilde{g} \in Sp(2n)$, we see that $C = 0$. 
On the other hand, $t \in \tilde{H}$ acts on $T_0\mathbf{C}^2 \oplus T_0\mathbf{C}^{2n-2}$ as a matrix 
$$\left(\begin{array}{cccccccc} 
 t  & 0 & 0 & ... &... & 0\\
 0 & t^{-1} & 0 & ... & ... & 0 \\
 0 & 0 & 1 & 0 & ... & 0  \\
 ... & ... & 0 & 1 & ... & 0 \\
 ... & ... & ... & ... & ... & ... \\
 0 & 0 & 0 & ... & 0 & 1
\end{array}\right).$$ 
Since $t \cdot \tilde{g} = \tilde{g} \cdot t$, we see that $$A = \left(\begin{array}{cccccccc} 
 a  & 0 \\
 0 &  d 
\end{array}\right)\:\:\: a, \: d \in \mathbf{C}, \:\: ad = 1,$$ which implies that 
$\tilde{g} \in \tilde{H}$. In particular $p(\tilde{g}) \in H^0$. This contradicts the choice of $\tilde{g}$. $\square$  
\vspace{0.2cm}

Let $\tau: X \to X/\hspace{-0.1cm}/H$ be the GIT quotient of $X$ by $H$. An open subset of $X$ (in the Euclidian topology) is called 
{\em saturated} if it is the inverse image of an open subset of $X/\hspace{-0.1cm}/H$. We will describe a saturated 
open neighborhood of $x \in X$. Identify $\mathbf{C}^2/\Gamma$ with 
$$\{(x_1, x_2, x_3) \in \mathbf{C}^3\: \vert \: x_1x_2 = x_3^m\}$$ and define a map $\nu$ by  
$$\nu: \mathbf{C}^2/\Gamma \to \mathbf{C}, \:\: (x_1, x_2, x_3) \to x_3.$$
Let $0 \in B_{\epsilon} \subset \mathbf{C}$ be a small open disc and put $W_{\epsilon} := \nu^{-1}(B_{\epsilon})$. 
Consider the affine space $\mathbf{C}^{2n-2}$ with coordinates $(t_1, ..., t_{2n-2})$ and take a sufficiently small 
disc $0 \in \Delta^{2n-2} \subset \mathbf{C}^{2n-2}$. Then $t \in \mathbf{C}^*$ acts on $W_{\epsilon} \times \Delta^{2n-2}$ by 
$$(x_1, x_2, x_3, t_1, ..., t_{2n-2}) \to (tx_1, t^{-1}x_2, x_3, t_1, ..., t_{2n-2}).$$
On the other hand, $H = \mathbf{C}^*$ acts on $X$.

\begin{Prop}\label{(2.10)}
There is a $\mathbf{C}^*$-equivariant open immersion $\Psi: W_{\epsilon} \times \Delta^{2n-2} \to X$ in such a way that $\Psi(0,0) = x$ and there is a 
Cartesian diagram
\begin{equation} 
\begin{CD} 
W_{\epsilon} \times \Delta^{2n-2} @>{\Psi \:\:\:(\subset )}>> X \\ 
@V{\nu \times id}VV  @V{\tau}VV  \\  
B_{\epsilon} \times \Delta^{2n-2} @>{\subset}>> X/\hspace{-0.1cm}/H.       
\end{CD} 
\end{equation}
\end{Prop}

{\em Proof}. 
As proved in Proposition \ref{(2.8)} there is a $\mathbf{C}^*$-equivariant isomorphism of complex analytic germs: 
$$((\mathbf{C}^2/\Gamma, 0) \times 
(\mathbf{C}^{2n-2}, 0), \:\: \omega_{\mathbf{C}^2/\Gamma} + \omega_{st}) \cong 
((X, x), \: \omega).$$ Then there 
exists an open neighborhood $0 \in U \subset \mathbf{C}^2/\Gamma$, an open neighborhood  
$x \in V \subset X$, and an isomorphism 
$$\Psi_{loc}: U \times \Delta^{2n-2} \stackrel{\cong}\to V$$ such that $\Psi_{loc}$ realizes the isomoprphism of the germs.  
Let us consider the map  $\nu: \mathbf{C}^2/\Gamma \to \mathbf{C}$.     
If we take $\epsilon$ sufficiently small,  then, for every $p \in W_{\epsilon}$, there is an element $t \in \mathbf{C}^*$ such that 
$t\cdot p \in U$. We then define 
$$\Psi (p; t_1, ..., t_{2n-2}) := t^{-1}\cdot \Psi_{loc}(t\cdot p; t_1, ..., t_{2n-1})$$ for 
$(p; t_1, ..., t_{2n-2}) \in W_{\epsilon} \times \Delta^{2n-2}$. Since $(t\cdot p; t_1, ..., t_{2n-1}) \in U \times \Delta^{2n-2}$, 
we have $\Psi_{loc}(t\cdot p; t_1, ..., t_{2n-1}) \in V$. Then $t^{-1} \in \mathbf{C}^*$ sends $\Psi_{loc}(t\cdot p; t_1, ..., t_{2n-1})$ to  
$t^{-1}\cdot \Psi_{loc}(t\cdot p; t_1, ..., t_{2n-1}) \in X$ by the $\mathbf{C}^*$-action on $X$. This $\Psi$ is a well-defined $\mathbf{C}^*$-equivariant map from $W_{\epsilon} \times \Delta^{2n-2}$ to $X$. In this situation we can apply an analytic version of Luna's Fundamental Lemma (cf. Theorem 1 of \cite{S}, Chapter 6, (1.2)) because, first  
$\Psi$ induces an isomorphism of neighborhoods of $(0; 0) \in W_{\epsilon} \times \Delta^{2n-2}$ and $x \in X$, next, 
$(0; 0) \in W_{\epsilon} \times \Delta^{2n-2}$ and $x \in X$ are both fixed points of the $\mathbf{C}^*$-action, and 
finally both $W_{\epsilon} \times \Delta^{2n-2}$ and $X$ have $\mathbf{C}^*$-linear embeddings in complex vector spaces. 
Then we see that $\Psi$ induces an isomorphism $W_{\epsilon} \times \Delta^{2n-2} \cong \Psi (W_{\epsilon} \times \Delta^{2n-2})$ and
$\Psi (W_{\epsilon} \times \Delta^{2n-2})$ is a saturated open subset of $X$.  $\square$
\vspace{0.2cm}

Recall that $Y$ is a symplectic leaf of $X$ of codimension $2$ with $x \in Y$. We choose the moment map $\mu: X \to (\mathfrak{t}^n)^*$ 
so that $\mu (x) = 0$. By definition, an open subset $U$ of $X$ is called $T^n$-saturated if $U$ is stable (as a set) by the 
$T^n$-action.  In the remainder we shall describe $\mu$ around a $T^n$-saturated open subset $U$ containing $x$.  
We write $T^n = H \times T^{n-1}$ with an $n-1$ dimensional subtorus $T^{n-1}$ of $T^n$ and let $(\theta_1, ..., \theta_{n-1})$ be the standard coordinates of $T^{n-1}$. Let $\Delta^{n-1}$ be a $n-1$ dimensional 
disc wth coordinates $(t_1, ..., t_{n-1})$. Let $H (= \mathbf{C}^*)$ act on $W_{\epsilon} \times \Delta^{n-1}$ by 
$$ (x_1, x_2, x_3, t_1, ..., t_{n-1}) \to (tx_1, t^{-1}x_2, x_3, t_1, ..., t_{n-1}), \:\:\: t \in \mathbf{C}^*.$$ 
Then $T^n = T^{n-1} \times H$ acts on $T^{n-1} \times (W_{\epsilon} \times \Delta^{n-1})$. Define a $T^n$-invariant 
symplectic form on $T^{n-1} \times (W_{\epsilon}^{reg} \times \Delta^{n-1})$ by 
$$\omega' := \mathrm{Res}(\frac{dx_1 \wedge dx_2 \wedge dx_3}{x_1x_2 - x_3^m}) +
dt_1 \wedge \frac{d\theta_1}{\theta_1} + ... + dt_{n-1} \wedge \frac{d\theta_{n-1}}{\theta_{n-1}}.$$ 

\begin{Thm}\label{(2.11)}
There are a $T^n$-saturated open subset $U$ of $X$ and a $T^n$- 
equivariant isomorphism of symplectic varieties 
$$\Phi: (T^{n-1} \times  (W_{\epsilon} \times \Delta^{n-1}), \omega')  \cong  (U, \omega\vert_U)$$ 
such that the moment maps commute    
\begin{equation} 
\begin{CD} 
T^n \times (W_{\epsilon} \times \Delta^{n-1}) @>{\Phi}>> U \\ 
@V{\mu'}VV  @V{\mu\vert_U}VV  \\  
(\mathfrak{t}^n)^* @>{id}>> (\mathfrak{t}^n)^*.    
\end{CD} 
\end{equation} 
The moment map $\mu'$ is given by 
$$\mu' (\theta_1, ..., \theta_{n-1}, x_1, x_2, x_3, t_1, ..., t_{n-1}) = (x_3, t_1, ..., t_{n-1}).$$
\end{Thm}

{\em Proof}. The orbit $T^nx$ is a smooth subvariety of $Y$ with dimension $n-1$. We have $$Y \cap (W_{\epsilon} \times \Delta^{2n-2}) =  \{0\} \times \Delta^{2n-2}.$$ 
As remarked just above,  we write $T^n = H \times T^{n-1}$ with an $n-1$-dimensional subtorus $T^{n-1}$.  
We may assume that the coordinates $(t_1, ..., t_{2n-2})$ 
of $\Delta^{2n-2}$ in Proposition \ref{(2.10)} are chosen such that $$T^nx \cap (W_{\epsilon} \times \Delta^{2n-2}) = \{0\} \times \{(0, ..., 0, t_n, ..., t_{2n-2}) \in \Delta^{2n-2}\},$$ 
where $t_n = \mathrm{log}\theta_1, ..., t_{2n-2} = \mathrm{log}\theta_{n-1}$. 
We write $$\Delta^{2n-2} = \Delta_1 \times \Delta_2$$ with $n-1$ dimensional discs  
$\Delta_1(t_1, ..., t_{n-1})$ and $\Delta_2(t_n, ..., t_{2n-2})$. Then    
$$T^nx \cap (W_{\epsilon} \times \Delta_1 \times \Delta_2) = \{0\} \times \{0\} \times \Delta_2.$$  
The group $H (= \mathbf{C}^*)$ acts on $W_{\epsilon} \times \Delta_1 \times \Delta_2$ by 
$$(x_1, x_2, x_3, t_1, ..., t_{2n-2}) \to (tx_1, t^{-1}x_2, x_3, t_1, ..., t_{2n-2}), \:\: t \in \mathbf{C}^*,$$ hence acts on 
$W_{\epsilon} \times \Delta_1 \times \{0\}$. 
Let us consider the $T^n$-variety $T^n \times^H (W_{\epsilon} \times \Delta_1 \times \{0\})$.  Then the inclusion 
$W_{\epsilon} \times \Delta_1 \times \{0\} \subset X$ induces a $T^n$-equivariant map 
$$\Phi: T^n \times^H (W_{\epsilon} \times \Delta_1 \times \{0\}) \to X.$$
The map $\Phi$ induces an isomorphism between a neighborhood of $[1, (0, 0, 0)] \in T^n \times^H (W_{\epsilon} \times \Delta_1 \times \{0\})$ and a neighborhood of $x \in X$. 

In fact, if we write $T^n = H \times T^{n-1}$ with an $n-1$-dimensional subtorus $T^{n-1}$ of 
$T^n$, then $$T^n \times^H (W_{\epsilon} \times \Delta_1 \times \{0\}) = T^{n-1} \times (W_{\epsilon} \times \Delta_1 \times \{0\}).$$
The germ $(T^{n-1} \times \{(0, 0, 0)\}, (1, (0,0,0)))$ is isomorphically mapped onto the germ $(\{0\} \times \{0\} \times \Delta_2, (0, 0, 0))$ by $\Phi$, and the germ $(\{1\} \times (W_{\epsilon} \times \Delta_1 \times \{0\}), (1, (0,0,0)))$ is isomorphically mapped onto 
the germ $(W_{\epsilon} \times \Delta_1 \times \{0\}, (0, 0, 0))$. Hence the tangential map 
$$d\Phi_{(1, (0,0,0))}: T_{(1, (0,0,0))}(T^n \times^H (W_{\epsilon} \times \Delta_1 \times \{0\})) \to T_xX$$ is an isomorphism. 
Let us take a linear $T^n$-embeddng of $X$ in some complex vector space $\mathbf{C}^N$. Then this means that the map 
$\Phi: T^n \times^H (W_{\epsilon} \times \Delta_1 \times \{0\}) \to \mathbf{C}^N$ is an embedding at $(1, (0,0,0))$. Since 
$$\dim T^n \times^H (W_{\epsilon} \times \Delta_1 \times \{0\}) = \dim X,$$ $\Phi$ induces an isomorphism between a 
neighborhood of $[1, (0, 0, 0)] \in T^n \times^H (W_{\epsilon} \times \Delta_1 \times \{0\})$ and a neighborhood of $x \in X$.   
Since the orbits $T^n[1, (0, 0, 0)]$ and $T^nx$ are both closed orbits with stabilizer group $H$, one can apply again an analytic version 
of Luna's fundamental lemma (cf. Theorem 1 of \cite{S}, Chapter 6, (1.2)). Then $\Phi$ is an open immersion and   
$U: = \Phi(T^n \times^H (W_{\epsilon} \times \Delta_1 \times \{0\}))$ is a saturated open subset of $X$ with respect to the $T^n$-action.

We note here that there is a $H$-equivariant open immersion  
$$\iota: W_{\epsilon} \times \Delta^{2n-2} \to T^{n-1} \times (W_{\epsilon} \times \Delta_1 \times \{0\})$$ 
given by $$(x_1, x_2, x_3, t_1, ..., t_{n-1}, t_n, ..., t_{2n-2}) \to (e^{t_n}, ..., e^{t_{2n-2}}, x_1, x_2, x_3, t_1, ..., t_{n-1}).$$
Let us consider the $T^n$-invariant symplectic form $$\omega' := \mathrm{Res}(\frac{dx_1 \wedge dx_2 \wedge dx_3}{x_1x_2 - x_3^m}) +
dt_1 \wedge \frac{d\theta_1}{\theta_1} + ... + dt_{n-1} \wedge \frac{d\theta_{n-1}}{\theta_{n-1}}$$ on $T^{n-1} \times (W_{\epsilon} \times \Delta_1 \times \{0\})$. 
Then we have $$\omega_{\mathbf{C}^2/\Gamma} + \omega_{st} = \iota^*\omega'.$$
The map $\Phi \circ \iota$ coincides with the open immersion $\Psi: W_{\epsilon} \times \Delta^{2n-2} 
\to X$ in Proposition \ref{(2.10)}. By Proposition \ref{(2.10)} we have $\Psi^*\omega = \omega_{\mathbf{C}^2/\Gamma} + \omega_{st}$. 
This implies that $\omega' = \Phi^*\omega$. Writing $\Delta^{n-1}$ for $\Delta_1$, we have the result.  $\square$ \vspace{0.2cm}

\section{}  Let $(X, \omega)$ be an affine symplectic variety of dimension $2n$ with an effective Hamiltonian action of an $n$-dimensional 
algebraic torus $T^n$.  In this section we impose an additional condition that $X$ is {\em conical},  that is, $(X, \omega)$ has a good 
$\mathbf{C}^*$-action, compatible with the $T^n$-action.  More precisely, the coordinate ring $R$ of $X$ is positively graded: $R = \oplus_{i \ge 0}R_i$ with $R_0 = \mathbf{C}$ and $\omega$ is homogeneous, that is, there is a positive integer $l$ such that $t^*\omega = t^l\omega$ for $t \in \mathbf{C}^*$. By definition, $X$ has a unique fixed point $0 \in X$ for the $\mathbf{C}^*$-action. This fixed point corresponds to 
the maximal ideal $\oplus_{i > 0}R_i$ of $R$. We take the moment map $\mu: X \to (\mathfrak{t}^n)^*$ in such a way that 
$\mu (0) = 0$.  In this section we first prove in Corollary \ref{(3.2)} that the moment map $\mu$ coincides with the GIT quotient of 
$X$ by $T^n$. In particular, $\mu$ is a $\mathbf{C}^*$-equivariant surjective map. The proof fully uses the conical $\mathbf{C}^*$-action 
on $X$. We next define the {\em discriminant divisor} $H \subset (\mathfrak{t}^n)^*$ of $\mu$ in \ref{(3.3)}, which turns out to be the union of a finte number of hyperplanes through $0 \in (\mathfrak{t}^n)^*$ with some weights. There is a closed subset $F_X \subset (\mathfrak{t}^n)^*$ 
with $\mathrm{Codim}_{(\mathfrak{t}^n)^*}F_X \geq 2$ such that, for $\eta \in (\mathfrak{t}^n)^* - F_X$, the fiber $\mu^{-1}(\eta)$ is singular  
exactly when $\eta$ is lying on $H$ (cf. \ref{(3.4)}). Write $H = m_1H_1+ ... + m_kH_k + H_{k+1} + ... + H_d$ with $m_i > 1$. Aside from $F_X$, the singularity of $X$ only appears over $H_1 \cup ... \cup H_k$.

\begin{Prop}\label{(3.1)} The moment map $\mu$ factors through $X/\hspace{-0.1cm}/T^n$: 
$$X \stackrel{\tau}\to X/\hspace{-0.1cm}/T^n \stackrel{\nu}\to (\mathfrak{t}^n)^*.$$ Both maps $\tau$ and $\nu$ are 
$\mathbf{C}^*$-equivariant. Here the $\mathbf{C}^*$-action on $(\mathfrak{t}^n)^*$ is given by the scaling action 
$\times t^l$ with $t \in \mathbf{C}^*$.   
\end{Prop}

{\em Proof}. We first prove that any $T^n$-orbit $O$ of $X$ is contained in a fiber of $\mu$. Such an orbit $O$ is contained in a symplectic leaf $Y$ of $X$. By Theorem \ref{(2.2)}, (1), we see that  $\mu\vert_Y : Y \to (\mathfrak{t}^n)^*$ factors through $\mathfrak{t}_Y^*$ and it coincides with the moment map for the $T^n$-action on $(Y, \omega_Y)$ (Note that Theorem (2.2), (1) holds true without the condition that $X$ has a symplectic resolution). By applying Lemma \ref{(1.1)} to $\mu_Y: Y \to \mathfrak{t}_Y^*$, we see that $O$ is contained in a fiber of $\mu_Y$; hence, $O$ is contained in a fiber of $\mu$. This fact means that $\mu$ factors through $X/\hspace{-0.1cm}/T^n$. 
 
Since the $T^n$-action and the $\mathbf{C}^*$-action commute, $\mathbf{C}^*$ acts on $X/\hspace{-0.1cm}/T^n$. We next prove that, for $t \in \mathbf{C}^*$, the following diagram commutes 
\begin{equation} 
\begin{CD} 
X @>{t}>> X \\ 
@V{\mu'}VV  @V{\mu}VV  \\  
(\mathfrak{t}^n)^* @>{\times t^l}>> (\mathfrak{t}^n)^*.
\end{CD} 
\end{equation} 
 
For a function $h$ on $X_{reg}$, we define a vector field $H_h$ on $X_{reg}$ so that $\omega (\cdot, H_h) = dh$.  
This correspondence determines a map $H: \Gamma (X_{reg}, \mathcal{O}_{X_{reg}}) \to \Gamma (X_{reg}, \Theta_{X_{reg}})$. 
The $T^n$-action on $X$ determines a map $\mathfrak{t}^n \to \Gamma (X_{reg}, \Theta_{X_{reg}})$.
By the definition of the moment map, this map is factorized as 
$$\mathfrak{t}^n \stackrel{\mu^*}\to \Gamma (X_{reg}, \mathcal{O}_{X_{reg}}) \stackrel{H}\to \Gamma (X_{reg}, \Theta_{X_{reg}}).$$ 
Take $f \in \mathfrak{t}^n$ and consider the vector field $H_{\mu^*f}$.  
Since the $T^n$-action and the $\mathbf{C}^*$-action commute, $H_{\mu^*f}$ is a $\mathbf{C}^*$-invariant vector field. In other words, 
we have $H_{\mu^*f} \in \Gamma (X_{reg}, \Theta_{X_{reg}})(0)$. We identfy $\Theta_{X_{reg}}$ with $\Omega^1_{X_{reg}}$ by $\omega$. 
Since $wt(\omega) = l$, we have $d(\mu^*f) \in \Gamma (X_{reg}, \Omega^1_{X_{reg}})(l)$.  
We have an exact sequence
$$0 \to \mathbf{C} \to \Gamma(X_{reg}, \mathcal{O}_{X_{reg}}) \stackrel{d}\to \Gamma (X_{reg}, \Omega^1_{X_{reg}}).$$      
The differential $d$ preserves the grading and induces a map $d^{(i)}: \Gamma(X_{reg}, \mathcal{O}_{X_{reg}})(i) \to 
\Gamma (X_{reg}, \Omega^1_{X_{reg}})(i)$ for each $i$. Then $\mathrm{Ker}(d^{(i)}) = 0$ for $i \ne 0$ and $\Gamma (X_{reg}, \mathcal{O}_{X_{reg}})(0) = \mathbf{C}$. It follows from these facts that $\mu^*f \in \mathbf{C} \oplus  \Gamma(X_{reg}, \mathcal{O}_{X_{reg}})(l)$. 
Since $f$ is a linear function on $(\mathfrak{t}^n)^*$, we have $\mu^*f(0) = 0$ and $$\mu^*f \in \Gamma (X_{reg}, \mathcal{O}_{X_{reg}})(l) 
= \Gamma (X, \mathcal{O}_X)(l).$$  
 Then the following diagram commutes 
\begin{equation} 
\begin{CD} 
\mathfrak{t}^n @>{\times t^l}>> \mathfrak{t}^n \\ 
@V{\mu^*}VV  @V{\mu^*}VV  \\  
\Gamma (X, \mathcal{O}_X)(l) @>{t^*}>> \Gamma (X, \mathcal{O}_X)(l).
\end{CD} 
\end{equation} 
In fact, we have $$t^*(\mu^*f) = t^l\mu^*f = \mu^*(t^l  f).$$ 
Therefore we have a commutative diagram 
\begin{equation} 
\begin{CD} 
\mathfrak{t}^n @>{\times t^l}>> \mathfrak{t}^n \\ 
@V{\mu^*}VV  @V{\mu^*}VV  \\  
\Gamma (X, \mathcal{O}_X) @>{t^*}>> \Gamma (X, \mathcal{O}_X).
\end{CD} 
\end{equation}  $\square$ 

\begin{Cor}\label{(3.2)}
Assume that $(X, \omega)$ has a projective symplectic resolution $\pi: (\tilde{X}, \omega_{\tilde X}) \to 
(X, \omega)$. Then $\nu: X/\hspace{-0.1cm}/T^n \to (\mathfrak{t}^n)^*$ is an isomorphism. In particular, the moment map $\mu$ is a surjection.
\end{Cor}

{\em Proof}. We first prove that $\nu$ is an etale map. The $T^n$-action on $(X, \omega)$ extends to a $T^n$-action on $(\tilde{X}, \omega_{\tilde X})$, which is a Hamiltonian action. Fix $x \in X$ and take $\tilde{x} \in \pi^{-1}(x)$.    
We take a $T^n$-invariant affine open subset $U$ of $\tilde{X}$ so that $\tilde{x} \in U$. The composite 
$$\mu_U: U \subset \tilde{X} \stackrel{\pi}\to X \to (\mathfrak{t}^n)^*$$ is a moment map for the $T^n$-action on $(U, \omega_{\tilde X}\vert_U)$. Then $\mu_U$ factors through $U/\hspace{-0.1cm}/T^n$: $$U \stackrel{{\tau}_U}\to U/\hspace{-0.1cm}/T^n \stackrel{\nu_U}\to (\mathfrak{t}^n)^*.$$ 
The map $\nu_U$ is etale by Theorem \ref{(1.3)}.  In particular, $$\nu_U^*: \hat{\mathcal{O}}_{(\mathfrak{t}^n)^*, \mu(x)} \to 
\hat{\mathcal{O}}_{U/\hspace{-0.1cm}/T^n, \tau_U(\tilde{x})}$$ is an isomorphism. Let us consider the commutative diagram:
\begin{equation} 
\begin{CD} 
U/\hspace{-0.1cm}/T^n @>{\nu_U}>> (\mathfrak{t}^n)^* \\ 
@VVV  @V{id}VV  \\  
X/\hspace{-0.1cm}/T^n @>{\nu}>> (\mathfrak{t}^n)^*.
\end{CD} 
\end{equation}

Then $\nu_U^*$ factors through $\hat{\mathcal{O}}_{X/\hspace{-0.1cm}/T^n, \tau(x)}$: 
$$\hat{\mathcal{O}}_{(\mathfrak{t}^n)^*, \mu(x)} \stackrel{\hat{\nu}^*}\to \hat{\mathcal{O}}_{X/\hspace{-0.1cm}/T^n, \tau(x)} \to   
\hat{\mathcal{O}}_{U/\hspace{-0.1cm}/T^n, \tau_U(\tilde{x})}.$$ Here the second map is an injection because $U/\hspace{-0.1cm}/T^n \to X/\hspace{-0.1cm}/T^n$ is a dominating map. 
Hence $\hat{\nu}^*$ is an isomorphism. This means that $\nu$ is an etale map. 

By Proposition \ref{(3.1)}, $\nu: X/\hspace{-0.1cm}/T^n \to (\mathfrak{t}^n)^*$ is $\mathbf{C}^*$-equivariant. We put $x = 0$, where $0$ is the origin of $X$. Note that the coordinate ring $\mathbf{C}[(\mathfrak{t}^n)^*]$ is the $\mathbf{C}$-subalgebra of $\hat{\mathcal{O}}_{(\mathfrak{t}^n)^*, 0}$ 
generated by $\mathbf{C}^*$-eigenvectors. Similarly,  the coordinate ring $\mathbf{C}[X/\hspace{-0.1cm}/T^n]$ is the $\mathbf{C}$-subalgebra of 
$\hat{\mathcal{O}}_{X/\hspace{-0.1cm}/T^n, \tau(0)}$ generated by $\mathbf{C}^*$-eigenvectors.. Then $\hat{\nu}^*$ induces the map $\nu^*:  \mathbf{C}[(\mathfrak{t}^n)^*]  \to  \mathbf{C}[X/\hspace{-0.1cm}/T^n]$, which is an isomorphism. Therefore $\nu$ is an isomorphism. $\square$  

\begin{Paragraph}\label{(3.3)} {\rm {\bf Hyperplane arrangements of $(\mathfrak{t}^n)^*$}. 

Let $Y_1, ..., Y_k$ be the symplectic leaves of $X$ of codimenson $2$.   
As above, we define 
$$H^{(i)} := \{t \in T^n\: \vert \: t \: \mathrm{acts}\: \mathrm{trivially} \: \mathrm{on}\: Y_i\}.$$ Then $H^{(i)}$ is a 1-dimensional (connected) subtorus of $T^n$ 
and $T_{Y_i} := T^n/H^{(i)}$ is an $(n-1)$-dimensional algebraic torus. 
Put $H_i := (\mathfrak{t}_{Y_i})^*$ for $1 \le i \le k$. Then $H_i$ is a hyperplane of 
$(\mathfrak{t}^n)^*$. We identify $\mathrm{Hom}_{alg.gp}(\mathbf{C}^*, T^n) \otimes_{\mathbf Z} \mathbf{C}$ with $\mathfrak{t}^n$.    
Then there is a primitive element $\mathbf{b}_i$ of $\mathrm{Hom}_{alg.gp}(\mathbf{C}^*, T^n) = \mathbf{Z}^n$ such that 
$$H_i = \{\eta \in (\mathfrak{t}^n)^*\: \vert \: \langle \mathbf{b}_i, \: \eta \rangle = 0\} \:\: (i = 1, ..., k)$$      
These hyperplanes are called of {\em the 1-st kind}.  
We next define hyperplanes of {\em the 2-nd kind}. Let 
$Y_{k+1}, ..., Y_r$ be the symplectic leaves of $X$ of codimension $\geq 4$. For these leaves we define similarly $T_{Y_i}$, which have 
dimension $\le n - 2$. We cover $(\mathfrak{t}^n)^* - \cup_{1 \le i \le r} (\mathfrak{t}_{Y_i})^*$ by a finite number of affine open subsets  
$V_j$ $(j \in J)$ and put $X_j := \mu^{-1}(V_j)$. Then $X_j$ is a $T^n$-invariant smooth affine open subset of $X$.   
Applying Theorem \ref{(1.3)}, (2) to $\mu\vert_{X_j}: X_j \to (\mathfrak{t}^n)^*$, we see that there are hyperplanes of 
$(\mathfrak{t}^n)^*$ such that $\mu\vert_{X_j}$ has singular fibers over these hyperplanes. Since these hyperplanes are stable under the 
$\mathbf{C}^*$-action on $(\mathfrak{t}^n)^*$, they all pass through $0 \in (\mathfrak{t}^n)^*$. We gather all such hyperplanes 
and form a set of hyperplanes, which we denote by $\{H_{k+1}, ..., H_d\}$. The hyperplanes $H_i$ 
$(k+1 \le i \le d)$ are those of the 2-nd kind. 
They are also defined by primitive vectors $\mathbf{b}_{k+1}, ..., \mathbf{b}_d \in \mathrm{Hom}_{alg.gp}(\mathbf{C}^*, T^n) = \mathbf{Z}^n$: 
$$H_i = \{\eta \in (\mathfrak{t}^n)^*\: \vert \: \langle \mathbf{b}_i, \: \eta \rangle = 0\} \:\: (i = k+1, ..., d).$$   
Assume that $X$ has Klein singularities of type $A_{m_i - 1}$ 
along $Y_i$ for each $1 \le i \le k$. Then the divisor 
$$m_1H_1 + ... + m_kH_k + H_{k+1} + ... + H_d$$ of $(\mathfrak{t}^n)^*$ is called the {\em discriminant divisor}. } 
\end{Paragraph}

\begin{Paragraph}\label{(3.4)} {\rm  {\bf The local structure of $\mu$}. 

We prove that there is a closed subset $F_X \subset (\mathfrak{t}^n)^*$ with $\mathrm{Codim}_{(\mathfrak{t}^n)^*}F_X \geq 2$,  which has 
the following properties. 

(i) For $\eta \in (\mathfrak{t}^n)^* - F_X$, the fiber $\mu^{-1}(\eta)$ is singular if and only if $\eta \in H_1 \cup ... \cup H_d$. 

(ii) Around each $\eta \in (\mathfrak{t}^n)^* - F_X$, the moment map $\mu$ has the normal form as given in  
Theorem \ref{(1.2)}, (2) or in Theorem \ref{(2.11)}.    
\vspace{0.2cm}

In fact, since $\mathrm{Sing}(X)$ is the union of the symplectic leaves of codimension $\geq 2$, the open subset $\mu^{-1}((\mathfrak{t}^n)^* - \cup_{1 \le i \le r} \mathfrak{t}_{Y_i}^*)$ of $X$ is smooth. For each $\eta \in (\mathfrak{t}^n)^* - \cup_{1 \le i \le r} \mathfrak{t}_{Y_i}^*$, 
we have a complete description of $\mu$ around $\eta$ by Theorem \ref{(1.3)}. Next consider a point $\eta \in H_i$ for $1 \le i \le k$. 
Then, by Theorem \ref{(2.11)}, there is a closed subset $F_i \subsetneq H_i$ such that $\mu$ has a normal form around each point 
$\eta \in H_i -F_i$. Put $$F_X := (\cup_{1 \le i \le k}F_i) \cup  (\cup_{k+1 \le i \le r}\mathfrak{t}_{Y_i}^*).$$ By definition we have $\mathrm{Codim}_{(\mathfrak{t}^n)^*}F_X \geq 2$. This $F_X$ has the desired properties. } 
\end{Paragraph}

\begin{Exam}\label{(3.5)} {\rm 
(toric hyperk\"{a}hler varieties)

Let $N$ and $n$ be positive integers such that $N \ge n$. 
Let $\mathbf{C}^{2N}$ be an affine space with coordinates $z_1$, ..., $z_N$, $w_1$, ..., $w_N$. 
An $N$ dimensional algebraic torus $T^N$ acts on $\mathbf{C}^{2N}$ by 
$$(z_1, ..., z_N, w_1, ..., w_N) \to (t_1z_1, ..., t_Nz_N, t_1^{-1}w_1, ..., t_N^{-1}w_N).$$
By an integer valued $(N-n) \times N$-matrix $A := (a_{ij})$, we determine a homomorphism of algebraic tori 
$\phi: T^{N-n} \to T^N$ by $$ (t_1, ..., t_{N-n}) \to (t_1^{a_{11}}\cdot\cdot\cdot t_{N-n}^{a_{N-n,1}}, ..., t_1^{a_{1,N}}\cdot\cdot\cdot t_{N-n}^{a_{N-n, N}}).$$ 
Then $T^{N-n}$ acts on $\mathbf{C}^{2N}$ by 
$$(z_1, ..., z_N, w_1, ..., w_N) \to$$ $$(t_1^{a_{11}}\cdot\cdot\cdot t_{N-n}^{a_{N-n,1}}z_1, ..., t_1^{a_{1,N}}\cdot\cdot\cdot t_{N-n}^{a_{N-n,n}}z_N, 
t_1^{-a_{11}}\cdot\cdot\cdot t_{N-n}^{-a_{N-n,1}}w_1, ..., t_1^{-a_{1N}}\cdot\cdot\cdot t_{N-n}^{-a_{N-n,n}}w_N).$$ 
The homomorphism $\phi$ induces a map of characters: $\phi^*: \mathrm{Hom}_{alg.gp}(T^N, \mathbf{C}^*) \to 
\mathrm{Hom}_{alg.gp}(T^N, \mathbf{C}^*)$. When we identify the character groups respectively with $\mathbf{Z}^N$ and 
$\mathbf{Z}^{N-n}$ in a natural way, $\phi^*$ is nothing but the homomorphism $\mathbf{Z}^N \stackrel{A}\to \mathbf{Z}^{N-n}$ 
determined by $A$. 

We assume that $A$ is surjective and unimodular, that is, any minor $(N-n) \times (N-n)$-matrix of $A$ has determinant $1$, $-1$ or $0$ and at least one of them has nonzero determinant. 
Let $B$ be an integer valued $N \times n$-matrix such that the following 
sequence is exact: 
$$0 \to \mathbf{Z}^n \stackrel{B}\to \mathbf{Z}^N \stackrel{A}\to \mathbf{Z}^{N-n} \to 0.$$  
Then $B$ is also unimodular. We assume that all row vectors of $B$ are nonzero.  If ncecessary, we change the coordinates 
$$(z_1, ..., z_N, w_1, ..., w_N) \to (z_{\sigma(1)}, ..., z_{\sigma(N)}, w_{\sigma(1)}, ..., w_{\sigma(N)}), \:\: \exists\sigma \in \mathfrak{S}_N$$ or  
$$(z_i, w_i) \to (-w_i, z_i)\:\:\: \exists i \in \{1, ..., N\},$$ so that the row vectors $\mathbf{b}_1$, ..., $\mathbf{b}_N$ of $B$ have the following properties:  There are integers $m_1> 1$, ..., $m_k > 1$ such that 
$$\mathbf{b}_1 = \cdot\cdot\cdot = \mathbf{b}_{m_1}, \:\: \mathbf{b}_{m_1+1} = \cdot\cdot\cdot = \mathbf{b}_{m_1 + m_2}, ... \:\: , \mathbf{b}_{m_1+ ... + m_{k-1} + 1} = \cdot\cdot\cdot =  
\mathbf{b}_{m_1 + ... + m_{k-1} + m_k}.$$ 
Moreover, $\mathbf{b}_{m_1}$, $\mathbf{b}_{m_1 + m_2}$, ..., 
$\mathbf{b}_{m_1+ ... + m_k}$, $\mathbf{b}_j$ ($m_1 + ... + m_k < j  \le N$) are mutually non-parallel vectors.   
 
Define a symplectic 2-form $\omega$ on $\mathbf{C}^{2N}$ by 
$$\omega_{\mathbf{C}^{2N}} := \sum_{1 \le i \le N} dw_i \wedge dz_i.$$  Then the $T^{N-n}$-action is a Hamiltonian action on $(\mathbf{C}^{2N}, \omega_{\mathbf{C}^{2N}})$. 
Writing $A = (\mathbf{a}_1, ..., \mathbf{a}_N)$ by the column vectors, the moment map $\mu: \mathbf{C}^{2N} \to \mathbf{C}^{N-n}$ is 
given by $$(z_1, ..., z_N, w_1, ..., w_N) \to \sum_{1 \le i \le N}\mathbf{a}_iz_iw_i.$$  

Note that $T^{N-n}$ acts on each fiber of $\mu$. 
Put $M := \mathrm{Hom}_{alg.gp}(T^{N-n}, \mathbf{C}^*)$. 
For $\alpha \in M$, we define $$Y(A, \alpha) := \mu^{-1}(0)
/\hspace{-0.1cm}/_{\alpha}T^{N-n}.$$ Note that $\dim Y(A, \alpha) = 2n$. 
The symplectic form $\omega_{\mathbf{C}^{2N}}$ on $\mathbf{C}^{2N}$ is reduced to a symplectic 2-form $\omega_{Y(A, \alpha)}$ on 
$Y(A, \alpha)_{reg}$. 
When $\alpha = 0$, $$Y(A, 0) = \mathrm{Spec}\: \mathbf{C}[\mu^{-1}(0)]^{T^{N-n}}.$$
On the other hand, if we take $\alpha$ general, then $Y(A, \alpha)$ is smooth and the map 
$$(Y(A, \alpha), \omega_{Y(A, \alpha)}) \to (Y(A, 0), \omega_{Y(A, 0)})$$ is a projective symplectic resolution. 
The affine variety $Y(A, 0)$ is a conical symplectic variety. In fact, the scaling $\mathbf{C}^*$-action on $\mathbf{C}^{2N}$ is 
restricted to a $\mathbf{C}^*$-action on $\mu^{-1}(0)$ and it descends to a conical action on $Y(A, 0)$. 
We made $Y(A, 0)$ by taking the quotient of $\mu^{-1}(0)$ by the subtorus $T^{N-n}$ of $T^N$. But the quotient torus 
$T^n := T^N/T^{N-n}$ still acts on $(Y(A, 0), \omega_{Y(A,0)})$, and it is a Hamiltonian action. Let 
$$\bar{\mu}: Y(A, 0) \to (\mathfrak{t}^n)^*$$ 
be the moment map with $\bar{\mu}(0) = 0$. $\square$} \end{Exam} 
 
In the situation of Example \ref{(3.5)}, define hyperplanes $H_i$ of $(\mathfrak{t}^n)^*$ by 
$$H_i = \{\eta \in (\mathfrak{t}^n)^*\: \vert \: \langle \mathbf{b}_i, \: \eta \rangle = 0\} \:\: (i = 1, ..., N).$$
By definition, there might possibly appear the same hyperplanes more than once.   

\begin{Prop}\label{(3.6)}   
The discriminant divisor of $\bar{\mu}$ is 
$$m_1H_{m_1} + m_2H_{m_1 + m_2} + ... + m_kH_{m_1 + ... +m_k} + \sum_{m_1+ ... + m_k < j  \le N}H_j .$$  
\end{Prop}

{\em Proof}. As in \cite{BD} one can view $Y(A, \alpha)$ as a hyperk\"{a}hler quotient of $\mathbf{H}^N$ by a compact torus 
$T^{N-n}_{\mathbf R} (= (S^1)^{N-n})$. The action of $(T_{\mathbf{R}})^{N - n}$ on $\mathbf{H}^N = \mathbf{C}^N \oplus (\mathbf{C}j)^N$ 
induces a hyperk\"{a}hler moment map 
$$\mu_{hk} := (\mu_I,\:\:  \mu_J + i \mu_K): \mathbf{H}^N \to (\mathfrak{t}^{N-n}_{\mathbf R})^* \times 
(\mathfrak{t}^{N-n})^*,$$ where $\mu_I(0, 0) = \mu_J(0, 0) = \mu_K(0, 0) = 0$.  The $\mu$ defined above coincides with 
$\mu_J + i\mu_K$. Regard $\alpha$ as an element of $(\mathfrak{t}^{N-n}_{\mathbf R})^*$. Then 
$$Y(A, \alpha) = \mu_{hk}^{-1}(\alpha, 0)/T_{\mathbf R}^{N-n}.$$ 
The action of $T^n_{\mathbf R} := T^N_{\mathbf R}/T^{N-n}_{\mathbf R}$ on $Y(A, \alpha)$ preserves the hyperk\"{a}hler structure on $Y(A, \alpha)$ and 
gives rise to a hyperk\"{a}hler moment map $$\bar{\mu}_{hk}:= (\bar{\mu}_I, \bar{\mu}_J+ i\bar{\mu}_K): Y(A, \alpha) \to (\mathfrak{t}^n_{\mathbf R})^* \times (\mathfrak{t}^n)^*,$$ where the $\bar{\mu}$ above coincides with $\bar{\mu}_J + i\bar{\mu}_K$. Take a lift 
$\tilde{\alpha} \in \mathrm{Hom}_{alg.gp}(T^N, \mathbf{C}^*) (= \mathbf{Z}^N)$ of $\alpha$ and write 
$$\tilde{\alpha} = \left(\begin{array}{cccc}
\tilde{\alpha}_1 \\ 
... \\ 
... \\ 
\tilde{\alpha}_N
\end{array}\right)$$
and define hyperplanes of $(\mathfrak{t}^n_{\mathbf R})^*$ 
$$H_j^{\tilde{\alpha}} := \{\eta \in (\mathfrak{t}^n_{\mathbf R})^*\: \vert \: \langle \mathbf{b}_j, \: \eta \rangle = \tilde{\alpha}_j\} \:\: (1 \le j  \le N).$$ 
By [\cite{BD}, Theorem 3.1], $\bar{\mu}_{hk}$ induces a homeomorphism $Y(A, \alpha)/T^n_{\mathbf R} \cong (\mathfrak{t}^n_{\mathbf R})^* \times (\mathfrak{t}^n)^*$ and, if $\eta \in (\mathfrak{t}^N_{\mathbf R})^* \times (\mathfrak{t}^n)^*$, then the $T_{\mathbf R}^n$-stabilizer of a point of $\bar{\mu}_{hk}^{-1}(\eta)$ is the torus whose Lie algebra is spanned by $\mathbf{b}_j$ for which $\eta \in H_j^{\tilde{\alpha}} \times H_j \subset  (\mathfrak{t}^n_{\mathbf R})^* \times (\mathfrak{t}^n)^*$. If we take $\alpha$ general, then we may assume that $H_j^{\tilde{\alpha}}$ are all different. Choose an integer $j_0$ so that  $$m_1 + \cdot\cdot\cdot + m_{i-1} + 1  \le j _0 \le m_1 + \cdot\cdot\cdot + m_{i-1} + m_i.$$
Take $\eta_{\mathbf C} \in H_{j_0}$ general so that $\eta_{\mathbf C} \notin H_j$ for any $H_j \ne H_{j_0}$.  Then $(\mathfrak{t}^n_{\mathbf R})^* \times \{\eta_{\mathbf C}\}$ intersects 
$\cup_{1 \le  j  \le N}(H_j^{\tilde{\alpha}} \times H_j )$ with  
$$(H^{\tilde{\alpha}}_{m_1 + \cdot\cdot\cdot + m_{i-1} + 1} \times \{\eta_{\mathbf C}\}) \:\:  \cup \cdot\cdot\cdot \cup \:\:  
(H^{\tilde{\alpha}}_{m_1 + \cdot\cdot\cdot + m_{i-1} + m_i} \times \{\eta_{\mathbf C}\}).$$
Let us consider the map  
$$\bar{\mu}_{\eta_{\mathbf C}} := \bar{\mu}_{hk}\vert_{\bar{\mu}^{-1}(\eta_{\mathbf C})}: \bar{\mu}^{-1}(\eta_{\mathbf C}) \to (\mathfrak{t}^n_{\mathbf R})^* \times \{\eta_{\mathbf C}\}.$$
Then a fiber of any point of  
$$(H^{\tilde{\alpha}}_{m_1 + \cdot\cdot\cdot + m_{i-1} + 1} \times \{\eta_{\mathbf C}\}) \:\:  \cup \cdot\cdot\cdot \cup \:\:  
(H^{\tilde{\alpha}}_{m_1 + \cdot\cdot\cdot + m_{i-1} + m_i} \times \{\eta_{\mathbf C}\})$$ 
is a $T^n_{\mathbf R}$-orbit with 1 dimensional stabilizer group and other fibers are all free $T^n_{\mathbf R}$-orbits.
This means that there is a symplectic leaf $Y$ of $Y(A, 0)$ of codimension 2 along which $Y(A, 0)$ has $A_{m_i - 1}$ singularities,  
such that  $(\mathfrak{t}_Y)^* = H_{j_0}$. Then the moment map for the $T^n$-action on $(Y(A, 0), \omega_{Y(A,0)})$ is locally described in Theorem \ref{(2.11)}. Let $f: \tilde{W}_{\epsilon} \to W_{\epsilon}$ be the minimal resolution. Let us consider the composite 
$x_3 \circ f: \tilde{W}_{\epsilon} \stackrel{f}\to W_{\epsilon} \stackrel{x_3}\to B_{\epsilon}$ and put $C := (x_3\circ f)^{-1}(0)$. $C$ consists  
of $m_i + 1$ irreducible components, $m_i - 1$ of which are exceptional divisors of $f$.     
The moment map for $(Y(A, \alpha), \omega_{Y(A, \alpha)})$ is locally written as  
$$T^n \times^H (\tilde{W}_{\epsilon} \times \Delta_1 \times \{0\}) \to (\mathfrak{t}^n)^*.$$  
Then $\bar{\mu}^{-1}(\eta_{\mathbf C})$ is isomorphic to $T^n \times^H(C \times \{0\} \times \{0\})$. Let $p_l$ $(l = 1, ..., m_i)$ be the double points of $C$. 
Then $T^n \times^H (\{p_l\} \times \{0\} \times \{0\})$ corresponds to 
$$(\bar{\mu}_{\eta_{\mathbf C}})^{-1}(H^{\tilde{\alpha}}_{m_1 + \cdot\cdot\cdot + m_{i-1} + l} \times \{\eta_{\mathbf C}\}).$$  
On the other hand, let $j_0$ be an integer such that $j_0 > m_1 + \cdot\cdot\cdot + m_k$ and take a general point of $\eta_{\mathbf C} 
\in H_{j_0}$. Then  
$(\mathfrak{t}^n_{\mathbf R})^* \times \{\eta_{\mathbf C}\}$ intersects 
$\cup_{1 \le  j  \le N}(H_j^{\tilde{\alpha}} \times H_j )$ with   
$H^{\tilde{\alpha}}_{j_0} \times \{\eta_{\mathbf C}\}$.  This means that $H_{j_0}$ is a hyperplane of the 2-nd kind. 
As a consequence, the discriminant divisor of $\bar{\mu}$ is 
$$m_1H_{m_1} + m_2H_{m_1 + m_2} + ... + m_kH_{m_1 + ... +m_k} + \sum_{m_1+ ... + m_k < j  \le N}H_j.$$
$\square$  \vspace{0.2cm}

\section{}  As in the previous section $(X, \omega)$ is a conical symplectic variety of dimension $2n$ with $wt(\omega) = l > 0$. 
We assume that $(X, \omega)$ admits a Hamiltonian $T^n$-action, compatible with the conical $\mathbf{C}^*$-action.   
Moreover, we assume that there is a projective symplectic resolution $\pi: (\tilde{X}, \omega_{\tilde X}) \to (X, \omega)$. 
The conical $\mathbf{C}^*$-action extends to a $\mathbf{C}^*$-action on $\tilde{X}$. 
The symplectic form $\omega_{\tilde X}$ determines a Poisson structure $\{\:, \:\}_{\tilde X}$ on $\tilde{X}$.
 
As explained in Introduction, we cannot directly prove that $X$ and $Y(A, 0)$ are isomorphic as $T^n$-Hamiltonian spaces over 
$(\mathfrak{t}^n)^*$. Instead we introduce a Poisson deformation $\mathcal{Z} \to \mathbf{C}^1$ of $X$ by using the universal Poisson deformation of a symplectic resolution $\tilde{X}$ of $X$, and construct its relative moment map $\mu_{\mathcal Z}: \mathcal{Z} \to (\mathfrak{t}^n)^* \times \mathbf{C}^1$. 
The main result of this section is Proposition \ref{(4.6)} which gives a local description of $\mu_{\mathcal Z}$. Proposition \ref{(4.6)} will be used in Proposition \ref{(5.1)}. 

When $X$ is smooth, $(X, \omega)$ is isomorphic to $(\mathbf{C}^{2n}, \omega_{st})$ as a $T^n$-space 
by Losev's result (see a footnote in Introduction). In the remainder we assume that $X$ is singular.
Put $r := b_2(\tilde{X})$. Since $\pi$ is projective, we have $r > 0$. 
Let $$f: (\tilde{\mathcal X}, \omega_{\tilde{\mathcal X}/\mathbf{C}^r}) \to \mathbf{C}^r$$ be the universal Poisson deformation. 
Here $\tilde{\mathcal X} \to \mathbf{C}^r$ is a smooth surjectve morphism whose central fiber $\tilde{\mathcal X}_0$ is identified with 
$\tilde{X}$ by an isomorphism $\phi: \tilde{X} \cong \tilde{\mathcal X}_0$, and 
$\omega_{\tilde{\mathcal X}/\mathbf{C}^r}$ is a relative symplectic form which determines a Poisson structure 
$\{\:, \:\}_{\tilde{\mathcal X}}$ over $\mathbf{C}^r$. The Poisson structure $\{\:, \:\}_{\tilde{\mathcal X}}$ is restricted to the original Poisson structure $\{\:, \:\}_{\tilde X}$ by $\phi$. There is a natural $\mathbf{C}^*$-action on $\tilde{\mathcal{X}}$ which is restricted to the 
$\mathbf{C}^*$-action on $\tilde{X}$ again by $\phi$. If we introduce a $\mathbf{C}^*$-action on $\mathbf{C}^r$ by the scaling  
$\times \sigma^l$ for $\sigma \in \mathbf{C}^*$, then the universal Poisson deformation is $\mathbf{C}^*$-equivariant 

The relative symplectic form $\omega_{\tilde{\mathcal X}/\mathbf{C}^r}$ determines the {\em period map}  
$p: \mathbf{C}^r \to H^2(\tilde{X}, \mathbf{C})$ as follows. We regard $\tilde{\mathcal X}$ and $\mathbf{C}^r$ as complex analytic spaces. 
Since the de Rham complex $\Omega^{\cdot}_{\tilde{\mathcal X}^{an}/\mathbf{C}^r}$ is a locally free resolution of 
$(f^{an})^{-1}\mathcal{O}_{\mathbf{C}^r}$, we have an isomorphism $$\mathbf{R}^2f_*\Omega^{\cdot}_{\tilde{\mathcal X}^{an}/\mathbf{C}^r} 
\cong R^2(f^{an})_*(f^{an})^{-1}\mathcal{O}_{\mathbf{C}^r}.$$ Since $f^{an}$ is a $C^{\infty}$-trivial fiber bundle with a typical fiber $\tilde{\mathcal{X}}^{an}_0$, 
we have an isomorphism $$R^2(f^{an})_*(f^{an})^{-1}\mathcal{O}_{\mathbf{C}^r} \cong H^2(\tilde{\mathcal X}_0, \mathbf{C}) \otimes_{\mathbf C}\mathcal{O}_{\mathbf{C}^r}.$$ For each $s \in \mathbf{C}^r$, we have the evaluation map 
$$ev_s: H^2(\tilde{\mathcal X}_0, \mathbf{C}) \otimes_{\mathbf C}\mathcal{O}_{\mathbf{C}^r} \to 
H^2(\tilde{\mathcal X}_0, \mathbf{C}) \otimes_{\mathbf C}k(s).$$ Composing these maps, we get 
$$\Gamma (ev_s): \Gamma (\mathbf{C}^r, \mathbf{R}^2(f^{an})_*\Omega^{\cdot}_{\tilde{\mathcal X}^{an}/\mathbf{C}^r}) \to 
H^2(\tilde{\mathcal X}_0, \mathbf{C}).$$ We regard $\omega_{\tilde{\mathcal X}/\mathbf{C}^r}$ as an element of 
$\Gamma (\mathbf{C}^r, \mathbf{R}^2(f^{an})_*\Omega^{\cdot}_{\tilde{\mathcal X}^{an}/\mathbf{C}^r})$. 
We are given an identification $\phi^*: H^2(\tilde{\mathcal X}_0, \mathbf{C}) \cong H^2(\tilde{X}, \mathbf{C})$. Then the period map is defined by $$p: \mathbf{C}^r \to H^2(\tilde{X}, \mathbf{C}), \:\:\: s \to \phi^* \circ \Gamma (ev_s)(\omega_{\tilde{\mathcal X}/\mathbf{C}^r}).$$
Then the period map $p$ turns out to be a $\mathbf{C}$-linear isomorphism. In fact, 
we introduce a $\mathbf{C}^*$-action on $H^2(\tilde{X}, \mathbf{C})$ by the scaling $\times \sigma^l$ for $\sigma \in \mathbf{C}^*$. 
Since $\omega_{\tilde{X}/\mathbf{C}^r}$ has weight $l$ for the $\mathbf{C}^*$-action on $\tilde{\mathcal X}$, we see that the period map $p$ is a $\mathbf{C}^*$-equivariant holomorphic map, which implies that $p$ is a $\mathbf{C}$-linear map. On the other hand, the 
tangential map $(dp)_0$ of $p$ at $0 \in \mathbf{C}^r$ coincides with the Poisson Kodaira-Spencer map $\kappa: T_0\mathbf{C}^r \to 
H^2(\tilde{X}, \mathbf{C})$, which is an isomoprphism by the universality of $f: (\tilde{\mathcal X}, \omega_{\tilde{\mathcal X}/\mathbf{C}^r}) \to \mathbf{C}^r$. Therefore $p$ is also an isomorphism.

\begin{Lem}\label{(4.1)}
The torus $T^n$ acts on $f: (\tilde{\mathcal X}, \omega_{\tilde{\mathcal X}/\mathbf{C}^r}) \to \mathbf{C}^r$ 
fiberwisely, that is, it acts trivially on the base $\mathbf{C}^r$.  
\end{Lem}

{\em Proof}. Since $f$ is the universal Poisson deformation of $\tilde{X}$, $T^n$ acts on $(\tilde{\mathcal X}, \omega_{\tilde{\mathcal X}/\mathbf{C}^r})$ and $\mathbf{C}^r$ in such a way that $f$ is $T^n$-equivariant. 
For $\lambda \in T^n$, we regard $\tilde{\mathcal X} \to \mathbf{C}^r$ as a Poisson deformation of $\tilde{X}$ by the identfication 
$\tilde{X} \stackrel{\lambda}\to \tilde{X} \stackrel{\phi}\to \tilde{\mathcal X}_0$. Let $p_{\lambda}: \mathbf{C}^r \to H^2(\tilde{X}^{an}, \mathbf{C})$ be the 
period map for this Poisson deformation. Since ${\lambda}^*: H^2(\tilde{X}, \mathbf{C}) \to H^2(\tilde{X}, \mathbf{C})$ is the identity map 
and ${\lambda}^*\omega_{\tilde{X}} = \omega_{\tilde{X}}$, we see that $p = p_{\lambda}$. This means that $T^n$ acts on the base $\mathbf{C}^r$ trivially. 
$\square$   \vspace{0.2cm}
  
Put $\mathcal{X} := \mathrm{Spec}\: \Gamma (\tilde{\mathcal X}, \mathcal{O}_{\tilde{\mathcal X}})$.
We then have a $\mathbf{C}^*$-equivariant Poisson deformation of $(X, \omega)$    
\begin{equation} 
\begin{CD} 
(X, \omega) @>>> ({\mathcal X}, \omega_{\mathcal{X}/\mathbf{C}^r}) \\ 
@VVV  @V{\bar{f}}VV  \\  
0 @>>> \mathbf{C}^r.
\end{CD} 
\end{equation}
The natural map $\tilde{\mathcal X} \stackrel{\Pi}\to \mathcal{X}$ induces a $\mathbf{C}^*$-equivariant commutative diagram of 
Poisson schemes
\begin{equation} 
\begin{CD} 
\tilde{\mathcal X} @>{\Pi}>> \mathcal{X} \\ 
@V{f}VV  @V{\bar{f}}VV  \\  
\mathbf{C}^r @>{id}>> \mathbf{C}^r.
\end{CD} 
\end{equation}
The map $\Pi$ is a birational projective morphism which induces birational morphisms of the fibers 
$\Pi_t: \tilde{\mathcal X}_t \to \mathcal{X}_t$, $t \in \mathbf{C}^r$. When $t = 0$, $\Pi_0 = \pi$ and 
when $t$ is general, $\Pi_t$ is an isomorphism. More precisely, 
there are a finite number of linear subspaces $\{L_i\}_{i \in I}$ of codimension $1$ in $\mathbf{C}^r$ such that, if 
$t \notin \cup L_i$, then $\Pi_t$ is an isomorphism (cf. \cite{Na 3}).  
Since $\mathcal{X}_t = \mathrm{Spec}\: \Gamma (\tilde{\mathcal X}_t, \mathcal{O}_{\tilde{\mathcal X}_t})$ for all $t$, 
the diagram above is $T^n$-equivariant.  In particular, we have 

\begin{Cor}\label{(4.2)} 
The torus $T^n$ acts on each fiber of the map $\bar{f}: \mathcal{X} \to \mathbf{C}^r$. $\square$
\end{Cor}

Moreover, we can prove: 

\begin{Prop}\label{(4.3)} 
The action of $T^n$ on each fiber of $\bar{f}$ is Hamiltonian, and there exists a relative moment map 
$$\mu_{{\mathcal X}/\mathbf{C}^r}: \mathcal{X} \to (\mathfrak{t}^n)^* \times \mathbf{C}^r,$$ which is a $\mathbf{C}^r$-morphism and 
$\mathbf{C}^*$-equivariant. Here the action of $\mathbf{C}^*$ on $(\mathfrak{t}^n)^* \times \mathbf{C}^r$ is the 
scaling action $\times \sigma^l$ for $\sigma \in \mathbf{C}^*$ with $l := wt(\omega)$. Moreover, $\mu_{\mathcal X}$ factorizes as  
$$\mathcal{X} \to \mathcal{X}/\hspace{-0.1cm}/T^n \stackrel{\nu_{\mathcal X}}\to (\mathfrak{t}^n)^* \times \mathbf{C}^r$$ and 
$\nu_{\mathcal X}$ is an isomorphism.  
\end{Prop}

{\em Proof}. Let $m \subset \mathcal{O}_{\mathbf{C}^r, 0}$ be the maximal ideal and put $A_k :=  \mathcal{O}_{\mathbf{C}^r, 0}/m^{k+1}$ and 
$S_k := \mathrm{Spec}\: A_k$. 
Set $\tilde{\mathcal X}_k := \tilde{\mathcal X} \times_{\mathbf{C}^r} S_k$. Consider the Lichenerowicz-Poisson complex  
$$\Theta^{\geq 1}_{\tilde{\mathcal X}_k/S_k} : \:\:\: \Theta_{\tilde{\mathcal X}_k/S_k} \stackrel{\delta_1}\to  
\wedge^2\Theta_{\tilde{\mathcal X}_k/S_k}  \stackrel{\delta_2}\to \wedge^3\Theta_{\tilde{\mathcal X}_k/S_k} \stackrel{\delta_3}\to \cdot\cdot\cdot. $$ 
Here we put $\wedge^i\Theta_{\tilde{\mathcal X}_k/S_k}$ on the degree $i$ part. Define $$P\Theta_{\tilde{\mathcal X}_k/S_k} := \mathrm{Ker}({\delta_1}).$$ Then it is easily checked that $$H^0(\tilde{X}, P\Theta_{\tilde{\mathcal X}_k/S_k}) \cong \mathbf{H}^1(\tilde{X}, 
\Theta^{\geq 1}_{\tilde{\mathcal X}_k/S_k}).$$ By $\omega_{\tilde{\mathcal X}/\mathbf{C}^r}$, we can identify  
$\Theta^{\geq 1}_{\tilde{\mathcal X}_k/S_k}$ with the truncated De Rham complex 
$$\Omega^{\geq 1}_{\tilde{\mathcal X}_k/S_k} : \:\:\: \Omega^1_{\tilde{\mathcal X}_k/S_k} \stackrel{d}\to 
\Omega^2_{\tilde{\mathcal X}_k/S_k}  \stackrel{d}\to \Omega^3_{\tilde{\mathcal X}_k/S_k} \stackrel{d}\to \cdot\cdot\cdot.$$  
If we put $$\Omega^{1, closed}_{\tilde{\mathcal X}_k/S_k} := \mathrm{Ker}(d) \subset \Omega^1_{\tilde{\mathcal X}_k/S_k},$$ then 
$$H^0(\tilde{X}, P\Theta_{\tilde{\mathcal X}_k/S_k}) \cong H^0(\tilde{X}, \Omega^{1, closed}_{\tilde{\mathcal X}_k/S_k}) \cong 
\mathbf{H}^1(\tilde{X}, \Omega^{\geq 1}_{\tilde{\mathcal X}_k/S_k}).$$
By Grothendieck's theorem, we have $$\mathbf{H}^i(\tilde{X}, \Omega^{\cdot}_{\tilde{\mathcal X}_k/S_k}) \cong H^i(\tilde{X}, A_k).$$ 
Now the exact triangle $$\Omega^{\geq 1}_{\tilde{\mathcal X}_k/S_k} \to \Omega^{\cdot}_{\tilde{\mathcal X}_k/S_k} \to 
\mathcal{O}_{\tilde{\mathcal X}_k} \to \Omega^{\geq 1}_{\tilde{\mathcal X}_k/S_k}[1]$$ yields 
an exact sequence $$\mathbf{H}^0(\tilde{X}, \Omega^{\cdot}_{\tilde{\mathcal X}_k/S_k}) \to H^0(\tilde{X}, \mathcal{O}_{\tilde{\mathcal X}_k}) \to \mathbf{H}^1(\tilde{X}, \Omega^{\geq 1}_{\tilde{\mathcal X}_k/S_k}) \to \mathbf{H}^1(\tilde{X}, \Omega^{\cdot}_{\tilde{\mathcal X}_k/S_k}).$$
The 1-st term is isomorphic to $A_k$ and the 4-th term vanishes because $H^1(\tilde{X}, A_k) = H^1(\tilde{X}, \mathbf{C})
\otimes_{\mathbf C}A_k = 0$. In particular, the map $$d: H^0(\tilde{X}, \mathcal{O}_{\tilde{\mathcal X}_k}) \to H^0(\tilde{X}, \Omega^{1, closed}_{\tilde{\mathcal X}_k/S_k})$$ is surjective. Define a map 
$$H_k:  H^0(\tilde{X}, \mathcal{O}_{\tilde{\mathcal X}_k}) \to H^0(\tilde{X}, P\Theta_{\tilde{\mathcal X}_k/S_k}), \:\:\: g \to H_g.$$  
Here $H_g$ is a Hamiltonian vector field uniquely defined by the property $\omega_{\tilde{\mathcal X}_k/S_k}(\cdot, H_g) = dg$. Note that 
$H_k$ is the composition of $d$ with the isomorphism $H^0(\tilde{X}, \Omega^{1, closed}_{\tilde{\mathcal X}_k/S_k}) \cong H^0(\tilde{X}, P\Theta_{\tilde{\mathcal X}_k/S_k})$. Therefore, $H_k$ is also surjective. Moreover, $\mathrm{Ker}(H_k) = A_k$.
Since $\omega_{\tilde{\mathcal X}/\mathbf{C}^r}$ has weight $l$, $H_k$ induces a surjection
$$H_k(l): H^0(\tilde{X}, \mathcal{O}_{\tilde{\mathcal X}_k})(l) \to H^0(\tilde{X}, P\Theta_{\tilde{\mathcal X}_k/S_k})(0).$$
Moreover, since $A_{k+1} \to A_k$ is a surjection for each $k$, we see that $$\mathrm{Ker}(H_{k+1}(l)) \to \mathrm{Ker}(H_k(l))$$ is also a surjection. 
On the other hand, the $T^n$-action on $\tilde{\mathcal X}$ determines a map 
$$\zeta: \mathfrak{t}^n \to H^0(\tilde{\mathcal X}, P\Theta_{\tilde{\mathcal X}/\mathbf{C}^r}).$$
Since the $T^n$-action commutes with the $\mathbf{C}^*$-action, we have 
$\mathrm{Im}(\zeta) \subset  H^0(\tilde{\mathcal X}, P\Theta_{\tilde{\mathcal X}/\mathbf{C}^r})(0)$. 
Let $e_1, ..., e_n$ be a basis of $\mathfrak{t}^n$ and we put $v_i := \zeta (e_i)$ for $1 \le i \le n$. These vector fields 
are restricted to vector fields $v_i^{(k)} \in H^0(\tilde{X}, P\Theta_{\tilde{\mathcal X}_k/S_k})(0)$. 
Note that $$\lim_{\leftarrow} H^0(\tilde{X}, \mathcal{O}_{\tilde{\mathcal X}_k})(l) \to   
\lim_{\leftarrow} H^0(\tilde{X}, P\Theta_{\tilde{\mathcal X}_k/S_k})(0)$$ is a surjection because, for each $k$,  
$H_k(l)$ is surjective and $\mathrm{Ker}(H_{k+1}(l)) \to \mathrm{Ker}(H_k(l))$ is surjective. Therefore we can find an eigen-element 
$$g_i: = \{g_i^{(k)}\} \in \lim_{\leftarrow} H^0(\tilde{X}, \mathcal{O}_{\tilde{\mathcal X}_k})(l)$$ for each $i$ such that 
$H_k(g_i^{(k)}) = v_i^{(k)}$ for all $k$. Note that $H^0(\tilde{X}, \mathcal{O}_{\tilde{\mathcal X}_k}) = H^0(X, \mathcal{O}_{\mathcal{X}_k})$. 
Since $g_i$ is an eigen-element, we have $g_i \in H^0(\mathcal{X}, \mathcal{O}_{\mathcal X})$. By the identification $H^0(\mathcal{X}, \mathcal{O}_{\mathcal X}) = 
H^0(\tilde{\mathcal{X}}, \mathcal{O}_{\tilde{\mathcal{X}}})$, we regard $g_i$ as an element of   
$H^0(\tilde{\mathcal{X}}, \mathcal{O}_{\tilde{\mathcal{X}}})(l)$. Now define a map 
$$\mu^*_{\tilde{\mathcal X}}: \:\: \mathfrak{t}^n \to H^0(\tilde{\mathcal{X}}, \mathcal{O}_{\tilde{\mathcal{X}}})(l)\:\:\:\:\:  \mathrm{by}\:\:\:\: 
\mu^*_{\tilde{\mathcal X}}(e_i) := g_i\:\: (i = 1, ..., n).$$ Then $\zeta$ factorizes as 
$$\zeta:\:\:  \mathfrak{t}^n  \stackrel{\mu^*_{\tilde{\mathcal X}}} \longrightarrow H^0(\tilde{\mathcal{X}}, \mathcal{O}_{\tilde{\mathcal{X}}})(l) \to 
H^0(\tilde{\mathcal X}, P\Theta_{\tilde{\mathcal X}/\mathbf{C}^r})(0).$$ 
Then $\mu^*_{\tilde{\mathcal X}}$ determines a moment map $$\mu_{\tilde{\mathcal X}}: \tilde{\mathcal X} \to (\mathfrak{t}^n)^*,$$ 
which factors through $\mathcal{X}$ and gives rise to a map $$\mu_{\mathcal X}:  \mathcal{X} \to (\mathfrak{t}^n)^*.$$
We define a relative moment map $$\mu_{\mathcal{X}/\mathbf{C}^r}: \mathcal{X} \to (\mathfrak{t}^n)^* \times \mathbf{C}^r$$ by 
$\mu_{\mathcal{X}/\mathbf{C}^r}  := \mu_{\mathcal{X}} \times \bar{f}$. The proof of the last statement is similar to the proof of Corollary 
\ref{(3.2)}.  $\square$ 

\begin{Rem}\label{(4.4)} {\rm
As is clear from the proof, the choice of $\mu^*_{\tilde{\mathcal X}}$ is not unique. 
Notice that  $H^0(\mathbf{C}^r, \mathcal{O}_{\mathbf{C}^r})(l)$ goes to zero by the map  $$H: H^0(\tilde{\mathcal{X}}, \mathcal{O}_{\tilde{\mathcal{X}}})(l) \to H^0(\tilde{\mathcal X}, P\Theta_{\tilde{\mathcal X}/\mathbf{C}^r})(0).$$ 
Note that $H^0(\mathbf{C}^r, \mathcal{O}_{\mathbf{C}^r})(l)$ is the space of linear functions on $\mathbf{C}^r$, namely, the dual space 
$(\mathbf{C}^r)^*$ of $\mathbf{C}^r$. Therefore, we can choose $\mu^*_{\tilde{\mathcal X}}$ up to an element of 
$$\mathrm{Hom}(\mathfrak{t}^n, (\mathbf{C}^r)^*) = \mathrm{Hom}(\mathbf{C}^r, (\mathfrak{t}^n)^*).$$} 
\end{Rem} 

\begin{Exam}\label{(4.5)} 
{\rm  Let $(Y(A, \alpha), \omega_{Y(A, \alpha)})$ be the same as in Example \ref{(3.5)}. If we take $\alpha$ general, then 
$(Y(A, \alpha), \omega_{Y(A, \alpha)}) \to (Y(A, 0), \omega_{Y(A, 0)})$ is a projective symplectic resolution. 
As $(X, \omega)$ we take $(Y(A, 0), \omega_{Y(A,0)})$, and as $(\tilde{X}, \omega_{\tilde X})$ we take $(Y(A, \alpha), \omega_{Y(A, \alpha)})$. 
We put $$X(A, \alpha) := \mathbf{C}^{2N}/\hspace{-0.1cm}/_{\alpha}T^{N-n}.$$  
The moment map $\mu: \mathbf{C}^{2N} \to (\mathfrak{t}^{N-n})^*$ factors through $X(A, \alpha)$ and gives rise to a map $X(A, \alpha) \to (\mathfrak{t}^{N-n})^*$. Since $Y(A, \alpha) = \mu^{-1}(0)/\hspace{-0.1cm}/_{\alpha}T^{N-n}$, 
$Y(A, \alpha)$ is nothing but the central fiber of this map. The symplectic 2-form $\omega_{\mathbf{C}^{2N}}$ descends to a relative symplectic 2-form $\omega_{X(A, \alpha)/(\mathfrak{t}^{N-n})^*}$ on the regular part $X(A, \alpha)_{reg}$ of $X(A, \alpha)$.
If we take $\alpha$ general, then $(X(A, \alpha), \omega_{X(A, \alpha)/(\mathfrak{t}^{N-n})^*}) \to (\mathfrak{t}^{N-n})^*$ is the universal Poisson deformation of $(Y(A, \alpha), \omega_{Y(A, \alpha)})$. Moreover, we have $X(A, 0) = \mathrm{Spec}\: \Gamma (X(A, \alpha), 
\mathcal{O}_{X(A, \alpha)})$ and the commutative diagram  
\begin{equation} 
\begin{CD} 
(Y(A, 0), \omega_{Y(A,0)}) @>>> (X(A, 0), \omega_{X(A,0)/(\mathfrak{t}^{N-n})^*}) \\ 
@VVV  @V{\bar{f}}VV  \\  
0 @>>> (\mathfrak{t}^{N-n})^* 
\end{CD} 
\end{equation}
corresponds to the commutative diagram 
\begin{equation}  
\begin{CD} 
(X, \omega) @>>> ({\mathcal X}, \omega_{\mathcal{X}/\mathbf{C}^r}) \\ 
@VVV  @V{\bar{f}}VV  \\  
0 @>>> \mathbf{C}^r
\end{CD} 
\end{equation} discussed above. 

The moment map $\mathbf{C}^{2N} \to (\mathfrak{t}^N)^* (= \mathbf{C}^{2N}/\hspace{-0.1cm}/_0T^N)$ factors through $X(A, 0) := \mathbf{C}^{2N}/\hspace{-0.1cm}/_0T^{N-n}$ and gives rise to a relative moment map 
$$\mu_{X(A, 0)/(\mathfrak{t}^{N-n})^*}: X(A, 0) \to (\mathfrak{t}^N)^*.$$
Here the surjection $(\mathfrak{t}^N)^* \stackrel{A}\to (\mathfrak{t}^{N-n})^*$ splits and one can write 
$(\mathfrak{t}^N)^* = (\mathfrak{t}^n)^* \times (\mathfrak{t}^{N-n})^*$. 
For $t \in (\mathfrak{t}^{N-n})^*$, we put $X(A, 0)_t := {\bar f}^{-1}(t)$ and $(\mathfrak{t}^N)^*_t := A^{-1}(t)$. 
Note that the relative moment map $\mu_{X(A, 0)/(\mathfrak{t}^{N-n})^*}$ induces a moment map 
$$\mu_{{X(A, 0)/(\mathfrak{t}^{N-n})^*}, t} : X(A, 0)_t  \to (\mathfrak{t}^N)^*_t$$ for each $t$. The discriminant divisor of 
$\mu_{{X(A, 0)/(\mathfrak{t}^{N-n})^*}, t}$ is explicitly described as follows. 
Let $e_1, ..., e_N \in \mathfrak{t}^N$ be the basis of $\mathfrak{t}^N$. We define divisors on $(\mathfrak{t}^N)^*$ by 
$$\mathcal{H}_i := \{\eta \in (\mathfrak{t}^N)^*\: \vert \: \langle e_i, \eta \rangle = 0\}, \:\:\: i = 1, 2, ..., N.$$
Then $$\sum_{1 \le i \le N} \mathcal{H}_{i, t} \subset (\mathfrak{t}^N)^*_t$$ is the discriminant divisor of  
$\mu_{{X(A, 0)/(\mathfrak{t}^{N-n})^*}, t}$. Note that, when $t = 0$, this is nothing but the discriminant divisor in Proposition \ref{(3.6)} 
because $e_i$ goes to $\mathbf{b_i}$ by the map $\mathfrak{t}^N \to \mathfrak{t}^n$. 
An important remark is that $\cap_{1 \le i \le N} \mathcal{H}_i = \{0\}$ because $e_1, ..., e_N$ are a basis of $\mathfrak{t}^N$. 
This means that $$\cap_{1 \le i \le N} \mathcal{H}_{i, t} = \emptyset \:\:\: \mathrm{for} \:\:\:  t \ne 0.$$
It is convenient to normalize the identification $(\mathfrak{t}^N)^* \cong (\mathfrak{t}^n)^* \times (\mathfrak{t}^{N-n})^*$ to 
see $\mathcal{H}_{i,t}$ in more details. 
To do so, we first take $\mathbf{b}_{i_1}$, ..., $\mathbf{b}_{i_n}$ so that they are a basis of $\mathfrak{t}^n$. We can take an 
identification in such a way that   
$$\mathcal{H}_{i_k, t} \cong \{\eta \in (\mathfrak{t}^n)^*\: \vert \: \langle \mathbf{b}_{i_k}, \eta \rangle = 0\} 
\times \{t\} \:\:\:\: k = 1, ..., n$$ 
for all $t \in (\mathfrak{t}^{N-n})^*$. 
For $1 \le i \le N$ with $i \notin \{i_1, ..., i_n\}$, we can write
$$\mathcal{H}_{i, t} \cong \{\eta \in (\mathfrak{t}^n)^*\: \vert \: \langle \mathbf{b}_i, \eta \rangle = f_i(t)\} \times \{t\}$$ 
with a linear function $f_i(t)$ on $(\mathfrak{t}^{N-n})^*$. We see that these $\{f_i(t)\}_{i \notin \{i_1, ..., i_n\}}$ are linearly independent. 
In fact, if they are linearly dependent, we find a nonzero $t_0 \in (\mathfrak{t}^{N-n})^*$ such that $f_i(t_0) = 0$ for all   
$i \notin \{i_1, ..., i_n\}$. Then $\cap_{1 \le i \le N}\mathcal{H}_{i, t_0} \ne \emptyset$, which is a contradiction. $\square$ } 
\end{Exam}

Consider the commutative diagram  
\begin{equation} 
\begin{CD} 
\tilde{\mathcal X} @>{\Pi}>> \mathcal{X} \\ 
@V{f}VV  @V{\bar{f}}VV  \\  
\mathbf{C}^r @>{id}>> \mathbf{C}^r.
\end{CD} 
\end{equation}
Recall that there are a finite number of linear subspaces $\{L_i\}_{i \in I}$ of codimension $1$ in $\mathbf{C}^r$ such that 
$\Pi_t$ is an isomorphism for $t \notin \cup L_i$ (\cite{Na 3}). Since $\tilde{\mathcal X}_t$ is smooth, $\mathcal{X}_t$ is also smooth 
for such $t$.  We take a line $\mathbf{C}^1 \subset \mathbf{C}^r$ passing through $0$ in such a way that $\mathbf{C}^1$ is not 
contained in any $L_i$, and pull back $\mathcal{X} \to \mathbf{C}^r$ to the line: 
\begin{equation} 
\begin{CD} 
\mathcal {Z} @>>> \mathcal{X} \\ 
@VVV  @VVV  \\  
\mathbf{C}^1 @>>> \mathbf{C}^r.
\end{CD} 
\end{equation}
Then the relative moment map $\mu_{\mathcal X}$ is restricted to the relative moment map over $\mathbf{C}^1$: 
$$\mu_{\mathcal Z}: \mathcal{Z} \to (\mathfrak{t}^n)^* \times \mathbf{C}^1.$$  
By the choice of the line $\mathbf{C}^1 \subset \mathbf{C}^r$, every fiber $\mathcal{Z}_t$ is a affine symplectic manifold 
with a Hamilton $T^n$-action for $t \in \mathbf{C}^1 - \{0\}$ and the moment map $\mu_{{\mathcal Z}_t}: \mathcal{Z}_t 
\to (\mathfrak{t}^n)^* \times \{t\}$ coincides with the quotient map $\mathcal{Z}_t \to \mathcal{Z}_t/\hspace{-0.1cm}/T^n$ by 
Proposition \ref{(4.3)}. Let us consider when $t = 1$. By Theorem \ref{(1.3)}, the discriminant divisor of  
$\mu_{{\mathcal Z}_1}$ is  $H_{1, 1} + \cdot\cdot\cdot  + H_{N, 1}$ in $(\mathfrak{t}^n)^* \times \{1\}$, where 
each $H_{i, 1}$ is written as 
$$H_{i, 1} := \{\eta \in (\mathfrak{t}^n)^*\: \vert \: \langle \mathbf{b}_i, \eta \rangle = \lambda_i\}$$ 
with a primitive vector $\mathbf{b}_i \in \mathbf{Z}^n = \mathrm{Hom}(\mathbf{C}^*, T^n) \subset \mathfrak{t}^n$ and 
$\lambda_i \in \mathbf{C}$. Moreover, they satisfy the following properties. 
\begin{itemize}
\item (P-1)  For any $n+1$ members $H_{i_1, 1}$, ..., $H_{i_n, 1}$ of them, we have $\cap_{1 \le k \le n} H_{i_k, 1} = \emptyset$.  
\item (P-2)  If $H_{i_1,1} \cap \cdot\cdot\cdot  \cap H_{i_m, 1} \ne \emptyset$, then $\{\mathbf{b}_{i_1}, ..., \mathbf{b}_{i_m}\}$ form a 
part of a basis of $\mathbf{Z}^n$.  
\end{itemize}
 
By the $\mathbf{C}^*$-action, the discriminant divisor of  
$\mu_{{\mathcal Z}_t}$ is  $H_{1, t} + \cdot\cdot\cdot  + H_{N, t}$ for $t \in \mathbf{C}^1 - \{0\}$, where 
$$H_{i, t} := \{\eta \in (\mathfrak{t}^n)^*\: \vert \: \langle \mathbf{b}_i, \eta \rangle = t\lambda_i\}.$$
In particular, all $H_{i, t}$ are different.  
For $t = 0$, we define 
$$H_{i, 0} := \{\eta \in (\mathfrak{t}^n)^*\: \vert \: \langle \mathbf{b}_i, \eta \rangle = 0\}.$$   

Let $H$ be an irreducible component of the discriminant divisor of $\mu (= \mu_{{\mathcal Z}_0})$. We shall study 
the local structure of the relative moment map $\mu_{\mathcal Z}: \mathcal{Z} \to (\mathfrak{t}^n)^* \times \mathbf{C}^1$ 
around $(\eta_0, 0) \in (\mathfrak{t}^n)^* \times \mathbf{C}^1$. 

For mutually different complex numbers $a_1$, ..., $a_m$, 
put $$\mathcal{Z}_{st} := \{(z_1, z_2, z_3, t_1, ..., t_{n-1}, t, \theta_1, ..., \theta_{n-1}) \in \mathbf{C}^2 \times \Delta^n \times \Delta^1 \times T^{n-1} \: \vert \: z_1z_2 = (z_3 + a_1t) \cdot \cdot \cdot (z_3 + a_mt) \}$$ and denote by $\mu_{\mathcal{Z}_{st}}: \mathcal{Z}_{st} \to \Delta^n \times \Delta^1$ the projection map $$(z_1, z_2, z_3, t_1, ..., t_{n-1}, t, \theta_1, ..., \theta_{n-1}) 
\to (z_3, t_1, ..., t_{n-1}, t)$$ and define $\mathcal{Z}_{st}^{reg} \subset \mathcal{Z}_{st}$ to be the locus where $\mu_{\mathcal{Z}_{st}}$ is smooth. Define a relative symplectic 2-form on $\mathcal{Z}_{st}^{reg}$ over $\Delta^1$ by 
$$\omega_{st} := \mathrm{Res}(\frac{dz_1 \wedge dz_2 \wedge dz_3}{z_1z_2 - \prod_{1 \le i \le m}(z_3 + a_it)}) + dt_1 \wedge \frac{d\theta_1}{\theta_1} + ... + dt_{n-1} \wedge \frac{d\theta_{n-1}}{\theta_{n-1}}.$$ 

\begin{Prop}\label{(4.6)}  
Let $H$ be an irreducible component of the discriminant divisor of $\mu (= \mu_{{\mathcal Z}_0})$. 
Take a general point $\eta_0$ on $H$ and a sufficiently small open neighborhood $(\eta_0, 0) \in \Delta^n \times\Delta^1 \subset 
(\mathfrak{t}^n)^* \times \Delta^1$.  Then, for a suitable choice of $a_1, ..., a_m$, there is a $T^n$-equivariant isomorphism $$(\mu^{-1}_{\mathcal Z}(\Delta^n \times \Delta^1), \omega_{\mathcal{Z}/\mathbf{C}^1}) \cong (\mathcal{Z}_{st}, \omega_{st})$$ such that the following diagram commutes
\begin{equation} 
\begin{CD} 
\mu^{-1}_{\mathcal Z}(\Delta^n \times \Delta^1) @>>> \mathcal{Z}_{st}\\ 
@V{\mu_{\mathcal Z}}VV  @V{\mu_{\mathcal{Z}_{st}}}VV  \\  
\Delta^{n} \times \Delta^1 @>{id}>> \Delta^{n} \times \Delta^1.
\end{CD} 
\end{equation}
\end{Prop}

{\em Proof}. Take a general point $\eta_0$ on $H$. Let $\eta_0 \in \Delta^n \subset (\mathfrak{t}^n)^*$ be a sufficiently  small 
disc. By Theorem \ref{(2.11)}  
$$\mu^{-1}(\Delta^n) = \{(x_1, x_2, x_3, \theta_1, ..., \theta_{n-1}, t_1, ..., t_{n-1} \in \mathbf{C}^3 \times (\mathbf{C}^*)^{n-1} \times 
\Delta^{n-1} \:\:\vert \:\:$$
$$ x_1x_2 = x_3^m,\:\: x_3 \in \Delta^1, \:\: (\theta_1, ..., \theta_{n-1}) \in (\mathbf{C}^*)^{n-1}, \:\: (t_1, ..., t_{n-1}) 
\in \Delta^{n-1}\}$$ and the map $\mu\vert_{\mu^{-1}(\Delta^n)}: \mu^{-1}(\Delta^n) \to \Delta^n$ is given by 
$$(x_1, x_2, x_3, \theta_1, ..., \theta_{n-1}, t_1, ..., t_{n-1}) \to (x_3, t_1, ..., t_{n-1}).$$   
$(\sigma, \sigma_1, ..., \sigma_{n-1}) \in T^n$ acts on $\mu^{-1}(\Delta^n)$ by 
$$(x_1, x_2, x_3, \theta_1, ..., \theta_{n-1}, t_1, ..., t_{n-1}) \to (\sigma x_1, \sigma^{-1}x_2, x_3, \sigma_1\theta_1, ..., \sigma_{n-1}\theta_{n-1}, 
t_1, ..., t_{n-1}).$$
Put $\mathbf{o} := (0,0,0; 1, ..., 1; 0, ..., 0) \in \mu^{-1}(\Delta^n)$. Then the stabilizer subgroup $T^n_{\mathbf o} \subset T^n$ for 
$\mathbf{o}$ is isomorphic 
to a 1-dimensional torus $T := \{\sigma, 1, ..., 1) \:\vert\: \sigma \in \mathbf{C}^*\}$. If we put 
$$V := \{(x_1, x_2, x_3, 1, ..., 1, t_1, ..., t_{n-1}) \in \mu^{-1}(\Delta^n)\},$$ then we have an identfication 
\begin{equation} 
\begin{CD} 
\mu^{-1}(\Delta^n) @>{\cong}>> T^n \times^T V\\ 
@VVV  @VVV  \\  
\Delta^{n} @>{\cong}>> V/\hspace{-0.1cm}/T.
\end{CD} 
\end{equation}
The map $V \to V/\hspace{-0.1cm}/T$ is given by $$(x_1 ,x_2 ,x_3, 1..., 1, t_1, ..., t_{n-1}) \to (x_3, t_1, ..., t_{n-1}),$$ hence all fibers are curves. 
The discrirminant locus of this map is $\{x_3 = 0\} \subset \Delta^n$. Namely, each fiber over the divisor is a nodal curve, but other fibers  
are smooth.    
 
On the other hand, let $(\eta_0, 0) \in \Delta^n \times \Delta^1 \subset (\mathfrak{t}^n)^* \times \mathbf{C}^1$ be a sufficiently small 
disc, and consider the $T^n$-variety $\mu_{\mathcal Z}^{-1}(\Delta^n \times \Delta^1)$. Applying an analytic version of Luna's slice theorem 
to the closed orbit $T^n\cdot \mathbf{o} \subset \mu_{\mathcal Z}^{-1}(\Delta^n \times \Delta^1)$, we can write 
$$\mu_{\mathcal Z}^{-1}(\Delta^n \times \Delta^1) = T^n \times^T \mathcal{V}$$ with a $T$-invariant subvariety $\mathcal{V} \subset 
\mu_{\mathcal Z}^{-1}(\Delta^n \times \Delta^1)$. Then we have an identification 
\begin{equation} 
\begin{CD} 
\mu_{\mathcal Z}^{-1}(\Delta^n \times \Delta^1) @>{\cong}>> T^n \times^T \mathcal{V}\\ 
@VVV  @VVV  \\  
\Delta^{n} \times \Delta^1 @>{\cong}>> \mathcal{V}/\hspace{-0.1cm}/T.
\end{CD} 
\end{equation}
Take the fiber product
\begin{equation} 
\begin{CD} 
\mathcal{V}_0 @>>> \mathcal{V}\\ 
@VVV  @VVV  \\  
\Delta^{n} \times \{0\} @>>> \Delta^{n} \times \Delta^1 = \mathcal{V}/\hspace{-0.1cm}/T.
\end{CD} 
\end{equation}
Then $\mu^{-1}(\Delta^n) = T^n \times^T \mathcal{V}_0$. This means that $\mathcal{V}_0 \cong V$ as a $T$-variety. 
In fact, put $T^{n-1} := \{(1, \theta_1, ..., \theta_{n-1} \in T^n \: \vert \: \theta_i \in \mathbf{C}^*\}$. Then 
$T^{n-1}$ acts on $\mu^{-1}(\Delta^n)$. The quotient $\mu^{-1}(\Delta^n)//T^{n-1}$ is naturally a $T$-variety. 
Since $\mu^{-1}(\Delta^n) \cong T^n \times^T V = T^{n-1} \times V$, we have $\mu^{-1}(\Delta^n)//T^{n-1} = V$. 
On the other hand, since $\mu^{-1}(\Delta^n) \cong T^n \times^T \mathcal{V}_0 = T^{n-1} \times \mathcal{V}_0$, we have 
$\mu^{-1}(\Delta^n)//T^{n-1} = \mathcal{V}_0$. Therefore $\mathcal{V}_0 \cong V$ as a $T$-variety.
Identify $\mathcal{V}_0$ with $V$ and consider the complex analytic germs at 
$\mathbf{o} := (0,0,0; 1, ..., 1; 0, ..., 0)$. We simply write $(x_1x_2 = 0, \mathbf{o})$ for the germ 
$$(\{ x_1x_2 = x_3 = 0,\:\:  \theta_1 = \cdot\cdot\cdot = \theta_{n-1} = 1, \:\: t_1= \cdot\cdot\cdot = t_{n-1} = 0\},\:\: \mathbf{o}),$$ which is the germ of a nodal curve. 

Then we have a commutative diagram 
\begin{equation} 
\begin{CD} 
(x_1x_2 = 0, \mathbf{o}) @>>> (V, \mathbf{o}) @>>> (\mathcal{V}, \mathbf{o})\\ 
@VVV @VVV  @VVV  \\  
(0,0) @>>> (\Delta^{n} \times \{0\}, (0,0)) @>>> (\Delta^{n} \times \Delta^1, (0,0)).
\end{CD} 
\end{equation}
They can be respectively regarded as $T$-equivariant flat deformations of the germ $(x_1x_2 = 0, \mathbf{o})$ of the nodal curve over the base spaces $(\Delta^{n} \times \{0\}, (0,0))$ and $(\Delta^{n} \times \Delta^1, (0,0))$.  
The $T$-equivariant semiuniversal deformation space for the germ $(x_1x_2 = 0, \mathbf{o})$ is given by 
$$\{(x_1, x_2, \lambda) \in (\mathbf{C}^3, 0)\: \vert \: x_1x_2 = \lambda\} \:\:\:\: \to \:\:\:\:  
\{\lambda \in (\mathbf{C}^1,0)\}.$$ 
The $T$ acts on the family by $$(x_1, x_2, \lambda) \to (\sigma x_1, \sigma^{-1}x_2, \lambda)$$ and 
$T$ acts trivially on the base space. 
The $T$-equivariant flat deformation $(\mathcal{V}, \mathbf{o}) \to (\Delta^n \times \Delta^1, (0,0))$ is obtained by pulling back 
the semiuniversal family by a map (cf. \cite{Pu})
$$\varphi: (\Delta^n \times \Delta^1, (0,0)) \to (\mathbf{C}^1, 0).$$ 
The relative moment map $\mu_{\mathcal Z}$ is the composite of the projection map $T^n \times^T \mathcal{V} = T^{n-1} \times \mathcal{V} \to \mathcal{V}$ and $\mathcal{V} \to \Delta^n \times \Delta^1$. 
In our case the discriminant divisor $\mathcal{H}$ of $\mu_{\mathcal Z}$ is given by 
$x_3^m + tg(x_3, t_1, ..., t_{n-1},t) = 0$ for some $g \in \mathcal{O}_{\Delta^n \times \Delta^1, 0}$. 
In particular, for each $t$, $H_t$ is a divisor of $\Delta^n \times \{t\}$. By our assumption, this means that 
$H = H_{i, 0} ( := \lim_{t \to 0} H_{i,t})$ for some $i$ and $\mathcal{H}$ has a form    
$$\prod_{1 \le i \le m} (x_3  + a_it) = 0,$$
with mutually different $a_1, ..., a_m \in \mathbf{C}$. In particular, there are exactly $m$ such $i$ that satisfy $H = H_{i,0}$.   
Therefore the map $\varphi$ must have the form 
$$\varphi (x_3, t_1, ..., t_{n-1}, t) = u(x_3, t_1, ..., t_{n-1}, t) (x_3 + a_1t) \cdot \cdot \cdot (x_3 + a_mt).$$
Here $u(x_3, t_1, ..., t_{n-1}, t)$ is a unit function around $0 := (0, 0, ..., 0, 0) \in \Delta^n \times \Delta^1$ such that 
$u(0, 0, ..., 0, 0) = 1$.  

Hereafter we put $z_1 := x_1u^{-1}$ and $z_2 := x_2$, and $z_3 := x_3$. Then we have a $T$-equivariant isomorphism 
$$(\mathcal{V}, \mathbf{o}) \cong \{(z_1, z_2, z_3, t_1, ..., t_{n-1}, t) \in (\mathbf{C}^{n+3}, 0)\: \vert \: 
z_1z_2 = (z_3 + a_1t) \cdot \cdot \cdot (z_3 + a_mt) \}$$ over $(\Delta^n \times \Delta^1, 0)$. Here 
the $T$-action on the right hand side is given by 
$$(z_1, z_2, z_3, t_1, ..., t_{n-1}, t) \to (\sigma z_1, \sigma^{-1}z_2, z_3, t_1, ..., t_{n-1}, t).$$
This isomorphism induces a $T^n$-equivariant isomorphism 
$$T^n \times^T (\mathcal{V}, \mathbf{0}) \cong \{(z_1, z_2, z_3, t_1, ..., t_{n-1}, t, \theta_1, ..., \theta_{n-1}) \in (\mathbf{C}^{n+3}, 0) \times T^{n-1} \: \vert \: 
z_1z_2 = (z_3 + a_1t) \cdot \cdot \cdot (z_3 + a_mt) \}$$ over $(\Delta^n \times \Delta^1, 0)$. 
By the $T^n$-action, this isomorpshim extends to a $T^n$-equivariant isomorphism 
$$T^n \times^T \mathcal{V} \cong 
\{(z_1, z_2, z_3, t_1, ..., t_{n-1}, t, \theta_1, ..., \theta_{n-1}) \in \mathbf{C}^2 \times \Delta^n \times \Delta^1 \times T^{n-1} \: \vert \: 
z_1z_2 = (z_3 + a_1t) \cdot \cdot \cdot (z_3 + a_mt) \}$$ over $\Delta^n \times \Delta^1$. 
By this isomorphism we regard $\omega_{{\mathcal Z}/\mathbf{C}^1}\vert_{T^n \times^T \mathcal{V}}$ as a 
relative symplectic 2-form on the right hand side. 
We write $\mu_{\mathcal Z}\vert_{T^n \times^T \mathcal{V}}$ for the restriction of the relative moment map $\mu_{\mathcal Z}$ to $T^n \times^T \mathcal{V}$. 
Then, under this isomorphism, the relative moment map  
$$\mu_{\mathcal Z}\vert_{T^n \times^T \mathcal{V}}\: : \: T^n \times^T \mathcal{V} =  T^{n-1} \times \mathcal{V} \to \Delta^n \times \Delta^1$$ is given by 
$$(\theta_1, ..., \theta_{n-1}, z_1, z_2, z_3, t_1, ..., t_{n-1}, t) \to (z_3, t_1, ..., t_{n-1}, t).$$ 

Let us consider the relative symplectic 2-form on  $T^n \times^T \mathcal{V}$: 
$$\omega_{st} := \mathrm{Res}(\frac{dz_1 \wedge dz_2 \wedge dz_3}{z_1z_2 - \prod_{1 \le i \le m}(z_3 + a_it)}) + dt_1 \wedge \frac{d\theta_1}{\theta_1} + ... + dt_{n-1} \wedge \frac{d\theta_{n-1}}{\theta_{n-1}}.$$
It is easily checked that $\mu_{\mathcal Z}\vert_{T^n \times^T \mathcal{V}}$ is the relative moment map 
for this relative symplectic 2-form.
We compare this 2-form with $\omega_{\mathcal{Z}/\mathbf{C}^1}\vert_{T^n \times^T \mathcal{V}}$. 
Since both $T^n$-invariant relative 2-forms have the same relative moment map, we can write 
$$\omega_{\mathcal{Z}/\mathbf{C}^1}\vert_{T^n \times^T \mathcal{V}} = 
\omega_{st} + (\mu_{\mathcal Z}\vert_{T^n \times^T \mathcal{V}})^*\eta$$ with a suitable $d$-closed relative 2-form  $\eta \in \Omega^2_{\Delta^n \times \Delta^1/\Delta^1}$. 

We shall prove that, after shrinking $\Delta^n$ and $\Delta^1$ further, there is a $T^n$-equivariant automorphism $\phi$ of $T^n \times^T \mathcal{V}$ over $\Delta^n \times \Delta^1$ such that $\phi \vert_{\mu_{\mathcal Z}^{-1}(0)} = id$ and  
$$\phi^*(\omega_{\mathcal{Z}/\mathbf{C}^1}\vert_{T^n \times^T \mathcal{V}}) 
= \omega_{st}.$$
Recall that the map $\mathcal{V} \to \Delta^n \times \Delta^1$ 
has a simultaneous crepant resolution $\tilde{\mathcal{V}} \to \mathcal{V}$, which induces a simultaneous resolution 
$\nu: T^n \times^T \tilde{\mathcal V} \to T^n \times^T \mathcal{V}$ of $T^n \times^T \mathcal{V} \to \Delta^n \times \Delta^1$.
We compare $\nu^*(\omega_{\mathcal{Z}/\mathbf{C}^1}\vert_{T^n \times^T \mathcal{V}})$ and $\nu^*\omega_{st}$. 

Under the identification  
$$T^n \times^T \mathcal{V} \cong 
\{(z_1, z_2, z_3, t_1, ..., t_{n-1}, t, \theta_1, ..., \theta_{n-1}) \in \mathbf{C}^2 \times \Delta^n \times \Delta^1 \times T^{n-1} \: \vert \: 
z_1z_2 = (z_3 + a_1t) \cdot \cdot \cdot (z_3 + a_mt) \}$$
we take a point $\mathbf{p} \in  T^n \times^T\mathcal{V}$ defined by 
$$z_1 = z_2 = z_3 = t_1 = ... = t_{n-1} = t = 0,\:\: \theta_1 = ... = \theta_{n-1} = 1.$$ 
Note that $\nu^{-1}(\mathbf{p})$ is a tree of $m-1$ smooth rational curves. Let us choose a nodal point $\mathbf{q} 
\in \nu^{-1}(\mathbf{p})$ and we regard $\mathbf{q}$ as a point of $T^n \times^T \tilde{\mathcal V}$.
 
We use the standard argument of Moser's proof of Darboux theorem. 
In order to do this, we can start with the situation where $\nu^*(\omega_{\mathcal{Z}/\mathbf{C}^1}\vert_{T^n \times^T \mathcal{V}})(\mathbf{q}) = \nu^*\omega_{st}(\mathbf{q})$.    
In fact,  we write $$\eta = \sum_{1 \le i \le n-1} f_idz_3 \wedge dt_i + \sum_{1 \le i < j \le n-1}g_{ij}dt_i \wedge dt_j$$ 
with functions $f_i$, $g_{ij}$ on $\Delta^n \times \Delta^1$. Put $b_{ij} := g_{ij}(0)$ for the origin $0 \in \Delta^n \times \Delta^1$. 
Then we may assume that $g_{ij}(0) = 0$ for all $i < j$ by taking   
the $T^n$-equivariant automorphism $\phi'$ of $T^n \times^T \mathcal{V}$ over $\Delta^n \times \Delta^1$ defined by 
$$ z_1 \to z_1, z_2 \to z_2,  z_3 \to z_3,\:\: t_i \to t_i (1 \le i \le n-1)$$ 
$$\theta_1 \to e^{b_{12}t_2 + \cdot\cdot\cdot  + b_{1, n-1}t_{n-1}}\theta_1,\:\:  ..., \theta_{n-2} \to e^{b_{n-2, n-1}t_{n-1}}\theta_{n-2}, \:\:
\theta_{n-1} \to \theta_{n-1}.$$
Since $(\mu_{\mathcal Z}\vert_{T^n \times^T \mathcal{V}} \circ \nu)^*(dz_3)(\mathbf{q}) = 0$, we see that 
$$(\mu_{\mathcal Z}\vert_{T^n \times^T \mathcal{V}} \circ \nu)^*(dz_3 \wedge dt_i)(\mathbf{q}) = 0, \:\: 1 \le i \le n-1.$$
Therefore $(\mu_{\mathcal Z}\vert_{T^n \times^T \mathcal{V}} \circ \nu)^*\eta (\mathbf{q}) = 0$ and hence,  
$\nu^*(\omega_{\mathcal{Z}/\mathbf{C}^1}\vert_{T^n \times^T \mathcal{V}})(\mathbf{q}) = \nu^*\omega_{st}(\mathbf{q})$. 
\vspace{0.3cm}

For simplicity we put $\omega_1 := \nu^*(\omega_{\mathcal{Z}/\mathbf{C}^1}\vert_{T^n \times^T \mathcal{V}})$ and 
$\omega_2 := \nu^*\omega_{st}$. Define  
$\omega (\lambda) := \lambda \omega_1 + (1 - \lambda)\omega_2$ for $\lambda \in [0,1]$. Then $\omega (\lambda) = 
\omega_2 + \lambda (\mu_{\mathcal Z}\vert_{T^n \times^T \mathcal{V}} \circ \nu)^*\eta$, hence  
$$\frac{d\omega(\lambda)}{d\lambda} =  (\mu_{\mathcal Z}\vert_{T^n \times^T \mathcal{V}} \circ \nu)^*\eta.$$ 
Since $\eta$ is a $d$-closed relative 2-form on $\Delta^n \times \Delta^1/\Delta^1$, one can write $\eta = d\gamma$ with a 
relative 1-form $\gamma$. We may assume that $\gamma (0) = 0$. Define a vector field $X_{\lambda} \in \Theta_{T^n \times^T \tilde{\mathcal{V}}/\Delta^1}$ by 
$$X_{\lambda} \rfloor \omega(\lambda) = - (\mu_{\mathcal Z}\vert_{T^n \times^T \mathcal{V}} \circ \nu)^*\gamma.$$
Then we have 
$$L_{X_{\lambda}}\omega(\lambda) = d({X_{\lambda}}\rfloor \omega(\lambda)) + {X_{\lambda}}\rfloor d\omega(\lambda) = d({X_{\lambda}}\rfloor \omega(\lambda)) = - (\mu_{\mathcal Z}\vert_{T^n \times^T \mathcal{V}} \circ \nu)^*\eta.$$
One can check that $X_0 \in \Theta_{T^n \times^T \tilde{\mathcal{V}}/\Delta^n \times \Delta^1}$ by a direct calculation. 
On the other hand, since $(\mu_{\mathcal Z}\vert_{T^n \times^T \mathcal{V}} \circ \nu)_*X_0 = 0$, we have 
$$X_0 \rfloor (\mu_{\mathcal Z}\vert_{T^n \times^T \mathcal{V}} \circ \nu)^*\eta = 0,$$ which implies that $X_{\lambda} = X_0$ 
for all $\lambda$. 
Moreover, since $\gamma (0) = 0$, we see that $X_{\lambda}$ vanishes along $(\mu_{\mathcal Z}\vert_{T^n \times^T \mathcal{V}} \circ \nu)^{-1}(0)$. 
Note that $X_{\lambda}$ is $T^n$-invariant because $\omega (\lambda)$ and $(\mu_{\mathcal Z}\vert_{T^n \times^T \mathcal{V}} \circ \nu)^*\gamma$ are both $T^n$-invariant. If necessary, shrinking $\Delta^n \times \Delta^1$ around the origin, the vector field $X_{\lambda}$ then defines a family of $T^n$-equivariant automorphisms 
$\phi_{\lambda}$ $(0 \le \lambda \le 1)$ of $T^n \times^T \tilde{\mathcal{V}}$ over $\Delta^n \times \Delta^1$ with $\phi_0 = id$. 

Since $$\frac{d}{d{\lambda}}\phi_{\lambda}^*\omega(\lambda) = \phi_{\lambda}^*(L_{X_{\lambda}}\omega(\lambda) + \frac{d\omega(\lambda)}{d\lambda})$$ 
$$= \phi_{\lambda}^*\{-(\mu_{\mathcal Z}\vert_{T^n \times^T \mathcal{V}} \circ \nu)^*\eta) + (\mu_{\mathcal Z}\vert_{T^n \times^T \mathcal{V}} \circ \nu)^*\eta \}  = 0,$$ we have 
$\phi_{\lambda}^*\omega(\lambda) = \omega_2$. In particular, 
when $\lambda = 1$, we have  $\omega_2 = \phi_1^*\omega_1$. The automorphism $\phi_1$ descends to a $T^n$-equivariant automorphism $\phi$ of $T^n \times^T \mathcal{V}$ over $\Delta^n \times \Delta^1$ with the desired property.  $\square$ 
\vspace{0.2cm}


By the proof of Proposition \ref{(4.6)}, we have:  

\begin{Cor}\label{(4.7)}
The discriminant divisor of $\mu (= \mu_{\mathcal{Z}_0})$ is 
$H_{1, 0} + \cdot\cdot\cdot + H_{N, 0}$. 
\end{Cor} 

{\em Remark}. It may possibly occur that $H_{i, 0} = H_{i', 0}$ even if $i \ne i'$.

\section{}  In this section we associate $X$ with a toric hyperk\"{a}hler variety $Y(A, 0)$ with $A$ unimodular. 
In the previous section we have constructed a Poisson deformation $\mathcal{Z} \to \mathbf{C}^1$ of $X$ and its relative moment map $\mu_{\mathcal Z}: \mathcal{Z} \to (\mathfrak{t}^n)^* \times \mathbf{C}^1$. 
In this section we shall construct a Poisson deformation $\mathcal{Z}' \to \mathbf{C}^1$ of $Y(A,0)$ and its relative moment map 
$\mu_{\mathcal{Z}'}: \mathcal{Z}' \to (\mathfrak{t}^n)^* \times \mathbf{C}^1$ in such a way that the discriminant divisors of $\mu_{\mathcal Z}$ and $\mu_{\mathcal{Z}'}$ coincide. Finally in Theorem \ref{(5.5)} we prove that $\mathcal{Z}$ is isomorphic to $\mathcal{Z}'$ as a $T^n$-Hamiltonian space over $(\mathfrak{t}^n)^* \times \mathbf{C}^1$. The main theorem is immediately obtained from Theorem \ref{(5.5)}.   
  
Let $(X, \omega)$ and $\mu_{\mathcal Z}: \mathcal{Z} \to (\mathfrak{t}^n)^* \times \mathbf{C}^1$ be the same as in \S 4. As explained in \S 4, each irreducible component of the discriminant divisor 
of $\mu_{\mathcal Z}$ is described in terms of  a primitive vector $\mathbf{b}_i$ in 
$\mathrm{Hom}(\mathbf{C}^*, T^n) \subset \mathfrak{t}^n$. 
The vectors $\mathbf{b}_1$, ..., $\mathbf{b}_N$ 
determine a homomorphism $B: \mathbf{Z}^n \to \mathbf{Z}^N$. We assume that 
$$(*) \:\:\: B\: \mathrm{is}\: \mathrm{an} \: \mathrm{injection}\:  \mathrm{and}\:  n < N.$$
By this assumption one can choose $n$ linearly independent vectors $\mathbf{b}_{i_1}, ..., \mathbf{b}_{i_n}$. Recall that we have remarked before Proposition (4.6) that the discriminant divisors $H_{i, 1}$ of $\mu_{\mathcal{Z}_1}$ have the properties (P-1) and (P-2).  Since $\mathbf{b}_{i_1}, ..., \mathbf{b}_{i_n}$ are linearly independent,  we have $H_{i_1, 1} \cap ... \cap H_{i_n, 1} \ne \emptyset$. Then by (P-2) we see that    
$\mathbf{b}_{i_1}, ..., \mathbf{b}_{i_n}$ form a basis of $\mathrm{Hom}(\mathbf{C}^*, T^n)$.  In particular, $B$ is unimodular.    
As we remarked in Remark \ref{(4.4)} we normalize the relative moment map $\mu_{\mathcal Z}$ in such a way 
that $\lambda_{i_1} = \cdot\cdot\cdot = \lambda_{i_n} = 0$. We define a divisor of $(\mathfrak{t}^n)^* \times \mathbf{C}^1$ by 
$$\mathcal{H} := \bigcup_{1 \le i \le N, \: t \in \mathbf{C}^1}H_{i, t}.$$   
Since $B$ is unimodular, there is an exact sequence 
$$0 \to \mathbf{Z}^n \stackrel{B}\to \mathbf{Z}^N \stackrel{A}\to \mathbf{Z}^{N-n} \to 0.$$
Here $A$ is also unimodular. For this $A$, we define a toric hyperk\"{a}hler variety $Y(A, 0)$ and its Poisson deformation 
$X(A, 0)$
\begin{equation} 
\begin{CD} 
(Y(A, 0), \omega_{Y(A,0)}) @>>>  (X(A, 0), \omega_{X(A,0)/\mathbf{C}^{N-n}})\\ 
@VVV @VVV  \\  
\{0\} @>>> \mathbf{C}^{N-n}.
\end{CD} 
\end{equation}

As in Example \ref{(4.5)} we choose a relative moment map  $\mu_{X(A, 0)/(\mathfrak{t}^{N-n})^*}: X(A, 0) \to (\mathfrak{t}^n)^* \times (\mathfrak{t}^{N-n})^*$  in such a way that the  discriminant divisor $\mathcal{H}^{X(A,0)} \subset (\mathfrak{t}^n)^* \times (\mathfrak{t}^{N-n})^*$ satisfies 
$$\mathcal{H}^{X(A,0)}_{i_k, t} \cong \{\eta \in (\mathfrak{t}^n)^*\: \vert \: \langle \mathbf{b}_{i_k}, \eta \rangle = 0\} 
\times \{t\} \:\:\:\: k = 1, ..., n$$ for all $t \in (\mathfrak{t}^{N-n})^*$, and,   
$$\mathcal{H}^{X(A,0)}_{i, t} \cong \{\eta \in (\mathfrak{t}^n)^*\: \vert \: \langle \mathbf{b}_i, \eta \rangle = f_i(t)\} \times \{t\}$$ 
with linear functions $f_i(t)$ on $(\mathfrak{t}^{N-n})^*$ for other $i$. As remarked in Example \ref{(4.5)}, these $N-n$ linear functions $\{f_i(t)\}$ are linearly independent.  
We can take a suitable line $\mathbf{C}^1 \to \mathbf{C}^{N-n}$ passing through $0$ and take the fiber product 
\begin{equation} 
\begin{CD} 
\mathcal{Z}' @>>>  X(A, 0)\\ 
@VVV @VVV  \\  
\mathbf{C}^1 @>>> \mathbf{C}^{N-n}
\end{CD} 
\end{equation}
so that the discriminant divisor $\mathcal{H}'$ of the relative moment map 
$$\mu_{\mathcal{Z}'}: \mathcal{Z}' \to (\mathfrak{t}^n)^* \times \mathbf{C}^1$$ 
satisfies $$\mathcal{H}' = \mathcal{H}.$$ 
Let $F_X \subset (\mathrm{t}^n)^*$ be the closed subset for $X$ defined in \ref{(3.4)} and let $F'_{Y(A,0)} \subset (\mathrm{t}^n)^*$ be the 
closed subset similarly defined for $Y(A, 0)$. We put $F := F_X \cup F_{Y(A,0)}$.  
In the remainder we regard $F$ as a closed subset 
of $(\mathrm{t}^n)^* \times \{0\}$ Since $\mathrm{Codim}_{(\mathrm{t}^n)^*}F \geq 2$,  
we have $$\mathrm{codim}_{(\mathfrak{t}^n)^* \times \mathbf{C}^1}F \geq 3.$$ 
Once the discriminant divisor $\mathcal{H}$ is fixed, both of the relative moment maps  $\mu_{\mathcal Z}$ and 
$\mu_{\mathcal{Z}'}$ have the same local form around each $(\eta, 0) \in (\mathfrak{t}^n)^* \times \mathbf{C}^1 - F$. 
More precisely, we get 

\begin{Prop}\label{(5.1)}
We have an open neighborhood 
$(\eta, t) \in U \subset (\mathfrak{t}^n)^* \times \mathbf{C}^1 - F$ such that 
there is a $T^n$-equivariant commutative diagram of the two families of symplectic varieties together with the relative moment maps: 
\begin{equation} 
\begin{CD} 
(\mu^{-1}_{\mathcal Z}(U), \omega_{{\mathcal Z}/\mathbf{C}^1}) @>{\Psi_U \:\: (\cong)}>> 
(\mu^{-1}_{{\mathcal Z}'}(U), \omega_{{\mathcal Z}'/\mathbf{C}^1})\\ 
@V{\mu_{\mathcal Z}}VV @V{\mu_{{\mathcal Z}'}}VV  \\  
U @>{id}>> U.
\end{CD} 
\end{equation} 
\end{Prop}

{\em Proof}. For $(\eta, 0) \in (\mathfrak{t}^n)^* \times \{0\} -F$, we have described the local form of the 
relative moment map in Proposition \ref{(4.6)}. 
For $(\eta, t) \in (\mathfrak{t}^n)^* \times \{t\}$ with $t \ne 0$, we put $U_t := U \cap ((\mathfrak{t}^n)^* \times \{t\})$. 
By the $\mathbf{C}^*$-action we see that  
$$(\mu^{-1}_{\mathcal Z}(U), \omega_{{\mathcal Z}/\mathbf{C}^1}) \to U, \:\:\mathrm{and}\:\:   
(\mu^{-1}_{{\mathcal Z}'}(U), \omega_{{\mathcal Z}'/\mathbf{C}^1}) \to U$$
are respectively trivial deformations of   
$$\mu_{\mathcal{Z}_t}: (\mu^{-1}_{{\mathcal Z}_t}(U_t), \omega_{{\mathcal Z}_t}) \to U_t, \:\:\mathrm{and}\:\:   
\mu_{\mathcal{Z}'_t}: (\mu^{-1}_{{\mathcal Z}'_t}(U_t), \omega_{{\mathcal Z}'_t}) \to U_t.$$
On the other hand, by Theorem \ref{(1.3)}, (2) we already know that
$\mu_{\mathcal{Z}_t}$ and $\mu_{\mathcal{Z}'_t}$ have the same local form.  $\square$ 
\vspace{0.2cm}

For simplicity, we put $S := (\mathfrak{t}^n)^* \times \mathbf{C}^1$ and $S^0 := (\mathfrak{t}^n)^* \times \mathbf{C}^1 - F$. 
We then define $$\mathcal{Z}^0 := \mu_{\mathcal Z}^{-1}(S^0), \:\:\: {\mathcal{Z}'}^0 := \mu_{{\mathcal Z}'}^{-1}(S^0).$$ 
Let ${\mathcal Aut}^{\mathcal{Z}/\mathbf{C}^1}$ be the sheaf on $(\mathfrak{t}^n)^* \times \mathbf{C}^1$ of Hamiltonian automorphisms 
of $(\mathcal{Z}, \omega_{\mathcal{Z}/\mathbf{C}^1})$. More precisely, for an open set $U \subset (\mathfrak{t}^n)^* \times \mathbf{C}^1$,  
the group ${\mathcal Aut}^{\mathcal{Z}/\mathbf{C}^1}(U)$ consists of the automorphisms $\tau$ of $\mu_{\mathcal Z}^{-1}(U)$ over $U$ such that $\tau$ preserve $\omega_{\mathcal{Z}/\mathbf{C}^1}\vert_{\mu_{\mathcal Z}^{-1}(U)}$ and are $T^n$-equivariant. 
 
We write ${\mathcal Aut}^{\mathcal{Z}^0/\mathbf{C}^1}$ for ${\mathcal Aut}^{\mathcal{Z}/\mathbf{C}^1}\vert_{S^0}.$ 
As in [\cite{Lo}, \S 3] we define a sheaf homomorphsm $\mathcal{O}_{S^0} \to {\mathcal Aut}^{\mathcal{Z}^0/\mathbf{C}^1}$. 
The following argument is almost the same as in \cite{Lo}. For $z \in \mathcal{Z}^0$, we denote by $T^n_z$ the stabilizer group of $z$ for
the $T^n$-action. If $\zeta \in \mathrm{Hom}_{alg. gp}(T, \mathbf{C}^*)\otimes \mathbf{R}$ 
is general enough, then $\zeta\vert_{\mathrm{Lie}(T^n_z)} \ne 0$ for all $z$. Then we define a Zariski open subset $(\mathcal{Z}^0)^{\zeta} \subset \mathcal{Z}^0$ as the set of points $z \in \mathcal{Z}$ such that \vspace{0.2cm}

(1) $T^n_z = \{1\}$, 

(2) for any $\lambda \in \mathrm{Hom}_{alg.gp}(\mathbf{C}^*, T^n)$ such that $\lim_{t \to 0} \lambda(t)\cdot z$ exists, 
we have $\langle \zeta, \lambda \rangle > 0$. \vspace{0.2cm}
  
Then, as in [\cite{Lo}, Lemma 2.6], $\mu_{\mathcal Z}\vert_{(\mathcal{Z}^0)^{\zeta}}: (\mathcal{Z}^0)^{\zeta} \to S^0$ is a principal $T^n$-bundle. 
  
Let $f$ be a holomorphic function on $U \subset S$. Regard $\mathcal{Z}$ as a Poisson variety over $S$ and denote by 
$\{\cdot, \cdot \}$ the $\mathcal{O}_S$-linear Poisson bracket. 
We define the Hamiltonian vector field $H_{\mu_{\mathcal Z}^*f}$ on $\mu_{\mathcal Z}^{-1}(U)$ by $\{\mu_{\mathcal Z}^*f, \cdot\}$. 
Then $H_{\mu_{\mathcal Z}^*f}$ preserves the Poisson structure. 
The $T^n$-action determines a vector field $\zeta_a$ on $\mathcal{Z}$ for  $a \in \mathfrak{t}^n$. By the first projection $(\mathfrak{t}^n)^* \times \mathbf{C}^1 \to (\mathfrak{t}^n)^*$, we regard $a$ as a linear function on $S$. 
By the definition of the relative moment map, 
we have $H_{\mu_{\mathcal Z}^*a} = \zeta_a$.  
Then we see that 
$$H_{\mu_{\mathcal Z}^*f}(d\mu_{\mathcal Z}^*a) = \{\mu_{\mathcal Z}^*f, \mu_{\mathcal Z}^*a\} = 
- \zeta_a (\mu_{\mathcal Z}^*(df)) = 0.$$ 
The last equality follows from the fact that every $T^n$-orbit is contained in a fiber of $\mu_{\mathcal Z}$. 
Therefore $H_{\mu_{\mathcal Z}^*f}$ is tangential to all fibers of $\mu_{\mathcal Z}$. 

Let us consider two fiber bundles $\mathfrak{t}^n \times S \to S$ 
and $T^n \times S \to S$ respectively with a typical fiber $\mathfrak{t}^n$ and a typical fiber $T^n$. 
Let $\mathfrak{t}^n_S$ and $T^n_S$ be the sheaves of holomorphic sections of these fiber bundles. 
There is a map of sheaves 
$$\mathrm{exp}(2\pi i): \mathfrak{t}^n_{S} \to T^n_{S}\:\:\:\:\: \xi \to \mathrm{exp}(2\pi i \xi ).$$ 
Moreover, $T^n_{S}(U)$ acts on $\mu_{\mathcal Z}^{-1}(U)$ by $\phi. z := \phi (\mu_{\mathcal Z}(z))\cdot z$, where 
$\phi (\mu_{\mathcal Z}(z)) \in T^n$ and $\cdot$ denotes the $T^n$-action on  
$\mu_{\mathcal Z}^{-1}(U)$.

Let $U \subset S^0$ be an open set of $S^0$. For a given $f \in \mathcal{O}_{S^0}(U)$, we construct an element 
of $T^n_{S^0}(U)$. In order to do this, we first consider 
$\mu_{\mathcal Z}\vert_{(\mathcal{Z}^0)^{\zeta}}^{-1}(U) = \mu_{\mathcal Z}^{-1}(U) \cap (\mathcal{Z}^0)^{\zeta}$. 
Note that each fiber of $\mu_{\mathcal Z}\vert_{(\mathcal{Z}^0)^{\zeta}}^{-1}(U)$ consists of a single free $T^n$-orbit. 
Then $H_{\mu_{\mathcal{Z}}^*f}$ determines a $T^n$-invariant vector field on each fiber. Therefore $H_{\mu_{\mathcal{Z}}^*f}$ 
is regarded as an element of $\mathfrak{t}^n_{S^0}(U)$. Then we have an element $\mathrm{exp}(2\pi i H_{\mu_{\mathcal Z}^*f} ) 
\in T^n_{S^0}(U)$ by the map $\mathrm{exp}(2\pi i): \mathfrak{t}^n_{S^0} \to T^n_{S^0}$. 
As $T^n_{S^0}(U)$ acts on $\mu_{\mathcal Z}^{-1}(U)$, $\mathrm{exp}(2\pi i H_{\mu_{\mathcal Z}^*f})$ determines an element 
of ${\mathcal Aut}^{\mathcal{Z}^0/\mathbf{C}^1}(U)$. As a consequence, we have a sheaf homomorphism 
$$\mathcal{O}_{S^0} \to {\mathcal Aut}^{\mathcal{Z}^0/\mathbf{C}^1}, \:\:\:\: f \to \mathrm{exp}(2\pi i H_{\mu_{\mathcal Z}^*f}).$$ 

Define $$X(T^n) := \mathrm{Hom}_{alg.gp}(T^n, \mathbf{C}^*).$$ Then an element of its dual $X(T^n)^*$ is a linear function 
on $(\mathfrak{t}^n)^*$. We regard $X(T^n)^*$ as a constant sheaf contained in $\mathcal{O}_{(\mathfrak{t}^n)^*}$. 
Let $p_1: S^0 \to (\mathfrak{t}^n)^*$ be the first projection. Then $p_1^{-1}(X(T^n)^*) \subset \mathcal{O}_{S^0}$ is also a constant sheaf on $S^0$, which we also denote by $X(T^n)^*$.  
By the almost same arguments as in [\cite{Lo}, Lemma 3.2, Lemma 3.3] we have 

\begin{Prop}\label{(5.2)} 
There is an exact sequence of abelian sheaves on $S^0$
$$ 0 \to p_2^{-1}\mathcal{O}_{\mathbf{C}^1} \oplus X(T^n)^* \to \mathcal{O}_{S^0} \to   
{\mathcal Aut}^{\mathcal{Z}^0/\mathbf{C}^1} \to 0,$$ where 
$p_2: S^0 \to \mathbf{C}^1$ is the second projection. 
\end{Prop}

\begin{Cor}\label{(5.3)} 
$$H^1(S^0, {\mathcal Aut}^{\mathcal{Z}^0/\mathbf{C}^1}) = 0.$$
\end{Cor}

{\em Proof}. We have an exact sequence 
$$H^1(S^0, \mathcal{O}_{S^0}) \to H^1(S^0, {\mathcal Aut}^{\mathcal{Z}^0/\mathbf{C}^1}) \to 
H^2(S^0, p_2^{-1}\mathcal{O}_{\mathbf{C}^1} \oplus X(T^n)^*) \to H^2(S^0, \mathcal{O}_{S^0}).$$ 
Note that $S^0 = S - F$ with $\mathrm{Codim}_SF \geq 3$. We see that 
$H^2_F(S, \mathcal{O}_S) = 0$ by the depth argument. Then $H^1(S, \mathcal{O}_S) \to H^1(S^0, \mathcal{O}_{S^0})$ 
is a surjection. Since $S$ is Stein, we have $H^1(S, \mathcal{O}_S) = 0$; hence, $H^1(S^0, \mathcal{O}_{S^0}) = 0$.   
Since $X(T)^* \cong \mathbf{Z}^{\oplus n}$ as a constant sheaf, the map $H^2(S, X(T)^*) \to H^2(S^0, X(T)^*)$ is an isomorphism 
because $\mathrm{Codim}_SF \geq 2$. On the other hand, $H^2(S, X(T)^*) = 0$; hence $H^2(S^0, X(T)^*) = 0$. 
Now we prove that $$\mathrm{Ker}[H^2(S^0, p_2^{-1}\mathcal{O}_{\mathbf{C}^1}) \to H^2(S^0, \mathcal{O}_{S^0})] = 0.$$ 
Since (analytic) de Rham complex 
$$0 \to p_2^{-1}\mathcal{O}_{\mathbf{C}^1}  \to \mathcal{O}_{S^0} \to \Omega^1_{S^0/\mathbf{C}^1} \to \Omega^2_{S^0/\mathbf{C}^1} \to \cdot\cdot\cdot $$ 
is exact, we have $$H^2(S^0, p_2^{-1}\mathcal{O}_{\mathbf{C}^1}) = \mathbf{H}^2(S^0, \Omega^{\cdot}_{S^0/\mathbf{C}^1}).$$
Consider the Hodge to de Rham spectral sequence 
$$E^{p, q}_1(S^0) := H^q(S^0, \Omega^p_{S^0/\mathbf{C}^1}) \Rightarrow \mathbf{H}^2(S^0, \Omega^{\cdot}_{S^0/\mathbf{C}^1}), \:\:\: 
p + q = 2.$$ 
The spectral sequence determines a decreasing filtration $F^{\cdot}$ on $\mathbf{H}^2(S^0, \Omega^{\cdot}_{S^0/\mathbf{C}^1})$ and 
we have $$\mathrm{Ker}[H^2(S^0, p_2^{-1}\mathcal{O}_{\mathbf{C}^1}) \to H^2(S^0, \mathcal{O}_{S^0})]  = F^1(\mathbf{H}^2(S^0, \Omega^{\cdot}_{S^0/\mathbf{C}^1})). $$ 
Similarly, we have the Hodge to de Rham spectral sequence for the de Rham complex $\Omega^{\cdot}_{S/\mathbf{C}^1}$ on $S$ 
and we get a decreasing filtraton $F^{\cdot}$ on $\mathbf{H}^2(S, \Omega^{\cdot}_{S/\mathbf{C}^1})$. 
We shall prove that $$F^1(\mathbf{H}^2(S, \Omega^{\cdot}_{S/\mathbf{C}^1})) = 
F^1(\mathbf{H}^2(S^0, \Omega^{\cdot}_{S^0/\mathbf{C}^1})).$$
Let us compute $E_{\infty}^{1,1}(S^0)$. 
By the complex 
$$E_1^{0, 1}(S^0) \stackrel{d_{0,1}}\to E_1^{1,1}(S^0) \stackrel{d_{1,1}} \to E_1^{2, 1}(S^0), $$
we have $E_2^{1,1}(S^0) := \mathrm{Ker}(d_{0,1})/\mathrm{Im}(d_{1,1})$. 
On the other hand, the complex 
$$E_1^{2,0}(S^0) \stackrel{d_{2,0}}\to E_1^{3,0}(S^0) \stackrel{d_{3,0}}\to E_1^{4,0}(S^0)$$ 
yields $E_2^{3,0}(S^0) = \mathrm{Ker}(d_{3,0})/\mathrm{Im}(d_{2,0})$. 
Finally $$E_{\infty}^{1,1}(S^0) = E_3^{1,1}(S^0) = \mathrm{Ker}[E_2^{1,1}(S^0) \to E_2^{3,0}(S^0)].$$
Since $\mathrm{Codim}_SF \geq 3$, we have $E_1^{i,1}(S^0) = E_1^{i,1}(S)$ and 
$E_1^{i,0}(S^0) = E_1^{i,0}(S)$ for all $i$. Hence $E_2^{1,1}(S^0) = E_2^{1,1}(S)$ and $E_2^{3,0}(S^0) = E_2^{3,0}(S)$. 
This implies that $E_{\infty}^{1,1}(S^0) = E_{\infty}^{1,1}(S)$. 

We next compute $E_{\infty}^{2,0}(S^0)$. 
By the complex 
$$E_1^{1,0}(S^0) \stackrel{d_{1,0}}\to E_1^{2,0}(S^0) \stackrel{d_{2,0}} \to E_1^{3, 0}(S^0), $$
we have $E_2^{2,0}(S^0) := \mathrm{Ker}(d_{2,0})/\mathrm{Im}(d_{1,0})$. 
On the other hand, 
$E_2^{0,1}(S^0) = \mathrm{Ker}[E_1^{0,1}(S^0) \to E_1^{1,1}(S^0)]$. 
Finally $$E_{\infty}^{2,0}(S^0) = E_3^{2,0}(S^0) = \mathrm{Coker}[E_2^{0,1}(S^0) \to E_2^{2,0}(S^0)].$$
Since $\mathrm{Codim}_SF \geq 3$, we have $E_1^{i,1}(S^0) = E_1^{i,1}(S)$ and 
$E_1^{i,0}(S^0) = E_1^{i,0}(S)$ for all $i$. Hence $E_2^{0,1}(S^0) = E_2^{0,1}(S)$ and $E_2^{2,0}(S^0) = E_2^{2,0}(S)$. 
This implies that  $E_{\infty}^{2,0}(S^0) = E_{\infty}^{2,0}(S)$.

These show that  
$$F^1(\mathbf{H}^2(S, \Omega^{\cdot}_{S/\mathbf{C}^1})) = 
F^1(\mathbf{H}^2(S^0, \Omega^{\cdot}_{S^0/\mathbf{C}^1})).$$ 

We next show that $H^2(S, p_2^{-1}\mathcal{O}_{\mathbf{C}^1}) = 0$.
In order to do, we apply the Leray spectral sequence 
$$E_2^{p,q} := H^p(S, R^q(p_2)_*p_2^{-1}\mathcal{O}_{\mathbf{C}^1}) \Rightarrow H^2(S, p_2^{-1}\mathcal{O}_{\mathbf{C}^1}).$$
Since $S = (\mathfrak{t}^n)^* \times \mathbf{C}^1$, we have an isomorphism 
$$R^p(p_2)_*p_2^{-1}\mathcal{O}_{\mathbf{C}^1} \cong R^p(p_2)_*\mathbf{C} \otimes_{\mathbf C}\mathcal{O}_{\mathbf{C}^1}.$$ 
Note that $R^p(p_2)_*\mathbf{C} = 0$ for $p > 0$ and $(p_2)_*\mathbf{C} = \mathbf{C}$. Since 
$H^2(\mathbf{C}^1, \mathcal{O}_{\mathbf{C}^1}) = 0$, we see that $H^2(S, p_2^{-1}\mathcal{O}_{\mathbf{C}^1}) = 0$.

In particular, we have $F^1(\mathbf{H}^2(S, \Omega^{\cdot}_{S/\mathbf{C}^1})) = 0$; hence, 
$$\mathrm{Ker}[H^2(S^0, p_2^{-1}\mathcal{O}_{\mathbf{C}^1}) \to H^2(S^0, \mathcal{O}_{S^0})] = 0.$$
By the first exact sequence in the proof, we see that 
$$H^1(S^0, {\mathcal Aut}^{\mathcal{Z}^0/\mathbf{C}^1}) = 0.$$ $\square$

\begin{Cor}\label{(5.4)} 
There is a $T^n$-equivariant isomorphism $$\Psi^0: (\mathcal{Z}^0, \omega_{\mathcal{Z}^0/\mathbf{C}^1}) 
\cong ((\mathcal{Z}')^0, \omega_{(\mathcal{Z}')^0/\mathbf{C}^1})$$ which makes the following diagram commutative 
\begin{equation} 
\begin{CD} 
(\mathcal{Z}^0, \omega_{\mathcal{Z}^0/\mathbf{C}^1}) @>{\Psi^0}>> ((\mathcal{Z}')^0, \omega_{(\mathcal{Z}')^0/\mathbf{C}^1})\\  
@V{\mu_{{\mathcal Z}}}VV @V{\mu_{{\mathcal Z}'}}VV  \\  
S^0 @>{id}>> S^0.
\end{CD} 
\end{equation}
\end{Cor}

{\em Proof}. By Proposition \ref{(5.1)} we cover $S^0$ by open sets $U_i$ $(i \in I)$ such that there are $T^n$-equivariant 
isomorphisms $\Psi_i: \mu^{-1}_{\mathcal Z}(U_i) \to 
\mu^{-1}_{{\mathcal Z}'}(U_i)$. Then $$\Psi_{ij} := (\Psi_i)^{-1} \circ \Psi_j \:\: \vert_{\mu^{-1}_{\mathcal Z}(U_i \cap U_j)}: 
\mu^{-1}_{\mathcal Z}(U_i \cap U_j) \to \mu^{-1}_{\mathcal Z}(U_i \cap U_j)$$ is a 1-cocycle in 
${\mathcal Aut}^{\mathcal{Z}^0/\mathbf{C}^1}$. By Corollary \ref{(5.3)}, if we choose each $U_i$ small enough, then there 
are Hamiltonian automorphisms $f_i$ of $\mu^{-1}_{\mathcal Z}(U_i)$ such that   
$\Psi_{ij} = f_i \circ f_j^{-1}\vert_{\mu^{-1}_{\mathcal Z}(U_i \cap U_j)}$. Replace $\Psi_i$ by $\Psi_i \circ f_i$ for each $i$. Then 
$\{\Psi_i\}$ glue together to give an isomorphism $\Psi^0$. $\square$ 

\begin{Thm}\label{(5.5)} 
There is a $T^n$-equivariant isomorphism $$\Psi: (\mathcal{Z}, \omega_{\mathcal{Z}/\mathbf{C}^1}) 
\cong (\mathcal{Z}', \omega_{\mathcal{Z}'/\mathbf{C}^1})$$ which makes the following diagram commutative 
\begin{equation} 
\begin{CD} 
(\mathcal{Z}, \omega_{\mathcal{Z}/\mathbf{C}^1}) @>{\Psi}>> (\mathcal{Z}', \omega_{\mathcal{Z}'/\mathbf{C}^1})\\  
@V{\mu_{{\mathcal Z}}}VV @V{\mu_{{\mathcal Z}'}}VV  \\  
S @>{id}>> S.
\end{CD} 
\end{equation}   
\end{Thm}

{\em Proof}. By Proposition \ref{(2.1)}, (2) each fiber of $\mu_{\mathcal Z}: \mathcal{Z} \to S$ has dimension $n$. 
Then we have $$\mathrm{Codim}_{\mathcal Z}(\mathcal{Z} - \mathcal{Z}^0) \geq  3$$ because $\mathrm{Codim}_S (S - S^0) \geq 3$. 
Similarly we have $$\mathrm{Codim}_{{\mathcal Z}'}(\mathcal{Z}' - (\mathcal{Z}')^0) \geq 3.$$
Since $\mathcal{Z}$ and $\mathcal{Z}'$ are normal, we get $$\Gamma (\mathcal{Z}, \mathcal{O}_{\mathcal{Z}}) = \Gamma (\mathcal{Z}^0, \mathcal{O}_{\mathcal{Z}^0}), \:\:\:\:  
\Gamma (\mathcal{Z}', \mathcal{O}_{\mathcal{Z}'}) = \Gamma ((\mathcal{Z}')^0, \mathcal{O}_{(\mathcal{Z}')^0}).$$
Hence we see that $$\Gamma (\mathcal{Z}, \mathcal{O}_{\mathcal{Z}}) = \Gamma (\mathcal{Z}', \mathcal{O}_{\mathcal{Z}'})$$ 
by Corollary \ref{(5.4)}. Since $\mathcal{Z}$ and $\mathcal{Z}'$ are both Stein spaces, it follows that $\mathcal{Z} \cong \mathcal{Z}'$. 
By the construction of the isomorphism, this is a $T^n$-equivariant and preserves the relative moment maps $\mu_{\mathcal Z}$ and 
$\mu_{{\mathcal Z}'}$. $\square$ 
\vspace{0.2cm}

If we restrict the commutative diagram in Theorem \ref{(5.5)} over $(\mathfrak{t}^n)^* \times \{0\} \subset S$, then we get: 
 
 \begin{Cor}\label{(5.6)} 
Assume that $X$ satisfies the condition (*) at the beginning of \S 5. Then there is a $T^n$-equivariant isomorphism $\varphi: (X, \omega) \to (Y(A, 0), \omega_{Y(A,0)})$ 
which makes the following diagram commutative  
\begin{equation} 
\begin{CD} 
(X, \omega) @>{\varphi}>> (Y(A, 0), \omega_{Y(A,0)}) \\  
@V{\mu}VV @V{\bar{\mu}}VV  \\  
(\mathfrak{t}^n)^* @>{id}>> (\mathfrak{t}^n)^*.
\end{CD} 
\end{equation}   
Moreover,  $\varphi (0_X) = 0_{Y(A, 0)}$. 
\end{Cor}

{\em Proof}. We need a proof for the last statement. Let us consider a singular point $p$ of $X$ fixed by all elements of $T^n$. 
Note that the origin $0_X \in X$ satisfies the property. We prove that $p = 0_X$.
Take the symplectic leaf $Y$ of $X$ containing $p$. Then the algebraic torus $T_Y$ of $1/2 \cdot \dim Y$ acts on $(Y, \omega_Y)$ by Theorem \ref{(2.2)}. By applying Theorem \ref{(1.3)} and Example \ref{(1.2)} to a $T_Y$-stable affine open neighborhood of $p \in Y$,  
we see that such points fixed by $T_Y$ are all isolated in $Y$. Assume that $p \ne 0_X$. 
For the conical $\mathbf{C}^*$-action on $X$, we consider the points 
$\sigma (p)$ ($\sigma \in \mathbf{C}^*$). They are all contained in $Y$. Moreover, they are fixed by $T_Y$ because the conical $\mathbf{C}^*$-action commutes with the $T^n$-action. This contradicts that such points are all isolated in $Y$.  
Therefore $p = 0_X$. Similarly, $0_{Y(A,0)}$ is the 
unique point of $Y(A,0)$ fixed by $T^n$. This means that $\varphi (0_X) = 0_{Y(A, 0)}$. $\square$

\begin{Rem}\label{(5.7)} 
{\rm $\varphi$ is a morphism of complex analytic varieties, not necessarily a morphism of algebraic varieties.}  
\end{Rem}

In the remainder we discuss what happens when the condition (*) does not hold. 

(1) The case when $n = N$ and $B$ is an injection: 

Since $B$ is unimodular,  $B: \mathbf{Z}^n \to \mathbf{Z}^n$ is an isomorphism.  
Then $Y(A, 0) = \mathbf{C}^{2n}$ and $\omega_{Y(A,0)}$ is the standard symplectic form $\omega_{st}$.   
In this case we do not have $\mathcal{Z}$. Instead we consider the sheaf ${\mathcal Aut}^X$ on $(\mathfrak{t}^n)^*$ 
(cf. \cite{Lo}). Then we have $H^1((\mathfrak{t}^n)^*, {\mathcal Aut}^X) = 0$ and we get a $T^n$-equivariant isomorphism 
$$\varphi: (X, \omega) \cong (\mathbf{C}^{2n}, \omega_{st})$$ and the moment map $\mu$ is given by 
$$\mathbf{C}^{2n} \to \mathbf{C}^n, \:\:\:\:\: (z_1, ..., z_n, w_1, ..., w_n) \to (z_1w_1, ..., z_nw_n).$$ 

(2) The case when $B$ is not an injection: 

We can write $\mathbf{Z}^n = \mathrm{Ker}(B) \oplus \mathbf{Z}\langle \mathbf{b}_1, ..., \mathbf{b}_N \rangle$.
We put $\mathbf{Z}^{n'} := \mathbf{Z}\langle \mathbf{b}_1, ..., \mathbf{b}_N \rangle$ and consider the injection 
$\mathbf{Z}^{n'} \to \mathbf{Z}^N$, which we denote again by $B$. Since $B$ is unimodular, we have an exact sequence 
$$0 \to \mathbf{Z}^{n'} \stackrel{B}\to \mathbf{Z}^N \stackrel{A}\to \mathbf{Z}^{N-n'} \to 0.$$
For the toric hyperk\"{a}hler variety $Y(A, 0)$ we construct its Poisson deformation $\mathcal{Z}' \to \mathbf{C}^1$. 
Now we consider the Poisson deformation $\mathcal{Z}' \times T^*{(\mathbf{C}^*)^{n-n'}} \to \mathbf{C}^1$ of 
$Y(A, 0) \times T^*{(\mathbf{C}^*)^{n-n'}}$. 
Let  
$$\mu_{\mathcal{Z}' \times T^*{(\mathbf{C}^*)^{n-n'}}}: \mathcal{Z}' \times T^*{(\mathbf{C}^*)^{n-n'}} \to (\mathfrak{t}^{n'})^* \oplus (\mathfrak{t}^{n-n'})^* \times \mathbf{C}^1$$ be the relative moment map. 
Then we see that 
$$\mathcal{H}_{\mathcal Z} = \mathcal{H}_{\mathcal{Z}' \times T^*{(\mathbf{C}^*)^{n-n'}}}.$$   
As in Theorem \ref{(5.5)} we have a $T^n$-equivariant isomorphism $\Psi: \mathcal{Z} \cong 
 \mathcal{Z}' \times T^*{(\mathbf{C}^*)^{n-n'}}$ which makes the following diagram commutative  
\begin{equation} 
\begin{CD} 
\mathcal{Z} @>{\Psi}>> \mathcal{Z}' \times T^*{(\mathbf{C}^*)^{n-n'}} \\  
@V{\mu_{{\mathcal Z}}}VV @V{\mu_{\mathcal{Z}' \times T^*{(\mathbf{C}^*)^{n-n'}}}}VV  \\  
(\mathfrak{t}^n)^* \times \mathbf{C}^1@>{id}>> (\mathfrak{t}^n)^* \times \mathbf{C}^1.
\end{CD} 
\end{equation}
Restricting this diagram above $(\mathfrak{t}^n)^* \times \{0\}$, we get an isomorphism 
$$X \cong Y(A, 0) \times T^*{({\mathbf C}^*)^{n-n'}}.$$
Since $X$ is a conical symplectic variety, it must be contractible as a topological space. However, the right hand side 
is not contractible. This is a contradiction. Hence the case (2) does not occur. 

As a conclusion we have 

\begin{Thm}\label{(5.8)} 
Let $(X, \omega)$ be a conical symplectic variety of dimension $2n$ which has a projective symplectic resolution. 
Assume that $X$ admits an effective Hamiltonian action of an $n$-dimensional algebraic torus $T^n$, compatible with the conical 
$\mathbf{C}^*$-action. Then there is a $T^n$-equivariant (complex analytic) isomorphism $\varphi: (X, \omega) \to (Y(A, 0), \omega_{Y(A,0)})$ 
which makes the following diagram commutative  
\begin{equation} 
\begin{CD} 
(X, \omega) @>{\varphi}>> (Y(A, 0), \omega_{Y(A,0)}) \\  
@V{\mu}VV @V{\bar{\mu}}VV  \\  
(\mathfrak{t}^n)^* @>{id}>> (\mathfrak{t}^n)^*.
\end{CD} 
\end{equation}   
Moreover,  $\varphi (0_X) = 0_{Y(A, 0)}$.  
\end{Thm}

\begin{Rem}\label{(5.9)} {\rm 
As the following example shows, the map $\varphi$ cannot be chosen to be $\mathbf{C}^*$-equivariant in general.  
For $n > 1$, we define $X := \{(x,y,z) \in \mathbf{C}^3\: \vert \: xy = z^n\}$ and 
put $\omega := \mathrm{Res}(\frac{dx \wedge dy \wedge dz}{xy - z^n})$. Then $\mathbf{C}^*$ acts on $X$ by 
$x \to \sigma^a x \:\: y \to \sigma^{2n-a}y, \:\: z \to \sigma^2z, \:\: \sigma \in \mathbf{C}^*$ for a positive integer $a$ with $0 < a < 2n$.  With this $\mathbf{C}^*$-action, $(X, \omega)$ becomes a conical symplectic variety with $wt(\omega) = 2$. Moreover, another 1-dimensional algebraic torus $T$ acts on $(X, \omega)$ by $x \to ty, \: y \to t^{-1}y, \: z \to z$ for $t \in T$. 

Let us consider the exact sequence 
$$0 \to \mathbf{Z} \stackrel{B}\to \mathbf{Z}^n \stackrel{A}\to \mathbf{Z}^{n-1} \to 0$$ 
with 
$A = \left(\begin{array}{cccccccc} 
 1 & 0 & ... & 0 & -1 \\
 0 & 1 & ... & 0 & -1\\ 
 ... & ... & ... & ...&  ... \\ 
 0 & ... & ... & 1 & -1 
\end{array}\right)$   
and 
$B = {}^t(1, 1, ..., 1)$. Let $(Y(A,0), \omega_{Y(A,0)})$ be the corresponding toric hyperk\"{a}hler variety. 
As explained above, $Y(A, 0)$ is a conical symplectic variety with a standard conical $\mathbf{C}^*$-action. 
Then we have a $T$-equivariant isomorphism $\varphi: (X, \omega) \cong (Y(A, 0), \omega_{Y(A, 0)})$. When $a = n$, this is also 
a $\mathbf{C}^*$-equivariant map. But, when $a \ne n$, there is no $\mathbf{C}^*$-equivariant isomorphism between $X$ and $Y(A, 0)$. 
This example suggests that there would be a {\em canonical one} among conical $\mathbf{C}^*$-actions (cf. \cite{Na-Od}). $\square$ } 
\end{Rem}

Assume that $wt(\omega) = 2$ for the conical $\mathbf{C}^*$-action on $X$. Note that $\mathbf{C}^* \times T^n$ acts on $X$. 
Take a homomorphism $\iota: \mathbf{C}^* \to \mathbf{C}^* \times T^n$ so that the composite $\mathbf{C}^* \stackrel{\iota}\to \mathbf{C}^* \times T^n \stackrel{p_1}\to \mathbf{C}^*$ is an isomorphism. Then we have $wt(\omega) = 2$ for the new $\mathbf{C}^*$-action.    
Now we can pose the following question: 
\vspace{0.2cm}

{\bf Question 5.10.} {\em Assume that $wt(\omega) = 2$. If necessary, after taking a suitable $\iota: \mathbf{C}^* \to \mathbf{C}^* \times T^n$ and replacing the original $\mathbf{C}^*$-action by the new 
conical $\mathbf{C}^*$-action,  can we take $\varphi: (X, \omega) \to (Y(A, 0), \omega_{Y(A, 0)})$ in a $\mathbf{C}^*$-equivariant way ?}



 



\begin{center}
Research Institute for Mathematical Sciences, Kyoto University, Oiwake-cho, Kyoto, Japan

E-mail address: namikawa@kurims.kyoto-u.ac.jp  
\end{center}


\begin{thebibliography}{99}

\bibitem[A-P]{A-P}  Arbo, M., Proudfoot, N.: Hypertoric varieties and zonotopal tilings, Int. Math. Res. Notices, {\bf 2016}, No.23 (2016), 7268 - 7301 


\bibitem[Be]{Be} Beauville, A.: Symplectic singularities, Invent. Math. {\bf 139} (2000), no.3, 541-549



 
\bibitem[B]{B} Bielawski, R.: Complete hyperk\"{a}hler 4n-manifolds with n commuting tri-Hamiltonian vector
fields, Math. Ann. {\bf 314} (1999), 505 - 528 
 
 
\bibitem[BD]{BD} Bielawski, R.,  Dancer, A.: 
The geometry and topology of toric hyperk\"{a}hler manifolds.
Comm. Anal. Geom. {\bf 8} (2000), no. 4, 727 - 760


 



\bibitem[Del]{Del} Delzant, T.: Hamiltoniens p\'{e}riodiques et images convexes de l'application moment. 
Bull. Soc. Math. France {\bf 116} (1988), no. 3, 315 - 339.


\bibitem[Go]{Go} Goto, R.: On toric hyper-K\"{a}hler manifolds given by the hyper-K\"{a}hler quotient method, in Infinite Analysis, World Scientific, 1992, pp. 317 - 338



\bibitem[G-S]{G-S} Guillemin, V., Sternberg, S.: Symplectic techniques in physics. Cambridge University Press, Cambridge, 1984. xi +468 pp.


\bibitem[HS]{HS} Hausel, T.; Sturmfels, B.: 
Toric hyperk\"{a}hler varieties. 
Doc. Math. {\bf 7} (2002), 495 - 534.

     

\bibitem[Ka]{Ka} Kaledin, D.: Symplectic singularities from the Poisson point of view.
J. Reine Angew. Math. {\bf 600} (2006), 135 - 156.



\bibitem[Ko]{Ko} Konno, H.: The geometry of toric hyperk\"{a}hler varieties. Toric topology, 241- 260,
Contemp. Math., {\bf 460}, Amer. Math. Soc., Providence, RI, 2008.


\bibitem[Lo]{Lo} Losev, I.:  Classification of multiplicity free Hamiltonian actions of algebraic tori on Stein manifolds. (English summary)
J. Symplectic Geom. {\bf 7} (2009), no. 3, 295 - 310.

\bibitem[Na 1]{Na 1} Namikawa, Y.: Flops and Poisson deformations of symplectic varieties. Publ. Res. Inst. Math. Sci. {\bf 44} (2008), no. 2, 259 - 314.



\bibitem[Na 2]{Na 2} Namikawa, Y.: Poisson deformations of affine symplectic varieties. Duke Math. J. {\bf 156} (2011), no. 1, 51 - 85.


\bibitem[Na 3]{Na 3} Namikawa, Y.: Poisson deformations and birational geometry.
J. Math. Sci. Univ. Tokyo {\bf 22} (2015), no. 1, 339 - 359.

 
 
\bibitem[Na-Od]{Na-Od} Namikawa, Y., Odaka, Y.: Canonical torus action on symplectic singularities, arXiv:2503.15791

 

\bibitem[Pr]{Pr} Proudfoot, N.: A survey of hypertoric geometry and topology. (English summary) Toric topology, 323 - 338, Contemp. Math. {\bf 460}, Amer. Math. Soc., Providence, RI, 2008


\bibitem[Pu]{Pu} Puerta, F.: 
D\'{e}formations semiuniverselles et germes d'espaces analytiques $\mathbf{C}^*$-\'{e}quivariantes. (French)
Algebraic geometry (La R\'{a}bida, 1981), 267 - 274,
Lecture Notes in Math., 961, Springer, Berlin, 1982.



\bibitem[S]{S} Schwarz, G.: Quotients of Compact and Complex Reductive Groups, in Th\'{e}orie des invariants \& G\'{e}ometrie des vari\'{e}t\'{e}s quotients, Travaux en Cours {\bf 61} (2000), Hermann et cie., Paris, 5 - 83.



\bibitem[Sn]{Sn} Snow, D.: Reductive group actions on Stein spaces.
Math. Ann. {\bf 259} (1982), no. 1, 79 - 97.

\end{thebibliography}
\end{document}